\documentclass[a4paper,reqno]{amsart}
\usepackage{amsmath}
\usepackage{amsfonts}
\usepackage{amsthm}
\usepackage{amssymb}
\usepackage{mathrsfs}
\usepackage{bm}
\usepackage{booktabs}
\usepackage{lipsum}
\usepackage{bussproofs}
\usepackage{hyperref}
\usepackage{tikz}
\usetikzlibrary{arrows,automata,positioning,decorations.markings}
\usepackage{python}

\usepackage{amsmath}
\usepackage{amsfonts}
\usepackage{amsthm}
\usepackage{amssymb}
\usepackage{mathrsfs}

\def\parentheses#1{{\left( {#1} \right)}}
\def\of{\parentheses}

\def\Fraction#1/#2{\frac{#1}{#2}}
\def\fraction#1/#2{\tfrac{#1}{#2}}
\def\solidus#1/#2{%
    {{#1} \!\left/\!\vphantom{{#1}{#2}}\!\right. {#2}}%
}

\def\operator{\operatorname}

\def\script{\mathscr}
\def\gothic{\mathfrak}
\def\blackboard{\mathbb}

\def\infinity{\infty}
\def\continuum{\gothic c}

\def\naturals{\blackboard N}

\def\set#1{{\left\lbrace {#1} \right\rbrace}}
\def\setOf#1:#2{%
    \set{{#1} \,\left\vert\vphantom{{#1}{#2}}\right. {#2}}%
}

\def\emptySet{\varnothing}

\def\sequence#1{{\left\langle {#1} \right\rangle}}
\def\sequenceOf#1:#2{%
    \sequence{{#1} \,\left\vert\vphantom{{#1}{#2}}\right. {#2}}%
}

\def\familyOf#1:#2{{\parentheses{#1}_{#2}}}

\def\superset{\supset}
\def\union{\cup}
\def\Union{\bigcup}
\def\intersect{\cap}
\def\Intersection{\bigcap}
\def\setMinus{\setminus}

\def\functionSpace#1->#2{{{#2}^{#1}}}
\def\cardinality#1{{\left\lvert {#1} \right\rvert}}

\def\function#1:#2->#3{{{#1} \colon {#2} \to {#3}}}
\def\mapsTo{\mapsto}
\def\map#1->#2{\parentheses{{#1} \mapsTo {#2}}}

\def\Restriction#1|#2{{\left.\! {#1} \right\vert_{#2}}} % \restriction is \upharpoonright.

\def\presentation#1{{\left\langle {#1} \right\rangle}}
\def\presentationOf#1:#2{%
    \presentation{{#1} \,\left\vert\vphantom{{#1}{#2}}\right. {#2}}%
}

\def\closure#1{\overline{#1}}

\def\interior#1{{\operator{int}\of{#1}}}

\def\collection{\script}

\def\interval#1#2,#3#4{{\left#1 {#2}, {#3} \right#4}}

\def\legal#1{{\operator{legal}\of{#1}}}
\def\illegal#1{{\operator{illegal}\of{#1}}}
\def\nondegenerate#1{{\operator{nondegenerate}\of{#1}}}
\def\degenerate#1{{\operator{degenerate}\of{#1}}}
\def\mahavier#1#2{{\bigstar^{#1}_{i=0} {#2}}}
\def\Graph#1{{\Gamma\of{#1}}}
\def\height#1{{\operator{height}\of{#1}}}
\def\heightIn#1#2{{\operator{height}_{#1}\of{#2}}}
\def\levelIn#1#2{{\operator{level}_{#2}\of{#1}}}
\def\transitiveIn#1#2{{\operator{trans}_{#2}\of{#1}}}
\def\intrans#1{{\operator{intrans}\of{#1}}}

\def\isolated#1{{\operator{isolated}\of{#1}}}

\def\p{\parentheses}
\newtheorem{presentment}{}[section]

\theoremstyle{plain}
\newtheorem{theorem}[presentment]{Theorem}
\newtheorem{proposition}[presentment]{Proposition}
\newtheorem{lemma}[presentment]{Lemma}
\newtheorem{corollary}[presentment]{Corollary}

\theoremstyle{definition}
\newtheorem{definition}[presentment]{Definition}
\newtheorem{example}[presentment]{Example}

\newtheorem{question}[presentment]{Question}
\newtheorem{observation}[presentment]{Observation}

\theoremstyle{remark}

\newtheorem*{note}{Note}

    \def\fraction#1/#2{\tfrac{#1}{#2}}
    \def\Fraction#1/#2{\frac{#1}{#2}}
\newcommand*{\BigFraction}{}
    \def\BigFraction#1/#2{\dfrac{\displaystyle{#1}}{\displaystyle{#2}}}

    \def\setOf#1:#2{\set{{#1} \,\left\vert\vphantom{{#1}{#2}}\right. {#2}}}

    \def\sequenceOf#1:#2{\sequence{{#1} \,\left\vert\vphantom{{#1}{#2}}\right. {#2}}}

    \def\function#1:#2->#3{{{#1} \colon {#2} \to {#3}}}

\usepackage{pst-node}
\usepackage{tikz-cd} 
\tikzset{node distance=2cm, auto}
\allowdisplaybreaks[0]
\title{Transitivity in CR-Dynamical Systems} 
\author{Sina Greenwood}


\email {s.greenwood@auckland.ac.nz}

\address
{Department of Mathematics, University of Auckland, 38 Princes Street, 1010
Auckland, New Zealand}

\author
{Andrew Wood}


\email {} \email {andrew.wood@stcatz.ox.ac.uk}

\address
{Mathematical Institute, University of Oxford, Andrew Wiles Building
Radcliffe Observatory Quarter 
Woodstock Road
Oxford
OX2 6GG }

\begin{document}

\maketitle
%---------------------------------------------------------------
\begin{abstract}
A CR-dynamical system is a pair $\p{X, G}$, where $X$ is a compact metric space and $G$ is a closed relation (CR) on $X$. In this paper, we introduce a new type of transitive point and transitivity in CR-dynamical systems. We develop a new tool called transitivity trees, which we use to determine the relationship between the different types of transitive points.  
\end{abstract}



%\footnote{Parts of this work appears in the MSc thesis of the author.}

\makeatletter
\renewcommand\@makefnmark{}
\makeatother
\footnotetext{Sina Greenwood is supported by the Marsden Fund Council from Government funding, administered by the Royal Society of New Zealand. This work appears in the MSc thesis of the second author.}

%---------------------------------------------------------------
\section{Introduction}
Ethan Akin's book ``The General Topology of Dynamical Systems''~\cite{akin} introduces the study of dynamics on 
compact metric spaces with closed relations.  Recently, this has led to the study of CR-dynamical systems~\cite{nagar_transitivity_CR, minimality_CR, transitivity_CR}, that is, compact metric spaces with closed relations.   Iztok Banic, Goran Erceg, Rene Gril Rogina, and Judy Kennedy's paper ``Minimal dynamical systems with closed relations''~\cite{minimality_CR} formally introduces CR-dynamical systems, generalising minimality from topological dynamical systems to CR-dynamical systems.  In this paper, we study transitive points, transitivity, and mixing, which have been generalised to CR-dynamical systems~\cite{nagar_transitivity_CR, transitivity_CR}. Iztok Banic et al. in~\cite{transitivity_CR} introduce three different types of transitive points for CR-dynamical systems. We introduce a new fourth type of transitive point. We also introduce two new types of transitivity, and one new type of mixing for CR-dynamical systems.

Set-valued dynamics is closely related to the study of CR-dynamical systems, since CR-dynamical systems generalise set-valued dynamical systems, and both generalise topological dynamical systems. Applications include the Christiano-Harrison model in macroeconomics~\cite{economics}. Raines and Stockman explored Devaney chaos, Li–Yorke chaos and distributional chaos in the Christiano-Harrison model~\cite{RAINES20101173}.

% In our paper, we focus on compact Hausdorff spaces with closed relations, generalising work by BaniÄ et al. in~\cite{transitivity_CR}.  In Section~\ref{section:Legally-Transitive-Points}, we find the answer to Problem $4.13$ in~\cite{minimality_CR} changes whether we allow such a generalisation. However, for generalising certain results, we require our spaces to be compact metric. We also focus on compact metric spaces with respect to examples. 

We provide the background necessary for CR-dynamical systems in Section~\ref{section:CR-dynamical-systems}. The orbit structure of a given point $x$ in a topological dynamical system $\p{X, f}$ can be thought of as its trajectory $\sequence{x, f\of{x}, f^2\of{x}, \ldots}$. In CR-dynamical systems $\p{X, G}$, points may have multiple trajectories, or even no trajectories at all. As the structure of these orbits can be complicated, and there are many different types of structures one may obtain, it is desirable to have a tool to help guide intuition. We introduce the concept of a transitivity tree in Section~\ref{section:transitivity-trees}, using connected trees of height at most $\omega$ as the foundation. We introduce a fourth type of transitivity point, and employ transitivity trees to determine the relationship between the four types.  Subsequent sections include: 
\begin{itemize}
\item $0$-transitive points (Section~\ref{section:0-trans-points}).
\item $2$-transitive and $3$-transitive points (Section~\ref{section:3-trans-not-2-trans}).
\item Dense orbit transitivity (Section~\ref{section:dense-orbit-transitivity}). 
\item Transitivity (Section~\ref{section:transitivity}). 
%\item Mixing (Section~\ref{section:mixing}). 
\end{itemize}

%---------------------------------------------------------------
\section{CR-dynamical systems}\label{section:CR-dynamical-systems}
In this section we provide definitions, notation and results that are required in the sequel.  

\begin{definition}
Given a compact metric space $X$, we denote by {\color{cyan} $\isolated{X}$} the set of isolated points in $X$. 
\end{definition}

\begin{definition}
A {\color{cyan} \emph{topological dynamical system}} is a pair $\p{X, f}$, where $f : X \to X$ is a continuous self-map on a compact metric space $X$. We let {\color{cyan} $\Graph{f} = \setOf{\p{x, y}}:{y = f\of{x}}$} denote the {\color{cyan} \emph{graph of the self-map $f$}}.  
\end{definition}

\begin{definition}
Let $X$ and $Y$ be compact metric spaces.  We denote by {\color{cyan} $2^X$} the collection of non-empty closed subsets of $X$. A {\color{cyan} \emph{set-valued function}}, is a function $f : X \to 2^Y$.  We say a set-valued function $f : X \to 2^Y$ is {\color{cyan} \emph{upper semi-continuous}}, abbreviated {\color{cyan} \emph{usc}}, if $\setOf{x \in X}:{f\of{x} \subseteq O}$ is open in $X$ for each open set $O$ of $Y$. We say a set-valued function $f : X \to 2^Y$ is {\color{cyan} \emph{lower semi-continuous}}, abbreviated ${\color{cyan} \emph{lsc}}$, if $\setOf{x \in X}:{f\of{x} \intersect O \neq \emptySet}$ is open in $X$ for each open set $O$ of $Y$. We say a set-valued function $f : X \to 2^Y$ is {\color{cyan} \emph{continuous}}, if $f$ is both usc and lsc. 
\end{definition}

\begin{definition}
A {\color{cyan} \emph{set-valued dynamical system}} is a pair $\p{X, F}$, where $F : X \to 2^X$ is an upper semi-continuous (usc) set-valued function on a compact metric space $X$. We let 
\[
{\color{cyan} \Graph{F} = \setOf{\p{x, y}}:{y \in F\of{x}}}
\]
denote the {\color{cyan} \emph{graph of our set-valued function $F$}}.
\end{definition}

%For each $x \in X$, we let 
%\[
%{\color{cyan} F^{-1}\of{x} = \setOf{y \in X}:{x \in F\of{y}}}. 
%\]

%The following result is well-known (see~\cite[Theorem $1.2$]{ingram}). 

%\begin{theorem}\label{theorem:ingram}
%Let $X$ be a compact metric space, and let $F : X \to 2^X$ be a set-valued function.  Then, $F : X \to 2^X$ is an usc set-valued function if, and only if, $\Graph{F}$ is closed in $X \times X$. 
%\end{theorem}

%By Theorem~\ref{theorem:ingram}, it is natural to consider non-empty closed subsets of $X \times X$ in place of continuous functions. Indeed, as noted in~\cite{transitivity_CR}, this is useful when one wants to consider $\p{X, F^{-1}}$, and $F^{-1}$ is not well-defined.  It is natural to consider $\Graph{F}^{-1} := \setOf{\p{x, y}}:{\p{y, x} \in \Graph{F}}$, but $\Graph{F}^{-1}$ need not correspond to the graph of an usc set-valued function on $X$ (although, it will correspond to the graph of an usc set-valued function on a closed subspace of $X$).  Observe, however, $\Graph{F}^{-1}$ is a non-empty closed subset of $X \times X$. To this end, we finally define CR-dynamical systems. 

\begin{definition}
A {\color{cyan} \emph{relation}} on $X$, is a non-empty subset of $X \times X$. A {\color{cyan} \emph{CR-dynamical system}} is a pair $\p{X, G}$, where $X$ is a compact metric space and $G$ is a closed relation (CR) on $X \times X$. If in addition there exists usc set-valued (SV) function $F : X \to 2^X$ such that $\Graph{F} = G$, we say $\p{X, G}$ is an {\color{cyan} \emph{SV-dynamical system}}. 
\end{definition}

\begin{observation}
Let $\p{X, G}$ be a CR-dynamical system. Then, $\p{X, G}$ is an SV-dynamical system if, and only if, there is a set-valued dynamical system $\p{X, F}$ such that $G = \Graph{F}$. 
\end{observation}

\begin{definition}
Suppose $G$ is a relation on a set $X$. Let {\color{cyan} $G^{-1}$} denote the set 
\[
G^{-1} = \setOf{\p{x, y}}:{\p{y, x} \in G}. 
\]

\end{definition}

\begin{observation}
If $\p{X, G}$ is a CR-dynamical system, then $\p{X, G^{-1}}$ is a CR-dynamical system. 
\end{observation}

For each $m \in \naturals$, we define {\color{cyan} $[m] = \set{0, \ldots, m}$}. 

\begin{definition}[{Definition $2.5$~\cite{transitivity_CR}}]
Suppose $\p{X, G}$ is a CR-dynamical system. For each non-negative integer $m$, we call
\[
\mahavier{m}{G} = \setOf{\sequence{x_0, x_1, \ldots, x_m} \in \prod_{n = 0}^m X}:{\text{for each } n \in [m - 1], \p{x_n, x_{n+1}} \in G}
\]
the {\color{cyan} \emph{$m$-th Mahavier product of $G$}}, and we call
\[
\mahavier{\infinity}{G} = \setOf{\sequenceOf{x_n}:{n \in \naturals} \in \prod_{n = 0}^\infinity X}:{\text{for each } n \in \naturals, \p{x_n, x_{n+1}} \in G}
\]
the {\color{cyan} \emph{infinite Mahavier product of $G$}}. 
\end{definition}

%The following results are well-known. 

%\begin{lemma}
%Suppose $\p{X, G}$ is a CR-dynamical system. Then, for each $m \in \naturals$, $\mahavier{m}{G}$ is closed in $\prod_{n = 0}^m X$. 
%\end{lemma}

%\begin{lemma}
%Suppose $\p{X, G}$ is a CR-dynamical system. Then, $\mahavier{\infinity}{G}$ is closed in $\prod_{n = 0}^\infinity X$. 
%\end{lemma}

The following two definitions are similar to Definition $2.7$ and Definition $2.8$ in~\cite{transitivity_CR}, except we start indexing at $0$ instead of $1$. 

\begin{definition}
Let $X$ be a compact metric space. For each $k \in \naturals$, we let $\pi_k : \prod_{n = 0}^\infinity X \to X$ denote the {\color{cyan} \emph{$k$-th standard projection}}. For each $n \in \naturals$, we also use $\pi_k : \prod_{i=0}^n X \to X$ to denote the {\color{cyan} \emph{$k$-th standard projection}}, for each $k \in [n]$.  
\end{definition}

\begin{definition}
Suppose $G$ is a relation on a compact metric space $X$. Let $x \in X$ and $A \subseteq X$. We define
\begin{itemize}
\item $G\of{x} = \setOf{y \in X}:{\p{x, y} \in G}$;
\item $G\of{A} = \Union_{y \in A} G\of{y}$.
\end{itemize}
Furthermore, for each positive integer $n$, we define 
\[
G^n\of{x} = \setOf{y \in X}:{\text{there exists } \bm{x} \in \mahavier{n}{G} \text{ such that } \pi_0\of{\bm{x}} = x \text{ and } \pi_n\of{\bm{x}} = y}
\]
and
\[
G^n\of{A} = \Union_{y \in A} G^n\of{y}. 
\]
We then define $G^n = \setOf{\p{x, y} \in X \times X}:{y \in G^n\of{x}}$ for each $n \geq 1$. We use the convention 
\begin{itemize}
\item {$G^0\of{x} = \set{x}$;}
\item {$G^0\of{A} = A$;}
\item {$G^0 = \Delta_X$;}
\end{itemize}
where $\Delta_X := \setOf{\p{x, x}}:{x \in X}$ denotes the {\color{cyan} \emph{diagonal}} of $X$.  Moreover, for each $n \in \naturals$, we define
\begin{itemize}
\item {$G^{-n}\of{x} = \p{G^{-1}}^n\of{x}$};
\item {$G^{-n}\of{A} = \p{G^{-1}}^n\of{A}$};
\item {$G^{-n}[A] = \setOf{y \in X}:{G^n\of{y} \subseteq A}$.}
\end{itemize}
\end{definition}

We now state the following observations from~\cite[Section $2$]{transitivity_CR}. 

\begin{observation}[{Observation $2.9$~\cite{transitivity_CR}}]
Let $\p{X, G}$ be a CR-dynamical system, $n \in \naturals$, and $x \in X$. Then, $G^{n+1}\of{x} = G\of{G^n\of{x}}$. 
\end{observation}

\begin{observation}[{Observation $2.10$~\cite{transitivity_CR}}]
Let $\p{X, G}$ be a CR-dynamical system, $n \in \naturals$, and $x, y \in X$. Then, the following statements are equivalent. 
\begin{enumerate}
\item $x \in G^n\of{y}$. 
\item There is $\sequence{x_0, \ldots, x_n} \in \mahavier{n}{G}$ such that $x_0 = y$ and $x_n = x$. 
\item There is $\sequence{x_0, \ldots, x_n} \in \mahavier{n}{G^{-1}}$ such that $x_0 = x$ and $x_n = y$. 
\end{enumerate}
\end{observation}

\begin{observation}[{Observation $2.11$~\cite{transitivity_CR}}]
Let $\p{X, G}$ be a CR-dynamical system, $n \in \naturals$, and $x, y \in X$. Then, the following statements are equivalent. 
\begin{enumerate}
\item $x \in G^{-n}\of{y}$. 
\item There is $\sequence{x_0, \ldots, x_n} \in \mahavier{n}{G}$ such that $x_0 = x$ and $x_n = y$. 
\item There is $\sequence{x_0, \ldots, x_n} \in \mahavier{n}{G^{-1}}$ such that $x_0 = y$ and $x_n = x$. 
\end{enumerate}
\end{observation}

We now make a further observation. 

%\begin{observation}
%Let $\p{X, G}$ be a CR-dynamical system, $n \in \omega$, and $x \in X$. Then, $G^{-\p{n+1}}[x] = G^{-1}[G^{-n}[x]]$. 
%\end{observation}

\begin{observation}
Let $\p{X, G}$ be a CR-dynamical system, $n \in \naturals$, and $x, y \in X$. Then, the following statements are equivalent. 
\begin{enumerate}
\item $x \in G^{-n}[\set{y}]$. 
\item If $\sequence{x_0, \ldots, x_n} \in \mahavier{n}{G}$ and $x_0 = x$, then $x_n = y$. 
\end{enumerate}
\end{observation}

We now recall the concepts of a trajectory and an orbit from topological dynamical systems.  We conclude this section by defining them for CR-dynamical systems as found in~\cite{transitivity_CR}. Refer to~\cite{transitivity_CR, minimality_CR} for basic results. 

\begin{definition}
The ${\color{cyan} \emph{trajectory}}$ of a point $x \in X$ in a topological dynamical system $\p{X, f}$, is the sequence {\color{cyan} $\sequence{x, f\of{x}, f^2\of{x}, \ldots}$}. The {\color{cyan} \emph{orbit}} of a point $x \in X$ in a topological dynamical system $\p{X, f}$ is the set {\color{cyan} $\set{x, f\of{x}, f^2\of{x}, \ldots}$}.
\end{definition}

\begin{definition}
A {\color{cyan} \emph{trajectory}} of $x \in X$ in a CR-dynamical system $\p{X, G}$, is a sequence $\bm{x} \in \mahavier{\infinity}{G}$ such that $\pi_0\of{\bm{x}} = x$. We let {\color{cyan} $T_G^+\of{x}$} denote the set of trajectories of $x$ in $\p{X, G}$. We say $x \in X$ {\color{cyan} \emph{is a legal point of $\p{X, G}$}}, if $T_G^+\of{x}$ is non-empty. Otherwise, we say $x \in X$ {\color{cyan} \emph{is an illegal point of $\p{X, G}$}}. We let {\color{cyan} $\legal{G}$} denote the set of legal points, and {\color{cyan} $\illegal{G}$} denote the set of illegal points in $\p{X, G}$. 
\end{definition}

\begin{definition}
Let $\p{X, G}$ be a CR-dynamical system and $x \in X$. Suppose $\bm{x} \in T_G^+\of{x}$. Then, we call {\color{cyan} $\mathcal{O}_G^\oplus\of{\bm{x}}$} {\color{cyan} \emph{a forward orbit of $x$}}, where
\[
\mathcal{O}_G^\oplus\of{\bm{x}} = \setOf{\pi_k\of{\bm{x}}}:{k \in \naturals}. 
\]
We denote by {\color{cyan} $\mathcal{U}^\oplus_G\of{x}$} the set
\[
\mathcal{U}^\oplus_G\of{x} = \Union_{\bm{y} \in T_G^+\of{x}} \mathcal{O}_G^\oplus\of{\bm{y}}. 
\]
\end{definition}

%---------------------------------------------------------------
\section{Transitivity trees}\label{section:transitivity-trees}
For a topological dynamical systems $\p{X, f}$, each point $x \in X$ has exactly one trajectory. In a CR-dynamical system $\p{X, G}$, a given point may have multiple trajectories, or no trajectory at all.  In this section we introduce transitivity trees which allow us to view all trajectories of a point simultaneously, and without losing information about the progression to any coordinate of a trajectory of $x$. 

We first recall definitions relating to trees. We refer the reader to~\cite{roitman_set_theory} for an introduction to set-theoretic trees.  We restrict our interest to connected trees of height at most $\omega$. We will use $\infinity$ in place of $\omega$.

\begin{definition}
Let $\p{T, \leq}$ be a partially ordered set. We say that $T$ is a {\color{cyan}\emph{connected tree}}, if there exists unique $r \in T$, called the {\color{cyan}\emph{root}} of $T$, such that for each $x \in T$:
\begin{itemize}
\item $r \leq x$; and 
\item $\p{\setOf{y \in T}:{y \leq x}, \leq}$ is well-ordered. 
\end{itemize}
\end{definition}

\begin{definition}
Suppose $\p{T, \leq}$ is a connected tree and $x \in T$.  We define 
\[
x^+ = \setOf{y \in T}:{x \leq y}. 
\]
Furthermore, we say $y$ is a {\color{cyan} \emph{successor}} of $x$, if $x < y$ and $x \leq z \leq y$ implies $z \in \set{x, y}$. Finally, a {\color{cyan} \emph{leaf of $T$}} is a point with no successor. 
\end{definition}

\begin{definition}
Suppose $\p{T, \leq}$ is a tree. If $B$ is a maximal well-ordered set in $T$, we call $B$ a {\color{cyan} \emph{branch of $T$}}. If $\cardinality{B} = n$, we say {\color{cyan} \emph{the height of $B$ is $n - 1$}}, and write {\color{cyan} $\heightIn{T}{B} = n - 1$}. If $n$ is finite, we say that the {\color{cyan} \emph{height of $B$ is finite}} and write {\color{cyan} $\heightIn{T}{B} < \infinity$}. Otherwise, we say the {\color{cyan} \emph{height of $B$ is infinite}} and write {\color{cyan} $\heightIn{T}{B} = \infinity$}. We denote by {\color{cyan} $\mathcal{B}\of{T}$} the collection of branches of $T$, and by {\color{cyan} $\mathcal{B}_\infinity\of{T}$} the collection of infinite branches of $T$. We further define the {\color{cyan} \emph{height}} of $\p{T, \leq}$ as 
\[
{\color{cyan} \height{T} = \sup\setOf{\heightIn{T}{B}}:{B \in \mathcal{B}\of{T}}}. 
\]
\end{definition}

Observe that every branch contains the root of the tree.  

\begin{figure}[h]
\[
 \begin{tikzpicture}[
       decoration = {markings,
                     mark=at position .5 with {\arrow{Stealth[length=2mm]}}},
       dot/.style = {circle, fill, inner sep=2.4pt, node contents={},
                     label=#1},
     dots/.style = {circle, fill, inner sep=0.0pt, node contents={},
                     label=#1},
every edge/.style = {draw, postaction=decorate}
                        ]
\node[cyan] (x0) [dot=below:$r$];
\node[orange] (x1) at (-1, 1) [dot];
\node[orange] (x2) at (1, 1) [dot];
\node[orange] (x3) at (0, 1) [dot];
\node[violet] (x4) at (-2, 2) [dot];
\node[violet] (x5) at (0, 2) [dot];
\node[violet] (x6) at (1, 2) [dot];
\node[purple] (x7) at (1, 3) [dot];
\node (x8) at (1, 4) {$\vdots$};
\node (x10) at (4, 4) {$\vdots$};
\node[red] (x9) at (1, 4.75) [dot];

\node (l0) at (4, 0) {$\begin{small}\color{cyan}\text{level$_0\of{T}$}\end{small}$};

\node (l1) at (4, 1) {$\begin{small}\text{\color{orange}level$_1\of{T}$}\end{small}$};

\node (l2) at (4, 2) {$\begin{small}\color{violet}\text{level$_2\of{T}$}\end{small}$};

\node (l3) at (4, 3) {$\begin{small}\color{purple}\text{level$_3\of{T}$}\end{small}$};

\node (l4) at (4, 4.75) {$\begin{small}\color{red}\text{level$_n\of{T}$}\end{small}$};

\draw (5.25, -.125) -- (5.25, 2.125);
\draw (5, -.125) -- (5.25, -.125);
\draw (5, 2.125) -- (5.25, 2.125);
\draw (5.25, 1) -- (5.5, 1); 

\node (b0) at (6.25, 1) {$\begin{small}\text{$\mathcal{L}_2\of{T}$}\end{small}$};

\draw   (x0) -- (x1);
\draw	(x0) -- (x2);   
\draw   (x0) -- (x3);
\draw   (x1) -- (x4);
\draw   (x1) -- (x5);
\draw   (x6) -- (x2);
\draw (x6) -- (x7);
\draw (x7) -- (1, 3.5);
\draw (1, 4.25) -- (x9);
 \end{tikzpicture}
\]
\caption{Levels in a tree}\label{figure:tree-levels}
\end{figure}
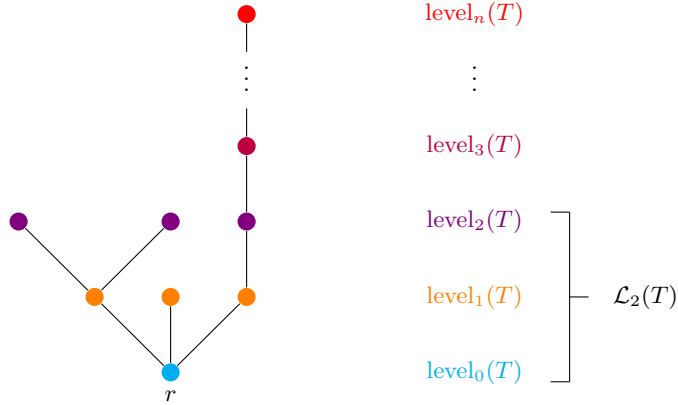

\begin{definition}
Let $\p{T, \leq}$ be a tree. Let $n \in \naturals$. Then we define {\color{cyan} $\levelIn{T}{n}$} by
\[
\levelIn{T}{n} = \setOf{x \in T}:{n = \cardinality{\setOf{y \in T}:{y \leq x}} - 1}. 
\]
We also define the set {\color{cyan} $\mathcal{L}_n\of{T}$} by
\[
\mathcal{L}_n\of{T} = \Union_{i = 0}^n \levelIn{T}{i}. 
\]
See Figure~\ref{figure:tree-levels}.
\end{definition}

%We now turn our attention back to CR-dynamical systems, introducing the concept of paths between points with respect to our closed relation.

We now define transitivity trees.

\begin{definition}
Let $\p{X, G}$ be a CR-dynamical system, $x, y \in X$ and $n \in \naturals$. We say {\color{cyan} \emph{$\gamma = x_0 \ldots x_n$ is an $x$-path in $G$}}, if
\begin{itemize}
\item $x_0 = x$; and
\item $\sequence{x_0, \ldots, x_n} \in \mahavier{n}{G}$. 
\end{itemize}
We denote by ${\color{cyan} x_\gamma}$ the endpoint of our path $\gamma$, i.e., $x_\gamma = x_n$. If in addition $y = x_\gamma$, we say {\color{cyan} \emph{$\gamma = x_0 \ldots x_n$ is a path from $x$ to $y$ in $G$}}. Furthermore, we say the {\color{cyan} \emph{length of $\gamma = x_0 \ldots x_n$ is $n$}}. 

We say an $x$-path $\gamma_1$ {\color{cyan} \emph{extends}} to an $x$-path $\gamma_2$, if $\gamma_2 = \gamma_1 \gamma_3$ for some path $\gamma_3$ in $G$. We denote by {\color{cyan}$T_G\of{x}$} the set of $x$-paths in $G$.
\end{definition}
\begin{note}
We observe $\p{T_G\of{x}, \leq}$ is a tree, where for each $\gamma_1, \gamma_2 \in T_G\of{x}$,
\[
\gamma_1 \leq \gamma_2 \Longleftrightarrow \text{$\gamma_1$ extends to $\gamma_2$, or $\gamma_1 = \gamma_2$.}
\]
It is easily checked $\p{T_G\of{x}, \leq}$ is a partially ordered set. Furthermore, $x$ is the root of $T_G\of{x}$, and if $\gamma = x_0 \ldots x_n \in T_G\of{x}$, then 
\[
\setOf{\gamma^\prime \in T_G\of{x}}:{\gamma^\prime \leq \gamma} = \set{x_0, x_0 x_1, \ldots, x_0 \ldots x_n}, 
\]
which is clearly well-ordered. 
\end{note}

\begin{definition}
Let $\p{X, G}$ be a CR-dynamical system. The {\color{cyan} \emph{transitivity tree of $\p{X, G}$ with respect to $x$}}, is the tree $\p{T_G\of{x}, \leq}$, where $T_G\of{x}$ is the set of $x$-paths in $G$ and $\leq$ is the path extension order. 
\end{definition}

Observe that for each $x_\gamma \in T_G\of{x}$, as the index $\gamma$ is the path from $x$ to $x_\gamma$, the index retains the information of the place $x_\gamma$ in the trajectory of $x$ containing it. 

%\begin{note}
%Observe there may exist $\gamma \neq \gamma^\prime \in T_G\of{x}$ such that $x_\gamma = x_{\gamma^\prime}$ in $X$ (see Example~\ref{ex:tails-same}).  To go from members of $T_G\of{x}$ to members of $X$, we define the following operation $\ast$. 
%\end{note}

\begin{definition}
Let $\p{X, G}$ be a CR-dynamical system and $x \in X$. For each $S \subseteq T_G\of{x}$, we denote by {\color{cyan} $S^\ast$} the set
\[
S^\ast = \setOf{x_\gamma \in X}:{\gamma \in S}. 
\]
\end{definition}

%\begin{figure}[h]
%\[
% \begin{tikzpicture}[
%       decoration = {markings,
%                     mark=at position .5 with {\arrow{Stealth[length=2mm]}}},
%       dot/.style = {circle, fill, inner sep=2.4pt, node contents={},
%                     label=#1},
%every edge/.style = {draw, postaction=decorate}
%                        ]
%\node[cyan] (x0) [dot=above:$0$];
%\node (x1) at (-2, 1) [dot=above:$0$];
%\node (y0) at (2, 1) [dot=above:$1$];
%\node (x2) at (-3.5, 2) [dot=above:$0$];
%\node (y2) at (-0.5, 2) [dot=above:$1$];
%\node (x3) at (3.5, 2) [dot=above:$1$];
%\node (y3) at (0.5, 2) [dot=above:$0$];

%\draw   (x0) -- (x1);
%\draw	(x0) -- (y0); 
%\draw (x1) -- (x2);
%\draw (x1) -- (y2); 
%\draw (y0) -- (x3);
%\draw (y0) -- (y3); 
% \end{tikzpicture}
% \]
% \caption{$\mathcal{L}_2\of{T_G\of{0}}$ in Example~\ref{ex:level-2-binary-tree}}\label{fig:level-2-binary-tree}
%\end{figure}

%\begin{example}\label{ex:level-2-binary-tree}
%Let $X = \set{0, 1}$ and $G = X \times X$.  Then, $\mathcal{L}_2\of{T_G\of{0}}$ is depicted in Figure~\ref{fig:level-2-binary-tree}. We note $T_G\of{0}$ consists of all binary strings starting at $0$, and $T_G\of{1}$ all binary strings starting at $1$. 
%\end{example}

We now make a few simple observations on transitivity trees. 

\begin{observation}
Let $\p{X, G}$ be a CR-dynamical system and $x \in X$. If $y \in T_G\of{x}^\ast$, then $T_G\of{y}^\ast \subseteq T_G\of{x}^\ast$. 
\end{observation}

\begin{observation}
Let $\p{X, G}$ be a CR-dynamical system and $x, y \in X$. Then, $y \in T_G\of{x}^\ast$ if, and only if, $x \in T_{G^{-1}}\of{y}^\ast$. 
\end{observation}

\begin{observation}\label{observation:branch-extension}
Let $\p{X, G}$ be a CR-dynamical system and $x \in X$. If $y \in T_{G^{-1}}\of{x}^\ast$ and $B_x \in \mathcal{B}_\infinity\of{T_G\of{x}}$, then there exists $B_y \in \mathcal{B}_\infinity\of{T_G\of{y}}$ such that $B_x^\ast \subseteq B_y^\ast$.  
\end{observation}

\begin{observation}
Suppose $\p{X, G}$ is a CR-dynamical system and $x \in X$. Let $\sequenceOf{x_n}:{n \in \naturals} \in \prod_{n \in \naturals} X$. Then the following are equivalent. 
\begin{enumerate}
    \item $\sequenceOf{x_n}:{n \in \naturals} \in T_G^+\of{x}$.
    \item For each $n \in \naturals$, there exists $\gamma_n = x_0 \ldots x_n \in T_G\of{x}$, and $B = \setOf{\gamma_n}:{n \in \naturals} \in \mathcal{B}_\infinity\of{T_G\of{x}}$. 
\end{enumerate}
Furthermore, there is a one-to-one correspondence between trajectories of $x$ in $\p{X, G}$ and infinite branches in $T_G\of{x}$. 
\end{observation}

\begin{observation}
Suppose $\p{X, G}$ is a CR-dynamical system and $x \in X$. Then, for each $n \in \naturals$, 
\[
\levelIn{T_G\of{x}}{n}^\ast = G^n\of{x}
\]
and 
\[
\mathcal{L}_n\of{T_G\of{x}}^\ast = \Union_{k \in [n]} G^k\of{x}.  
\]
Similarly, for each $n \in \naturals$, 
\[
\levelIn{T_{G^{-1}}\of{x}}{n}^\ast = G^{-n}\of{x}
\]
and 
\[
\mathcal{L}_n\of{T_{G^{-1}}\of{x}}^\ast = \Union_{k \in [n]} G^{-k}\of{x}.  
\]
\end{observation}
\begin{note}
By~\cite[Lemma $3.13$]{transitivity_CR}, $\levelIn{T_G\of{x}}{n}^\ast$,  $\mathcal{L}_n\of{T_G\of{x}}^\ast$, $\levelIn{T_{G^{-1}}\of{x}}{n}^\ast$, and $\mathcal{L}_n\of{T_{G^{-1}}\of{x}}^\ast$ are closed in $X$, for each $n \in \naturals$. 
\end{note}

\begin{observation}
Suppose $\p{X, G}$ is a CR-dynamical system and $x \in X$. Then, if $B \in \mathcal{B}_\infinity\of{T_G\of{x}}$ is the infinite branch corresponding to $\bm{x} \in T_G^+\of{x}$, then 
\[
B^\ast = \mathcal{O}_G^\oplus\of{\bm{x}}. 
\]
Furthermore,
\begin{align*}
&& \Union_{B \in \mathcal{B}_\infinity\of{T_G\of{x}}} B^\ast &= \mathcal{U}^\oplus_G\of{x}, && \\
&& T_G\of{x}^\ast &= \Union_{n \in \naturals} G^n\of{x}, && \\
\text{and}
&& T_{G^{-1}}\of{x}^\ast &= \Union_{n \in \naturals} G^{-n}\of{x}. &&
\end{align*}
Finally, we note $T_G\of{x}^\ast \intersect \legal{G} = \Union_{B \in \mathcal{B}_\infinity\of{T_G\of{x}}} B^\ast$. 
\end{observation}

\begin{figure}
\[
 \begin{tikzpicture}[
       decoration = {markings,
                     mark=at position .5 with {\arrow{Stealth[length=2mm]}}},
       dot/.style = {circle, fill, inner sep=2.4pt, node contents={},
                     label=#1},
every edge/.style = {draw, postaction=decorate}
                        ]
\node (x0) [dot];
\node (x1) at (0,1) [dot];
\node (y1) at (1, 1) [dot];
\node (y2) at (1, 2) [dot];
\node (z1) at (2, 1) [dot];
\node (z2) at (2, 2) [dot];
\node (z3) at (2, 3) [dot];
\node (d) at (3, 2) {$\ldots$};

\draw   (x0) -- (x1);
\draw	(x0) -- (y1);   
\draw   (x0) -- (z1);
\draw   (y1) -- (y2);
\draw  (z1) -- (z2);
\draw  (z2) -- (z3);
 \end{tikzpicture}
\]
\caption{An infinite tree with no infinite branches}\label{fig:omega-tree-no-infinite-branches}
\end{figure}
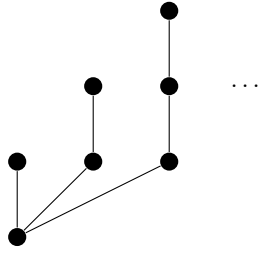

There exist infinite trees with no infinite branches (see Figure~\ref{fig:omega-tree-no-infinite-branches}). By the following observation, such a tree is not a transitivity tree.  

\begin{observation}
Suppose $\p{X, G}$ is a CR-dynamical system and $x \in X$. Then, the following are equivalent.  
\begin{enumerate}
\item $x \in \legal{G}$;
\item $\height{T_G\of{x}} = \infinity$; and
\item $\mathcal{B}_\infinity\of{T_G\of{x}} \neq \emptySet$.
\end{enumerate}
\end{observation}
\begin{note}
Follows from~\cite[Theorem $3.12$]{transitivity_CR}. 
\end{note}

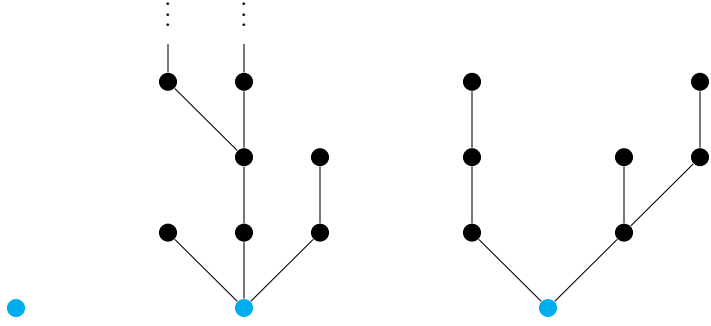
\begin{figure}
\[
 \begin{tikzpicture}[
       decoration = {markings,
                     mark=at position .5 with {\arrow{Stealth[length=2mm]}}},
       dot/.style = {circle, fill, inner sep=2.4pt, node contents={},
                     label=#1},
every edge/.style = {draw, postaction=decorate}
                        ]
\node[cyan] (x0) [dot];

\node[cyan] (x00) at (3, 0) [dot];
\node (x1) at (3,1) [dot];
\node (y0) at (3, 2) [dot];
\node (y1) at (3, 3) [dot];
\node (w0) at (3, 4) {\vdots};
\node (x) at (4, 1) [dot];
\node (y) at (4, 2) [dot];
\node (z) at (2, 3) [dot];
\node (l1) at (2, 4) {\vdots};
\node (a) at (2, 1) [dot];

\draw   (x00) -- (x1);
\draw	(x1) -- (y0);   
\draw   (y0) -- (y1);
\draw   (x00) -- (x);
\draw   (x) -- (y);
\draw   (y0) -- (z);
\draw   (x00) -- (a);
\draw (z) -- (2, 3.5);
\draw (y1) -- (3, 3.5);

\node[cyan] (x000) at (7, 0) [dot];
\node (x11) at (8,1) [dot];
\node (y11) at (6, 1) [dot];
\node (x22) at (8, 2) [dot];
\node (y22) at (6, 2) [dot];
\node (z22) at (9, 2) [dot];
\node (z33) at (9, 3) [dot];
\node (y33) at (6, 3) [dot];

\draw   (x000) -- (x11);
\draw	(x000) -- (y11);   
\draw   (x11) -- (x22);
\draw   (y11) -- (y22);
\draw   (y22) -- (y33);
\draw  (x11) -- (z22);
\draw  (z22) -- (z33);
 \end{tikzpicture}
\]
\caption{Types of transitivity trees}\label{fig:types-of-transitivity-trees}
\end{figure}

\begin{definition}
Suppose $\p{X, G}$ is a CR-dynamical system and $x \in X$. We say that {\color{cyan} \emph{$x$ is degenerate in $\p{X, G}$}}, if $T_G\of{x}$ is a singleton (see the first tree in Figure~\ref{fig:types-of-transitivity-trees}). Otherwise, we say that {\color{cyan} \emph{$x$ is non-degenerate in $\p{X, G}$}}. We denote by {\color{cyan} $\degenerate{G}$} the set of degenerate points in $\p{X, G}$. We also denote by  {\color{cyan} $\nondegenerate{G}$} the set of non-degenerate points in $\p{X, G}$. 
\end{definition}

Observe:
\begin{itemize}
\item The leaves of a transitivity tree correspond to degenerate points in $\p{X, G}$. 
\item{$\legal{G} \subseteq \nondegenerate{G}$ and $\degenerate{G} \subseteq \illegal{G}$. In Figure~\ref{fig:types-of-transitivity-trees}, see the second tree for what the transitivity tree of a legal point could look like, and the third tree for what the transitivity tree of an illegal non-degenerate point could look like.}
\item If every point is non-degenerate, then every point is legal.
\item $X = \nondegenerate{G} \union \degenerate{G}$. 
\end{itemize}

\begin{proposition}\label{prop:nondegenerate-closed}
Let $\p{X, G}$ be a CR-dynamical system. Then, $\nondegenerate{G}$ is closed in $X$. Equivalently, $\degenerate{G}$ is open in $X$. 
\end{proposition}
\begin{proof}
Observe $\nondegenerate{G} = \pi_0\of{G}$, and that $\pi_0$ is a closed map. 
\end{proof}

\begin{proposition}[{Theorem $3.15$~\cite{transitivity_CR}}]
Let $\p{X, G}$ be a CR-dynamical system. Then, $\legal{G}$ is closed in $X$, and equivalently $\illegal{G}$ is open in $X$. 
\end{proposition}

\begin{definition}
Let $\p{X, G}$ be a CR-dynamical system.  A trajectory of the form 
\[
\sequence{x_0, \ldots, x_n, x_0, \ldots, x_n, x_0, \ldots, x_n, \ldots}
\]
is {\color{cyan} \emph{cyclic}}, and has {\color{cyan} \emph{order $n$}}. 
\end{definition}

We conclude this section by defining a useful equivalence relation, as a way to distinguish points with cyclic trajectories.

\begin{definition}\label{def:transitivity-tree-relation}
Let $\p{X, G}$ be a CR-dynamical system. We define a relation {\color{cyan} $\sim_G$} on $X$. For each $x, y \in X$, $x \sim_G y$ if, and only if, $x \in T_G\of{y}^\ast$ and $y \in T_G\of{x}^\ast$. 
\end{definition}

\begin{lemma}
Let $\p{X, G}$ be a CR-dynamical system. The relation $\sim_G$ on $X$ is an equivalence relation. 

Denote by $[x]_G$, the equivalence class of $x$ with respect to $\sim_G$. Then the following hold. 
\begin{enumerate}
\item For each $x \in X$, $[x]_G = T_G\of{x}^\ast \intersect T_{G^{-1}}\of{x}^\ast = [x]_{G^{-1}}$. 
\item $x \sim_G y$ if, and only if, $T_G\of{x}^\ast = T_G\of{y}^\ast$ and $T_{G^{-1}}\of{x}^\ast = T_{G^{-1}}\of{y}^\ast$.
\end{enumerate}
\end{lemma}
\begin{proof}
Both reflexivity and symmetry are clear. To see why $\sim_G$ is transitive, suppose $x \sim_G y$ and $y \sim_G z$. Then there is a path from $x$ to $y$ in $G$, and a path from $y$ to $z$ in $G$. Hence, there is a path from $x$ to $z$ in $G$. On the other hand, there is a path from $z$ to $y$ in $G$, and a path from $y$ to $x$ in $G$. Thus, there is a path from $z$ to $x$ in $G$, and therefore $x \sim_G z$. 

$\p{1}$ is clear from the definition. To prove $\p{2}$, if $x \sim_G y$, then $x \in T_G\of{y}^\ast$ implies $T_G\of{x}^\ast \subseteq T_G\of{y}^\ast$. Similarly, $y \in T_G\of{x}^\ast$ implies $T_G\of{y}^\ast \subseteq T_G\of{x}^\ast$. As $x \sim_{G^{-1}} y$, it follows $T_{G^{-1}}\of{x}^\ast = T_{G^{-1}}\of{y}^\ast$. The converse holds by our observation $[x]_G = T_G\of{x}^\ast \intersect T_{G^{-1}}\of{x}^\ast$. 
\end{proof}

%Let $\p{X, G}$ be a CR-dynamical system. We define an equivalence relation $\sim_G$ on $X$. For each $x, y \in X$, 
%\[
%x \sim_G y \Longleftrightarrow x \in T_G\of{y}^\ast \text{ and } y \in T_G\of{x}^\ast. 
%\]
%Indeed, both reflexivity and symmetry are clear. To see why $\sim_G$ is transitive, suppose $x \sim_G y$ and $y \sim_G z$. Then there is a path from $x$ to $y$ in $G$, and a path from $y$ to $z$ in $G$. Hence, there is a path from $x$ to $z$ in $G$. On the other hand, there is a path from $z$ to $y$ in $G$, and a path from $y$ to $x$ in $G$. Thus, there is a path from $z$ to $x$ in $G$, and therefore $x \sim_G z$. 

%We denote by $[x]_G$, the equivalence class of $x$ with respect to $\sim_G$. Observe 
%\[
%[x]_G = T_G\of{x}^\ast \intersect T_{G^{-1}}\of{x}^\ast = [x]_{G^{-1}}
%\]
%for each $x \in X$. 

%We also note that $x \sim_G y$ if, and only if, $T_G\of{x}^\ast = T_G\of{y}^\ast$ and $T_{G^{-1}}\of{x}^\ast = T_{G^{-1}}\of{y}^\ast$. For if $x \sim_G y$, then $x \in T_G\of{y}^\ast$ implies $T_G\of{x}^\ast \subseteq T_G\of{y}^\ast$. Similarly, $y \in T_G\of{x}^\ast$ implies $T_G\of{y}^\ast \subseteq T_G\of{x}^\ast$. As $x \sim_{G^{-1}} y$, it follows $T_{G^{-1}}\of{x}^\ast = T_{G^{-1}}\of{y}^\ast$. The converse holds by our observation $[x]_G = T_G\of{x}^\ast \intersect T_{G^{-1}}\of{x}^\ast$. 

%---------------------------------------------------------------
\section{$0$-transitive points}\label{section:0-trans-points}
When studying transitivity properties of topological dynamical systems, it often reduces to the study of transitive points; for instance, a well-known result is that a topological dynamical system $\p{X, f}$ on a compact metric space $X$ with no isolated points is transitive (see Definition~\ref{definition:transitive-top-dyn-system} in Section~\ref{section:transitivity}) if, and only if, there exists a transitive point~\cite{akin2016variationsconcepttopologicaltransitivity}. Classically, transitive points are defined as follows. 

\begin{definition}
Let $\p{X, f}$ be a topological dynamical system. We say $x \in X$ is a {\color{cyan} \emph{transitive point in $\p{X, f}$}}, if its orbit $\set{x, f\of{x}, f^2\of{x}, \ldots}$ is dense in $X$. We denote by {\color{cyan} $\text{tr}\of{f}$} the set of transitive points in $\p{X, f}$. 
\end{definition}

In this section, we study transitive points in CR-dynamical systems $\p{X, G}$. Transitive points are generalised from topological dynamical systems to CR-dynamical systems on compact metric spaces by Banic et al. in~\cite{transitivity_CR}, where they introduced three different types of transitive points.  In this section we introduce a new type of transitive point (\emph{$0$-transitive point}). We now provide an equivalent definition of~\cite[Definition $3.18$]{transitivity_CR} below, making use of transitivity trees. 

\begin{definition}[{Definition $3.18$~\cite{transitivity_CR}}]\label{def:transitive-points}
Let $\p{X, G}$ be a CR-dynamical system, and $T_G\of{x}$ be the transitivity tree of $x \in \legal{G}$. We say that
\begin{itemize}
    \item $x$ is a {\color{cyan} \emph{$1$-transitive point}} in $\p{X, G}$, if for each $B \in \mathcal{B}_\infinity\of{T_G\of{x}}$, $\closure{B^\ast} = X$. We denote by {\color{cyan} $\transitiveIn{G}{1}$} the set of $1$-transitive points in $\p{X, G}$. 
    \item $x$ is a {\color{cyan} \emph{$2$-transitive point}} in $\p{X, G}$, if there exists $B \in \mathcal{B}_\infinity\of{T_G\of{x}}$ such that $\closure{B^\ast} = X$. We denote by {\color{cyan} $\transitiveIn{G}{2}$} the set of $2$-transitive points in $\p{X, G}$.
    \item $x$ is a {\color{cyan} \emph{$3$-transitive point}} in $\p{X, G}$, if $\closure{\Union_{B \in \mathcal{B}_\infinity\of{T_G\of{x}}} B^\ast} = X$. We denote by {\color{cyan} $\transitiveIn{G}{3}$} the set of $3$-transitive points in $\p{X, G}$.
    \item $x$ is an {\color{cyan} \emph{intransitive point}} in $\p{X, G}$, if $x$ is not $3$-transitive. We denote by {\color{cyan} $\intrans{G}$} the set of intransitive points in $\p{X, G}$. 
\end{itemize}
\end{definition}

We now restate Observations $3.19$, $3.20$, and $3.21$ from~\cite{transitivity_CR} below, and then give the definition for a $0$-transitive point. 

\begin{observation}[{Observations $3.19$, $3.20$, and $3.21$~\cite{transitivity_CR}}]\label{observation:3.19-3.20-3.21}
Let $\p{X, G}$ be a CR-dynamical system, and let $\p{X, f}$ be a topological dynamical system.  Then, 
\begin{enumerate}
\item $\transitiveIn{G}{1} \subseteq \transitiveIn{G}{2} \subseteq \transitiveIn{G}{3}$;
\item $\legal{G} = \transitiveIn{G}{3} \union \intrans{G}$; and
\item $\text{tr}\of{f} = \transitiveIn{\Graph{f}}{1} = \transitiveIn{\Graph{f}}{2} = \transitiveIn{\Graph{f}}{3}$. 
\end{enumerate}
Therefore, $\transitiveIn{G}{3} \union \intrans{G}$ is closed.  
\end{observation}

\begin{definition}\label{def:0-transitive-point}
Let $\p{X, G}$ be a CR-dynamical system, and $T_G\of{x}$ be the transitivity tree of $x \in \legal{G}$. We say that $x$ is a {\color{cyan} \emph{$0$-transitive point}} in $\p{X, G}$, if for each non-empty open set $U$ in $X$ there exists $n \in \naturals$ such that $\levelIn{T_G\of{x}}{n}^\ast \subseteq U$. We denote by {\color{cyan} $\transitiveIn{G}{0}$} the set of $0$-transitive points in $\p{X, G}$. 
\end{definition}

We now make an observation, analogous to $\p{1}$ and $\p{3}$ in Observation~\ref{observation:3.19-3.20-3.21}.  

\begin{observation}\label{observation:how-transitive-points-contained}
Let $\p{X, G}$ be a CR-dynamical system, and let $\p{X, f}$ be a topological dynamical system.  Then, 
\begin{enumerate}
\item $\transitiveIn{G}{0} \subseteq \transitiveIn{G}{1} \subseteq \transitiveIn{G}{2} \subseteq \transitiveIn{G}{3}$; and
\item $\text{tr}\of{f} = \transitiveIn{\Graph{f}}{k}$ for each $k \in \set{0, 1, 2, 3}$. 
\end{enumerate}
\end{observation}

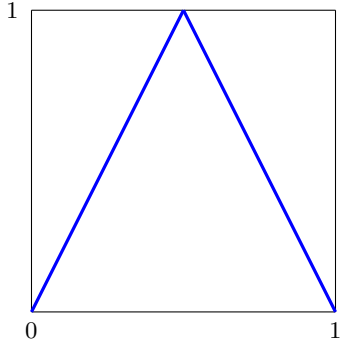
\begin{figure}[h]
\[
 \begin{tikzpicture}[
       decoration = {markings,
                     mark=at position .5 with {\arrow{Stealth[length=2mm]}}},
       dot/.style = {circle, fill, inner sep=1pt, node contents={},
                     label=#1},
every edge/.style = {draw, postaction=decorate}
                        ]
%\node (a) [dot];
%\node (b) at (5, 0) [dot];
%\node (c) at (5, 5) [dot]; 
%\node (d) at (0, 5) [dot];
%\node (e) at (8, 3) {$\color{red}\text{$G = X \times \set{0} \union \set{0} \times X$}$};
%\node (g) at (6.5, 4) {${\color{red} \text{$X = [0, 1]$}}$};
%\node (f) at (6.68, 2) {${\color{red} \text{$G\of{\Fraction 1 / 2} = \set{0}$}}$};
%\node (h) at (6.55, 1) {${\color{red} \text{$G^2\of{\Fraction 1 / 2} = X$}}$}

%\draw[dashed] (0, 2) -- (4, 2);
%\draw[dashed] (2, 0) -- (2, 4);
%\draw[dashed] (0, 2.75) -- (4, 2.75);
\draw (4, 0) -- (4, 4);
\draw (4, 4) -- (0, 4);
\draw (0, 0) -- (4, 0);
\draw (0, 0) -- (0, 4);
\draw[blue, very thick] (0, 0) -- (2, 4);
\draw[blue, very thick] (2, 4) -- (4, 0);
\node (a) at (0, -.25) {$0$};
\node (b) at (4, -.25) {$1$};
\node (c) at (-.25, 4) {$1$};
%\draw[blue, very thick] (0, 4) -- (4, 2);
%\draw[blue, very thick] (0, 0) -- (4, 2);
%\draw[blue, very thick] (2, 2) -- (3, 4);
%\draw[blue, very thick] (3, 4) -- (4, 2);
%\draw[blue, very thick] (4, 2) -- (4, 0);
%\node[blue] (a) at (0, 2.75) [dot=left:$x_2$];
 \end{tikzpicture}
\]
\caption{The Tent Map}\label{fig:tent map}
\end{figure}

By Observation \ref{observation:how-transitive-points-contained} (2), any function for which $\text{tr}\of{f}$ is dense admits a 0-transitive point. For example,
let $X = [0, 1]$, and $T : X \to X$ be the tent map (depicted in Figure~\ref{fig:tent map}), i.e., 
\[
T\of{x} = \begin{cases}
          2x & 
          \text{if $x \in [0, \Fraction 1 / 2]$;} \\
          2 - 2x & 
          \text{if $x \in [\Fraction 1 / 2, 1]$.}
          \end{cases}
\]
It is well-known $\text{tr}\of{T}$ is dense in $X$. 

We give an example of a CR-dynamical system with a $1$-transitive point that is not $0$-transitive.

\begin{figure}[h]
\[
 \begin{tikzpicture}[
       decoration = {markings,
                     mark=at position .5 with {\arrow{Stealth[length=2mm]}}},
       dot/.style = {circle, fill, inner sep=2.4pt, node contents={},
                     label=#1},
every edge/.style = {draw, postaction=decorate}
                        ]
\node[cyan] (x0) [dot=below:$x$];
\node (x1) at (-1,1) [dot=left:$T\of{x}$];
\node (y1) at (1, 1) [dot=right:$T^2\of{x}$];
\node (x2) at (-1, 2) [dot=left:$T^2\of{x}$];
\node (y2) at (1, 2) [dot=right:$T^3\of{x}$];
\node (x3) at (-1, 3) [dot=left:$T^3\of{x}$];
\node (y3) at (1, 3) [dot=right:$T^4\of{x}$];
\node (d1) at (-1, 4) {$\vdots$};
\node (d2) at (1, 4) {$\vdots$};

\draw   (x0) -- (x1);
\draw	(x0) -- (y1);   
\draw   (x1) -- (x2);
\draw   (x2) -- (x3);
\draw  (y1) -- (y2);
\draw  (y2) -- (y3);
\draw (y3) -- (1, 3.5);
\draw (x3) -- (-1, 3.5);
 \end{tikzpicture}
\]
\caption{Transitivity tree $T_G\of{x}$ in Example~\ref{ex:1-transitive-not-0-transitive}}\label{fig:1-transitive-not-0-transitive}
\end{figure}
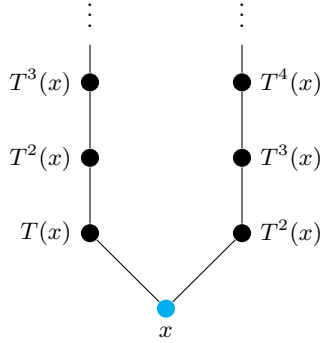

\begin{example}\label{ex:1-transitive-not-0-transitive}
Define $U := \p{\Fraction 1 / 4, \Fraction 1 / 2}$. Let $x \in \text{tr}\of{T} \setMinus U$. Define  
\[
G := \Gamma\of{T} \union \set{\p{x, T^2\of{x}}}. 
\]
Then, $\p{X, G}$ is a CR-dynamical system, and the transitivity tree $T_G\of{x}$ is depicted in Figure~\ref{fig:1-transitive-not-0-transitive}. Indeed, the transitivity tree is depicted accurately, because $x$ is a transitive point in the tent map (namely, $x$ is not periodic in $\p{X, T}$). Now, we firstly observe $x \in \transitiveIn{G}{1}$. For we have $x, T^2\of{x} \in \text{tr}\of{T}$, which implies both infinite branches $B \in \mathcal{B}_\infinity\of{T_G\of{x}}$ have $\closure{B^\ast} = X$. To derive a contradiction, suppose $x \in \transitiveIn{G}{0}$. Then, there exists $n \in \naturals$ such that $\levelIn{T_G\of{x}}{n}^\ast \subseteq U$. Since $x \notin U$, it follows $n \geq 1$. Note $\levelIn{T_G\of{x}}{n}^\ast = \set{T^n\of{x}, T^{n+1}\of{x}}$. However, $T^n\of{x} \in U$ implies $T^{n+1}\of{x} = 2 T^n\of{x}$, where $2T^n\of{x} \in \p{\Fraction 1 / 2, 1}$ (which is disjoint from $U$). This contradicts $T^{n+1}\of{x} \in U$, and thus $x \notin \transitiveIn{G}{0}$, as desired. 
\end{example}

We give two examples of CR-dynamical systems with a $0$-transitive point that has more than one trajectory.

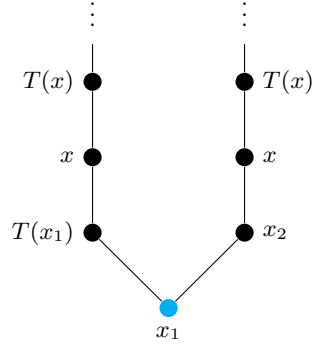
\begin{figure}[h]
\[
 \begin{tikzpicture}[
       decoration = {markings,
                     mark=at position .5 with {\arrow{Stealth[length=2mm]}}},
       dot/.style = {circle, fill, inner sep=2.4pt, node contents={},
                     label=#1},
every edge/.style = {draw, postaction=decorate}
                        ]
\node[cyan] (x0) [dot=below:$x_1$];
\node (x1) at (-1,1) [dot=left:$T\of{x_1}$];
\node (y1) at (1, 1) [dot=right:$x_2$];
\node (x2) at (-1, 2) [dot=left:$x$];
\node (y2) at (1, 2) [dot=right:$x$];
\node (x3) at (-1, 3) [dot=left:$T\of{x}$];
\node (y3) at (1, 3) [dot=right:$T\of{x}$];
\node (d1) at (-1, 4) {$\vdots$};
\node (d2) at (1, 4) {$\vdots$};

\draw   (x0) -- (x1);
\draw	(x0) -- (y1);   
\draw   (x1) -- (x2);
\draw   (x2) -- (x3);
\draw  (y1) -- (y2);
\draw  (y2) -- (y3);
\draw (y3) -- (1, 3.5);
\draw (x3) -- (-1, 3.5);
 \end{tikzpicture}
\]
\caption{Transitivity tree $T_G\of{x_1}$ in Example~\ref{ex:0-transitive-more-one-trajectory}}\label{fig:0-transitive-more-one-trajectory}
\end{figure}

\begin{example}\label{ex:0-transitive-more-one-trajectory}
Let $x \in \text{tr}\of{T}$. Set $x_1 := \Fraction x / 4$ and $x_2 := 1 - \Fraction x / 2$. Note $x_1 \neq x_2$ are both transitive points in $\p{X, T}$, and that we have 
\[
T^2\of{x_1} = x = T\of{x_2}. 
\]
Define 
\[
G := \Gamma\of{T} \union \set{\p{x_1, x_2}}. 
\]
Then, $\p{X, G}$ is a CR-dynamical system, and the transitivity tree $T_G\of{x_1}$ is depicted in Figure~\ref{fig:0-transitive-more-one-trajectory}. Indeed, $x_1$ has two trajectories,
\[
\sequence{x_1, T\of{x_1}, T^2\of{x_1}, \ldots}
\]
and
\[
\sequence{x_1, x_2, T\of{x_2}, \ldots}. 
\]
Observe $\levelIn{T_G\of{x_1}}{n}^\ast = \set{T^{n-2}\of{x}}$ for each $n \geq 2$. We claim $x_1 \in \transitiveIn{G}{0}$. To see why, suppose $U$ is a non-empty open set in $X$. Since $x \in \text{tr}\of{T}$, there exists $n \in \naturals$ such that $T^n\of{x} \in U$. Hence, $\levelIn{T_G\of{x_1}}{n+2}^\ast = \set{T^n\of{x}} \subseteq U$, which proves our claim. 
\end{example}

With Example~\ref{ex:0-transitive-more-one-trajectory}, it is clear we can modify our example by adding finitely many points to the closed relation to obtain even more trajectories. Moreover, we can change it so that the trajectories coincide at a later level, ensuring the point is $0$-transitive. Each of these types of examples have that the trajectories eventually coincide.  We construct the following example, where the trajectories do not coincide. 

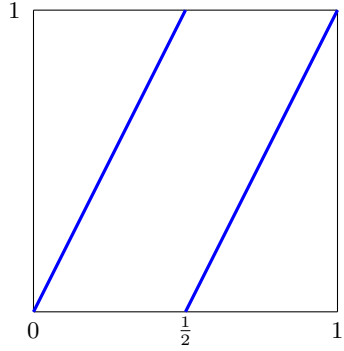
\begin{figure}[h]
\[
 \begin{tikzpicture}[
       decoration = {markings,
                     mark=at position .5 with {\arrow{Stealth[length=2mm]}}},
       dot/.style = {circle, fill, inner sep=1pt, node contents={},
                     label=#1},
every edge/.style = {draw, postaction=decorate}
                        ]
%\node (a) [dot];
%\node (b) at (5, 0) [dot];
%\node (c) at (5, 5) [dot]; 
%\node (d) at (0, 5) [dot];
%\node (e) at (8, 3) {$\color{red}\text{$G = X \times \set{0} \union \set{0} \times X$}$};
%\node (g) at (6.5, 4) {${\color{red} \text{$X = [0, 1]$}}$};
%\node (f) at (6.68, 2) {${\color{red} \text{$G\of{\Fraction 1 / 2} = \set{0}$}}$};
%\node (h) at (6.55, 1) {${\color{red} \text{$G^2\of{\Fraction 1 / 2} = X$}}$}

%\draw[dashed] (0, 2) -- (4, 2);
%\draw[dashed] (2, 0) -- (2, 4);
%\draw[dashed] (0, 2.75) -- (4, 2.75);
\draw (4, 0) -- (4, 4);
\draw (4, 4) -- (0, 4);
\draw (0, 0) -- (4, 0);
\draw (0, 0) -- (0, 4);
\draw[blue, very thick] (0, 0) -- (2, 4);
\draw[blue, very thick] (2, 0) -- (4, 4);
\node (a) at (0, -.25) {$0$};
\node (b) at (4, -.25) {$1$};
\node (c) at (-.25, 4) {$1$};
\node (d) at (2, -.25) {$\Fraction 1 /2$};
%\draw[blue, very thick] (0, 4) -- (4, 2);
%\draw[blue, very thick] (0, 0) -- (4, 2);
%\draw[blue, very thick] (2, 2) -- (3, 4);
%\draw[blue, very thick] (3, 4) -- (4, 2);
%\draw[blue, very thick] (4, 2) -- (4, 0);
%\node[blue] (a) at (0, 2.75) [dot=left:$x_2$];
 \end{tikzpicture}
\]
\caption{The Doubling Map}\label{fig:doubling map}
\end{figure}

\begin{example}\label{ex:0-transitive-interesting-example}
Let $X = [0, 1]$ and let $D : X \to 2^X$ be the doubling map (depicted in Figure~\ref{fig:doubling map}), i.e., for each $x \in X$, 
\[
D\of{x} = \begin{cases}
          \set{2x} & 
          \text{if $x \in [0, \Fraction 1 / 2)$;} \\
          \set{0, 1} &
          \text{if $x = \Fraction 1 / 2$;} \\
          \set{2x - 1} & 
          \text{if $x \in (\Fraction 1 / 2, 1]$.}
          \end{cases}
\]
Let $C = \prod_{i = 1}^\infinity \set{0, 1}$ be the Cantor space. Let $h : C \to X$ be the map representing members of $X$ in base $2$, i.e., 
\[
h\of{c_1, c_2, \ldots} = \sum_{i = 1}^\infinity \Fraction c_i / {2^i}. 
\]
Let $\sigma : C \to C$ denote the shift map, i.e., 
\[
\sigma\of{c_1, c_2, c_3, \ldots} = \p{c_2, c_3, \ldots}. 
\]
Then, for each $x \in C$, 
\[
h\of{\sigma\of{x}} = D\of{h\of{x}}, 
\]
with the exception of when $h\of{x} = \Fraction 1 / 2$. 

We now construct a binary sequence. First  write all possible blocks of length $1$, then all possible blocks of length $2$, then all possible blocks of length $3$, etc.: 
\[
\sequence{0}, \sequence{1}, \sequence{0, 0}, \sequence{0, 1}, \sequence{1, 0}, \sequence{1, 1}, \sequence{0, 0, 0}, \ldots. 
\]%sort it out here
Then we remove the brackets and let $s$ be this sequence in $C$. Now, let $s_0 \in C$ be such that it is $s$, except between each block place a $0$. Similarly, let $s_1 \in C$ be such that it is $s$, except between each block place a $1$. Observe $h\of{s_0} \neq h\of{s_1}$ are irrational, and that their trajectories in the doubling map do not coincide (their orbits are disjoint). 

Let $x_1 = \Fraction h\of{s_0} / 2$ and $x_2 = h\of{s_1}$. Define 
\[
G := \Gamma\of{D} \union \set{\p{x_1, x_2}}. 
\]
Then, $\p{X, G}$ is a CR-dynamical system, and the transitivity tree $T_G\of{x_1}$ is depicted in Figure~\ref{fig:0-transitive-interesting-example}. Notice $\levelIn{T_G\of{x_1}}{1}^\ast = \set{h\of{s_0}, h\of{s_1}}$. We claim $x_1 \in \transitiveIn{G}{0}$. To see this, let $y \in X$ and $\epsilon > 0$ be arbitrary. We show there exists $N \in \naturals$ such that $\levelIn{T_G\of{x_1}}{N+1}^\ast \subseteq B\of{y, \epsilon}$. Let 
\[
y := 0.b_1 b_2 b_3 \ldots 
\]
be the base $2$ binary expansion of $y$. Let $n \in \naturals$ such that $2^{-n} < \epsilon$. By construction of $s_0$ and $s_1$, there exists $N \in \naturals$ such that both $\sigma^N\of{s_0}$ and $\sigma^N\of{s_1}$ start with $b_1 b_2 \ldots b_n$. Hence, it follows 
\[
\levelIn{T_G\of{x_1}}{N+1}^\ast = \set{h\of{\sigma^N\of{s_0}}, h\of{\sigma^N\of{s_1}}} \subseteq B\of{y, \epsilon}, 
\]
and thus $x_1 \in \transitiveIn{G}{0}$ as desired. 
\end{example}

\begin{figure}[h]
\[
 \begin{tikzpicture}[
       decoration = {markings,
                     mark=at position .5 with {\arrow{Stealth[length=2mm]}}},
       dot/.style = {circle, fill, inner sep=2.4pt, node contents={},
                     label=#1},
every edge/.style = {draw, postaction=decorate}
                        ]
\node[cyan] (x0) [dot=below:$x_1$];
\node (x1) at (-1,1) [dot=left:$h\of{s_0}$];
\node (y1) at (1, 1) [dot=right:$h\of{s_1}$];
\node (x2) at (-1, 2) [dot=left:$h\of{\sigma\of{s_0}}$];
\node (y2) at (1, 2) [dot=right:$h\of{\sigma\of{s_1}}$];
\node (x3) at (-1, 3) [dot=left:$h\of{\sigma^2\of{s_0}}$];
\node (y3) at (1, 3) [dot=right:$h\of{\sigma^2\of{s_1}}$];
\node (d1) at (-1, 4) {$\vdots$};
\node (d2) at (1, 4) {$\vdots$};

\draw   (x0) -- (x1);
\draw	(x0) -- (y1);   
\draw   (x1) -- (x2);
\draw   (x2) -- (x3);
\draw  (y1) -- (y2);
\draw  (y2) -- (y3);
\draw (y3) -- (1, 3.5);
\draw (x3) -- (-1, 3.5);
 \end{tikzpicture}
\]
\caption{Transitivity tree $T_G\of{x_1}$ in Example~\ref{ex:0-transitive-interesting-example}}\label{fig:0-transitive-interesting-example}
\end{figure}

We pose the following questions.

\begin{question}
Does there exist a CR-dynamical system $\p{X, G}$, such that there is $x \in \transitiveIn{G}{0}$ with $\cardinality{\mathcal{B}_\infinity\of{T_G\of{x}}} = \aleph_0$?
\end{question}

\begin{question}
Does there exist a CR-dynamical system $\p{X, G}$, such that there is $x \in \transitiveIn{G}{0}$ with $\cardinality{\mathcal{B}_\infinity\of{T_G\of{x}}} > \aleph_0$? 
\end{question}

\begin{question}
Does there exist a CR-dynamical system $\p{X, G}$, such that there is $x \in \transitiveIn{G}{0}$ with $\cardinality{\levelIn{T_G\of{x}}{n}} \geq \continuum$ for some $n \in \naturals$?
\end{question}

We now review results from~\cite[Section $3$]{transitivity_CR} on transitive points in CR-dynamical systems,  Theorem \ref{thm:legally-transitive-point-motivation} to Proposition \ref{prop:no-isolated-implies-transitive-dense},
 and include similar results for $0$-transitive points when possible.   

We start with Theorem $3.25$ in~\cite{transitivity_CR}, which states that if there are transitive points in a CR-dynamical system $\p{X, G}$, then $\p{X, G}$ is an SV-dynamical system. Note $\p{X, G}$ is an SV-dynamical system if, and only if, $\illegal{G} = \emptySet$. 

%A notion of transitive points for CR-dynamical systems that are not SV-dynamical systems is desirable, which is the focus of Section~\ref{section:legally-transitive-point}. 

\begin{theorem}\label{thm:legally-transitive-point-motivation}
Let $\p{X, G}$ be a CR-dynamical system, and let $k \in \set{0, 1, 2, 3}$. If $\transitiveIn{G}{k} \neq \emptySet$, then $\illegal{G} = \emptySet$. 
\end{theorem}
\begin{proof}
Suppose $\transitiveIn{G}{k} \neq \emptySet$. By Observation~\ref{observation:how-transitive-points-contained}, $\transitiveIn{G}{k} \subseteq \transitiveIn{G}{3}$, implying $\transitiveIn{G}{3} \neq \emptySet$. The result now follows from~\cite[Theorem $3.25$]{transitivity_CR}. 
\end{proof}

%We now generalise results on transitive points from~\cite[Section $3$]{transitivity_CR} to legally transitive points. 

\begin{lemma}[{Lemma $3.27$~\cite{transitivity_CR}}]\label{lem:dense}
Let $X$ be a compact Hausdorff space, $A \subseteq X$ be such that $A$ is not dense in $X$, and $x \in X$. If $A \union \set{x}$ is dense in $X$, then $x \in \isolated{X}$. 
\end{lemma}

\begin{proposition}\label{prop:path-isolated-trans2}
Let $\p{X, G}$ be a CR-dynamical system, $x \in \isolated{X}$ and $y \in X$. Then,
\begin{enumerate}
    \item if $y \in \transitiveIn{G}{2}$, then there is a path from $y$ to $x$ in $G$;
    \item if $y \in \transitiveIn{G}{3}$, then there is a path from $y$ to $x$ in $G$;
    \item if $X$ is infinite and $\p{x, y} \in G$, then $y \notin \transitiveIn{G}{1}$. 
\end{enumerate}
\end{proposition}
\begin{proof}
Both $\p{1}$ and $\p{3}$ are proved in~\cite[Theorem $3.34$]{transitivity_CR}, so we only prove $\p{2}$. To this end, suppose $y \in \transitiveIn{G}{3}$. Then, $\Union_{B \in \mathcal{B}_\infinity\of{T_G\of{y}}} B^\ast$ is dense in $X$. In particular, 
\[
x \in \isolated{X} \subseteq \Union_{B \in \mathcal{B}_\infinity\of{T_G\of{y}}} B^\ast.
\]
Hence, it clearly follows there is a path from $y$ to $x$ in $G$. 
\end{proof}

\begin{proposition}\label{proposition:isolated-non-empty-then-transitive-isolated-points-exist}
Let $\p{X, G}$ be a CR-dynamical system such that $\isolated{X}$ is non-empty. Then the following hold. 
\begin{enumerate}
    \item If $\transitiveIn{G}{0} \neq \emptySet$, then $\transitiveIn{G}{0} \intersect \isolated{X} \neq \emptySet$. Furthermore, if $X$ is infinite and $\transitiveIn{G}{0} \neq \emptySet$, then $\transitiveIn{G}{0}$ contains exactly one isolated point of $X$. 
    \item If $\transitiveIn{G}{1} \neq \emptySet$, then $\transitiveIn{G}{1} \intersect \isolated{X} \neq \emptySet$. Furthermore, if $X$ is infinite and $\transitiveIn{G}{1} \neq \emptySet$, then $\transitiveIn{G}{1}$ contains exactly one isolated point of $X$. 
    \item If $\transitiveIn{G}{2} \neq \emptySet$, then $\transitiveIn{G}{2} \intersect \isolated{X} \neq \emptySet$. 
\end{enumerate}
\end{proposition}
\begin{proof}
Both $\p{2}$ and $\p{3}$ are proved in~\cite[Theorem $3.31$]{transitivity_CR}, so we only prove $\p{1}$. To this end, suppose $\transitiveIn{G}{0}$ is non-empty. We show $\transitiveIn{G}{0}$ contains an isolated point of $X$. Let $x \in \transitiveIn{G}{0}$. Since $x$ is $0$-transitive, there is $k \in \naturals$ such that $\levelIn{T_G\of{x}}{k}^\ast \subseteq \isolated{X}$. Therefore, there exists minimal $n \in \naturals$ such that $\levelIn{T_G\of{x}}{n}^\ast \intersect \isolated{X} \neq \emptySet$. Let $\gamma = y_0 \ldots y_n$ be a path from $x$ to $y$ in $G$, such that $y \in \isolated{X}$. Suppose $U$ is a non-empty open set in $X$. By choice of $n$, $U \setMinus \set{y_0, \ldots, y_{n-1}}$ is non-empty and open. Then, there exists $m \in \naturals$ such that $\levelIn{T_G\of{x}}{m}^\ast \subseteq U \setMinus \set{y_0, \ldots, y_{n-1}}$. Notice $m \geq n$. Hence, $\levelIn{T_G\of{y}}{m - n}^\ast \subseteq U$. Thus, $y \in \transitiveIn{G}{0}$, implying $\transitiveIn{G}{0} \intersect \isolated{X} \neq \emptySet$.

Now, suppose $X$ is infinite and $\transitiveIn{G}{0} \neq \emptySet$. By the above, $\transitiveIn{G}{0} \intersect \isolated{X} \neq \emptySet$. By Observation~\ref{observation:how-transitive-points-contained}, it follows $\transitiveIn{G}{1}$ is non-empty because $\transitiveIn{G}{0} \neq \emptySet$.  By $\p{2}$, $\transitiveIn{G}{1}$ contains exactly one isolated point. Since $\transitiveIn{G}{0} \subseteq \transitiveIn{G}{1}$, it follows $\transitiveIn{G}{0}$ contains exactly one isolated point, and we are done. 
\end{proof}

%[Theorem $3.35$~\cite{transitivity_CR}\footnote{There is a typo in $\p{1}$ of~\cite[Theorem $3.35$]{transitivity_CR}. It is meant to say that if $x \in \intrans{G}$, then for each $y \in \legal{G}$, $\p{x, y} \in G$ implies $y \in \intrans{G}$. A counter-example to how it is written is as follows. Let $X = \set{0, 1}$ and $G = \set{\p{0, 0}, \p{0, 1}}$. Then, $0$ is an intransitive point of $\p{X, G}$, but $1 \in \illegal{G}$ and so cannot be intransitive.}]

\begin{proposition}\label{proposition:transitive-branches}
Let $\p{X, G}$ be a CR-dynamical system and let $x \in X$. Then the following hold. 
\begin{enumerate}
    \item If $x \in \intrans{G}$, then for each $y \in \legal{G}$,
    \[
    \p{x, y} \in G \Longrightarrow y \in \intrans{G}. 
    \]
    \item If $x \in \intrans{G}$, then for each $B \in \mathcal{B}_\infinity\of{T_G\of{x}}$, $B^\ast \subseteq \intrans{G}$.
    \item If $x \notin \isolated{X}$ and if $x \in \transitiveIn{G}{0}$, then for each $y \in X$, 
    \[
    \p{x, y} \in G \Longrightarrow y \in \transitiveIn{G}{0}.
    \]
    \item If $\isolated{X} = \emptySet$ and if $x \in \transitiveIn{G}{0}$, then for each $B \in \mathcal{B}_\infinity\of{T_G\of{x}}$, $B^\ast \subseteq \transitiveIn{G}{0}$. 
    \item If $x \notin \isolated{X}$ and if $x \in \transitiveIn{G}{1}$, then for each $y \in X$,
    \[
    \p{x, y} \in G \Longrightarrow y \in \transitiveIn{G}{1}. 
    \]
    \item If $\isolated{X} = \emptySet$ and if $x \in \transitiveIn{G}{1}$, then for each $B \in \mathcal{B}_\infinity\of{T_G\of{x}}$, $B^\ast \subseteq \transitiveIn{G}{1}$. 
    \item If $x \notin \isolated{X}$ and if $x \in \transitiveIn{G}{2}$, then there exists $y \in X$ such that
    \[
    \p{x, y} \in G \text{ and } y \in \transitiveIn{G}{2}. 
    \]
    \item If $\isolated{X} = \emptySet$ and if $x \in \transitiveIn{G}{2}$, then there is $B \in \mathcal{B}_\infinity\of{T_G\of{x}}$ such that $\closure{B^\ast} = X$ and $B^\ast \subseteq \transitiveIn{G}{2}$. 
\end{enumerate}
\end{proposition}
\begin{proof}
It suffices to prove $\p{3}$ and $\p{4}$. The rest can be found in Theorem $3.35$~\cite{transitivity_CR}. To this end, suppose $x \notin \isolated{X}$ and $x \in \transitiveIn{G}{0}$. Further suppose $y \in X$ and $\p{x, y} \in G$. Let $U$ be a non-empty open set in $X$. Since $x$ is not isolated in $X$, $U \setMinus \set{x}$ is open and non-empty in $X$. It follows there exists $n \in \naturals$ such that $\levelIn{T_G\of{x}}{n}^\ast \subseteq U \setMinus \set{x}$. Notice $n \geq 1$. Hence, $\levelIn{T_G\of{y}}{n - 1}^\ast \subseteq U$. It follows $y \in \transitiveIn{G}{0}$, yielding $\p{3}$.  We note $\p{4}$ follows by inductively applying $\p{3}$. 
\end{proof}

\begin{proposition}\label{prop:no-isolated-implies-transitive-dense}
Let $\p{X, G}$ be a CR-dynamical system and $k \in \set{0, 1, 2}$. If we have $\isolated{X} = \emptySet$ and $\transitiveIn{G}{k} \neq \emptySet$, then $\transitiveIn{G}{k}$ is dense in $X$. 
\end{proposition}
\begin{proof}
Both $k = 1$ and $k = 2$ are proved in Theorem $3.38$~\cite{transitivity_CR}, so we prove for $k = 0$. Suppose $\isolated{X} = \emptySet$ and $\transitiveIn{G}{0} \neq \emptySet$. By Proposition~\ref{proposition:transitive-branches}, there exists $B \in \mathcal{B}_\infinity\of{T_G\of{x}}$ with $\closure{B^\ast} = X$ such that $B^\ast \subseteq \transitiveIn{G}{0}$. Thus, $\transitiveIn{G}{0}$ is dense in $X$. 
\end{proof}

The above result, adapted from~\cite[Theorem $3.38$]{transitivity_CR}, does not hold for $\transitiveIn{G}{3}$ (see~\cite[Example $3.37$]{transitivity_CR}). We will give a sufficient condition for $\transitiveIn{G}{3}$ to be dense. Firstly, we require the following observation. 

\begin{proposition}\label{prop:backwards}
Suppose $\p{X, G}$ is a CR-dynamical system and $k \in \set{2, 3}$. If $y \in \transitiveIn{G}{k}$ and $\p{x, y} \in G$, then $x \in \transitiveIn{G}{k}$, and consequently $T_{G^{-1}}\of{y}^\ast \subseteq \transitiveIn{G}{k}$.
\end{proposition}
\begin{proof}
Suppose $y \in \transitiveIn{G}{2}$ and $\p{x, y} \in G$. Then, there exists $B_y \in \mathcal{B}_\infinity\of{T_G\of{y}}$ such that $B_y^\ast$ is dense in $X$. By Observation~\ref{observation:branch-extension}, there exists $B_x \in \mathcal{B}_\infinity\of{T_G\of{x}}$ such that $B_y^\ast \subseteq B_x^\ast$. It follows $B_x^{\ast}$ is dense in $X$. Hence, $x \in \transitiveIn{G}{2}$. 

Suppose $y \in \transitiveIn{G}{3}$ and $\p{x, y} \in G$. Then, $\Union_{B \in \mathcal{B}_\infinity\of{T_G\of{y}}} B^\ast$ is dense in $X$.  By Observation~\ref{observation:branch-extension}, 
\[
\Union_{B_y \in \mathcal{B}_\infinity\of{T_G\of{y}}} B_y^\ast \subseteq \Union_{B_x \in \mathcal{B}_\infinity\of{T_G\of{x}}} B_x^\ast.
\]
It follows $\Union_{B \in \mathcal{B}_\infinity\of{T_G\of{x}}} B^\ast$ is dense in $X$. Thus, $x \in \transitiveIn{G}{3}$. 
\end{proof}

\begin{proposition}\label{prop:sufficient-for-ltrans3-dense}
Let $\p{X, G}$ be a CR-dynamical system and $k \in \set{2, 3}$. If $\transitiveIn{G}{k} \intersect \transitiveIn{G^{-1}}{k}$ is non-empty, then 
\begin{enumerate}
    \item $\transitiveIn{G}{k}$ is dense in $X$;
    \item $\transitiveIn{G^{-1}}{k}$ is dense in $X$; 
    \item $\isolated{X} \subseteq \transitiveIn{G}{k} \intersect \transitiveIn{G^{-1}}{k}$. 
\end{enumerate} 
\end{proposition}
\begin{proof}
Suppose there exists $x \in \transitiveIn{G}{k} \intersect \transitiveIn{G^{-1}}{k}$. Then, by Proposition~\ref{prop:backwards}, $T_G\of{x}^\ast \subseteq \transitiveIn{G^{-1}}{k}$, and $T_{G^{-1}}\of{x}^\ast \subseteq \transitiveIn{G}{k}$. It follows 
\[
\closure{\transitiveIn{G^{-1}}{k}} = X
\]
and
\[
\closure{\transitiveIn{G}{k}} = X. 
\]
From this, $\p{1}$, $\p{2}$, and $\p{3}$ clearly follow.
\end{proof}

\begin{example}
We show Proposition~\ref{prop:sufficient-for-ltrans3-dense} does not hold in general for $k \in \set{0, 1}$. Let $X = \set{1, 2}$ and $G = \set{\p{1, 2}, \p{2, 1}, \p{2, 2}}$. Then, $G = G^{-1}$, and 
\[
\transitiveIn{G}{0} = \transitiveIn{G}{1} = \set{1}, 
\]
which is not dense in $X$. Indeed, this example also shows Proposition~\ref{prop:backwards} fails for $k \in \set{0, 1}$ as well. 
\end{example}

We conclude this section by generalising the following well-known result to SV-dynamical systems. %Indeed, we find the result generalises best to $0$-transitive points, with the result also generalising (albiet less trivially) to $1$-transitive points, and a partial generalisation to $3$-transitive points.

\begin{theorem}\label{thm:transitive-points-dense-Gdelta}
Let $\p{X, f}$ be a topological dynamical system, where $\collection{U}$ is a countable base for $X$. Then, 
\[
\emph{tr}\of{f} = \Intersection_{U \in \collection{U}} \p{\Union_{k = 0}^\infinity f^{-k}\of{U}}. 
\]
In particular, $\emph{tr}\of{f}$ is a $G_\delta$ set. 
\end{theorem}

\begin{theorem}\label{thm:0-transitive-points-dense-Gdelta}
Let $\p{X, G}$ be an SV-dynamical system, where $\collection{U}$ is a countable base for $X$. Then,
\[
\transitiveIn{G}{0} = \Intersection_{U \in \collection{U}} \p{\Union_{k = 0}^\infinity G^{-k}[U]}. 
\]
In particular, $\transitiveIn{G}{0}$ is a $G_\delta$ set. 
\end{theorem}
\begin{proof}
Suppose $x \in \transitiveIn{G}{0}$ and $U \in \collection{U}$. Then, there exists $n \in \naturals$ such that $\levelIn{T_G\of{x}}{n}^\ast \subseteq U$. It follows $G^n\of{x} \subseteq U$, which implies $x \in G^{-n}[U]$.  Thus, 
$
\transitiveIn{G}{0} \subseteq \Intersection_{U \in \collection{U}} \p{\Union_{k = 0}^\infinity G^{-k}[U]}. 
$

Now, suppose $x \in \Intersection_{U \in \collection{U}} \p{\Union_{k = 0}^\infinity G^{-k}[U]}$. Suppose $O$ is a non-empty open set in $X$. There exists $U \in \collection{U}$ such that $U \subseteq O$. There exists $n \in \naturals$ such that $x \in G^{-n}[U]$. That is to say, 
$
\levelIn{T_G\of{x}}{n}^\ast = G^n\of{x} \subseteq U \subseteq O. 
$
Hence, $x \in \transitiveIn{G}{0}$,  and so
$ \Intersection_{U \in \collection{U}} \p{\Union_{k = 0}^\infinity G^{-k}[U]}\subseteq\transitiveIn{G}{0}$ as required.

Now, we show $\transitiveIn{G}{0}$ is a $G_\delta$ set in $X$. Firstly, if $\transitiveIn{G}{0} = \emptySet$, we are done. Consequently, we assume $\transitiveIn{G}{0} \neq \emptySet$. Since $\p{X, G}$ is an SV-dynamical system, $G$ is the graph of an usc set-valued function on $X$. It follows $G^{-k}[U]$ is open for each $k \in \naturals$ and $U \in \collection{U}$. Therefore, $\Union_{k = 0}^\infinity G^{-k}[U]$ is open for each $U \in \collection{U}$. Since $\collection{U}$ is countable, we obtain $\transitiveIn{G}{0}$ is $G_\delta$ (as desired).  
\end{proof}

\begin{proposition}\label{prop:ltrans1-Gdelta}
Let $\p{X, G}$ be an SV-dynamical system, where $\collection{U}$ is a countable base for $X$. Then,
\[
\transitiveIn{G}{1} = \Intersection_{U \in \collection{U}} \setOf{x \in X}:{\forall B \in \mathcal{B}_\infinity\of{T_G\of{x}}, B^\ast \intersect U \neq \emptySet}. 
\]
In particular, $\transitiveIn{G}{1}$ is a $G_\delta$ set in $X$. 
\end{proposition}
\begin{proof}
Suppose $x \in \transitiveIn{G}{1}$ and $U \in \collection{U}$. Then, $B^\ast \intersect U \neq \emptySet$ for each $B \in \mathcal{B}_\infinity\of{T_G\of{x}}$, since $\closure{B^\ast} = X$ for each $B \in \mathcal{B}_\infinity\of{T_G\of{x}}$, yielding the $\p{\subset}$ inclusion. 

On the other hand, suppose 
\[
x \in \Intersection_{U \in \collection{U}} \setOf{y \in X}:{\forall B \in \mathcal{B}_\infinity\of{T_G\of{y}}, B^\ast \intersect U \neq \emptySet}. 
\]
Let $B \in \mathcal{B}_\infinity\of{T_G\of{x}}$. Suppose $O$ is a non-empty open set in $X$. Then, there exists $U \in \collection{U}$ such that $U \subseteq O$. It follows $B^\ast \intersect U \neq \emptySet$, implying $B^\ast \intersect O \neq \emptySet$. Hence, $B^\ast$ is dense in $X$, and thus it follows $x \in \transitiveIn{G}{1}$.

Now, we show $\transitiveIn{G}{1}$ is a $G_\delta$ set in $X$. Notice we need only show 
\[
\setOf{y \in X}:{\forall B \in \mathcal{B}_\infinity\of{T_G\of{y}}, B^\ast \intersect U \neq \emptySet}
\]
is open in $X$ for each $U \in \collection{U}$. Fix $U \in \collection{U}$.  It suffices to show 
\[
C = \setOf{y \in X}:{\exists B \in \mathcal{B}_\infinity\of{T_G\of{y}}, B^\ast \subseteq X \setMinus U} 
\]
is sequentially closed. To this end, suppose $\sequenceOf{x_n}:{n \in \naturals}$ is a sequence in $C$ converging to $x \in X$. For each $n \in \naturals$, there exists $B_n \in \mathcal{B}_\infinity\of{T_G\of{x_n}}$ such that $B_n^\ast \subseteq X \setMinus U$. Each infinite branch $B_n$ corresponds to $\bm{x}_n \in \mahavier{\infinity}{G}$. As $\mahavier{\infinity}{G}$ is compact, we may assume (passing to subsequences if necessary) that $\bm{x}_n$ converges to some $\bm{x} \in \mahavier{\infinity}{G}$. Indeed, $\bm{x}$ corresponds to an infinite branch $B \in \mathcal{B}_\infinity\of{T_G\of{x}}$, and since we have $B_n^\ast \subseteq X \setMinus U$ it follows by continuity of the projection maps that $B^\ast \subseteq X \setMinus U$ (since $X \setMinus U$ is closed). Thus, $x \in C$, and we are done.
\end{proof}

\begin{proposition}\label{prop:3trans-characterisation}
Let $\p{X, G}$ be an SV-dynamical system, where $\collection{U}$ is a countable base for $X$. Then, 
\[
\transitiveIn{G}{3} = \Intersection_{U \in \collection{U}} \p{\Union_{k = 0}^\infinity G^{-k}\of{U}} = \Intersection_{U \in \collection{U}} \p{\Union_{k = 0}^\infinity G^{-k}\of{\closure{U}}} = \Intersection_{U \in \collection{U}} \p{\Union_{x \in U} T_{G^{-1}}\of{x}^\ast}. 
\]
In particular, $\transitiveIn{G}{3}$ is the countable intersection of $F_\sigma$ sets. 
\end{proposition}
\begin{proof}
Firstly, it is easy to see that
\[
\Intersection_{U \in \collection{U}} \p{\Union_{k = 0}^\infinity G^{-k}\of{U}} = \Intersection_{U \in \collection{U}} \p{\Union_{k = 0}^\infinity G^{-k}\of{\closure{U}}} = \Intersection_{U \in \collection{U}} \p{\Union_{x \in U} T_{G^{-1}}\of{x}^\ast}, 
\]
so we only show $\transitiveIn{G}{3} = \Intersection_{U \in \collection{U}} \p{\Union_{k = 0}^\infinity G^{-k}\of{U}}$. 

Furthermore, $\Union_{k = 0}^\infinity G^{-k}\of{\closure{U}}$ is clearly an $F_\sigma$ set for each $U \in \collection{U}$. Therefore, it follows $\transitiveIn{G}{3}$ is the countable intersection of $F_\sigma$ sets. 

Suppose $x \in \transitiveIn{G}{3}$ and $U \in \collection{U}$. Then, $\Union_{B \in \mathcal{B}_\infinity\of{T_G\of{x}}} B^\ast$ is dense in $X$, implying $G^n\of{x}$ meets $U$ for some $n \in \naturals$. That is to say, $x \in G^{-n}\of{U}$, which establishes the $\p{\subset}$ inclusion. 

On the other hand, suppose $x \in \Intersection_{U \in \collection{U}} \p{\Union_{k = 0}^\infinity G^{-k}\of{U}}$. Let $O$ be a non-empty open set in $X$. There exists $U \in \collection{U}$ such that $U \subseteq O$. There exists $n \in \naturals$ such that $U$ meets $G^n\of{x}$. Therefore, $O \intersect G^n\of{x} \neq \emptySet$. It clearly follows that $O$ meets $\Union_{B \in \mathcal{B}_\infinity\of{T_G\of{x}}} B^\ast$. Thus, $x \in \transitiveIn{G}{3}$, which establishes the $\p{\superset}$ inclusion. 
\end{proof}

%---------------------------------------------------------------
\section{$2$-transitive and $3$-transitive points}\label{section:3-trans-not-2-trans}

In this section we consider  $2$-transitive and $3$-transitive points.  To date, there are no examples of $3$-transitive points that are not $2$-transitive, unless $\transitiveIn{G}{2} = \emptySet$. This observation motivates us to ask if $\transitiveIn{G}{2} \neq \emptySet$ implies $\transitiveIn{G}{2} = \transitiveIn{G}{3}$.  We show that the implication does not hold (see Example~\ref{ex:stormy-campsite-map}).  We start by showing that when there are $2$-transitive points and isolated points, $\transitiveIn{G}{2} = \transitiveIn{G}{3}$.

\begin{theorem}\label{thm:fundamental-isolated}
Let $\p{X, G}$ be a CR-dynamical system such that $\isolated{X} \neq \emptySet$. If $\transitiveIn{G}{2} \neq \emptySet$, then there exists $x \in \isolated{X} \intersect \transitiveIn{G}{2}$ such that 
\[
T_{G^{-1}}\of{x}^\ast = \transitiveIn{G}{2} = \transitiveIn{G}{3}. 
\]
\end{theorem}
\begin{proof}
Suppose $\transitiveIn{G}{2} \neq \emptySet$. Since $\isolated{X}$ is non-empty, there exists $x \in \isolated{X} \intersect \transitiveIn{G}{2}$ by Proposition~\ref{proposition:isolated-non-empty-then-transitive-isolated-points-exist}. We know 
\[
T_{G^{-1}}\of{x}^\ast \subseteq \transitiveIn{G}{2} \subseteq \transitiveIn{G}{3}
\]
by Proposition~\ref{prop:backwards} and Observation~\ref{observation:3.19-3.20-3.21}. On the other hand, if $y \in \transitiveIn{G}{3}$, there is a path from $y$ to $x$ in $G$ by Proposition~\ref{prop:path-isolated-trans2}. Therefore, there is a path from $x$ to $y$ in $G^{-1}$, which implies $y \in T_{G^{-1}}\of{x}^\ast$. Thus, $T_{G^{-1}}\of{x}^\ast = \transitiveIn{G}{2} = \transitiveIn{G}{3}$. 
\end{proof}

By Theorem~\ref{thm:fundamental-isolated}, we obtain the following. 

\begin{corollary}
Let $\p{X, G}$ be a CR-dynamical system such that $\transitiveIn{G}{2} \neq \emptySet$ and $\isolated{X} \neq \emptySet$. Then, $\transitiveIn{G}{2} = \transitiveIn{G}{3}$. 
\end{corollary}

Of course, when studying $3$-transitive points that are not $2$-transitive, it is natural to assume $\transitiveIn{G}{2} \neq \emptySet$. For if $\transitiveIn{G}{2} = \emptySet$, every $3$-transitive point of $\p{X, G}$ is not $2$-transitive. Hence, by the above result, we may assume there are no isolated points. As we also assume $\transitiveIn{G}{2} \neq \emptySet$, Proposition~\ref{prop:no-isolated-implies-transitive-dense} tells us $\transitiveIn{G}{2}$ and $\transitiveIn{G}{3}$ are dense in $X$.  

To construct our example of a CR-dynamical system with $\transitiveIn{G}{2} \neq \emptySet$ and $\transitiveIn{G}{2} \neq \transitiveIn{G}{3}$, we make a few observations on the equivalence relation $\sim_G$ (see Definition~\ref{def:transitivity-tree-relation} in Section~\ref{section:transitivity-trees}).

\begin{observation}
Let $\p{X, G}$ be a CR-dynamical system.  If $x \in \intrans{G}$, then $[x]_G \subseteq \intrans{G}$.
\end{observation}

\begin{observation}\label{observation:if-2/3-trans-so-is-equiv-class}
Let $\p{X, G}$ be a CR-dynamical system and $k \in \set{2, 3}$. If $x \in \transitiveIn{G}{k}$, then $[x]_G \subseteq \transitiveIn{G}{k}$. 
\end{observation}
\begin{note}
Follows directly from Proposition~\ref{prop:backwards}. 
\end{note}

\begin{observation}
Let $\p{X, G}$ be a CR-dynamical system such that $\isolated{X} = \emptySet$, and $k \in \set{0, 1, 2, 3}$. If $x \in \transitiveIn{G}{k}$, then $[x]_G \subseteq \transitiveIn{G}{k}$. 
\end{observation}
\begin{note}
Follows directly from Proposition~\ref{prop:backwards} and Proposition~\ref{proposition:transitive-branches}. 
\end{note}

\begin{proposition}\label{prop:equiv-class-not-dense-if-not-2-trans}
Let $\p{X, G}$ be a CR-dynamical system. Then, $[x]_G$ is not dense in $X$ for each $x \in X \setMinus \transitiveIn{G}{2}$.
\end{proposition}
\begin{proof}
Suppose $[x]_G$ is dense in $X$ for some $x \in X \setMinus \transitiveIn{G}{2}$. Let $\setOf{U_n}:{n \in \naturals}$ be a base for $X$. For each $n \in \naturals$, let $x_n \in [x]_G \intersect U_n$. Then, there is a path from $x$ to $x_0$ in $G$, and a path from $x_n$ to $x_{n+1}$ in $G$ for each $n \in \naturals$. But this means there is a dense infinite branch of $T_{G}\of{x}$, contradicting the fact $x \notin \transitiveIn{G}{2}$. 
\end{proof}

\begin{corollary}\label{cor:[x]G-dense-imply-contained-2trans}
Let $\p{X, G}$ be a CR-dynamical system, and let $x \in X$. If $[x]_G$ is dense in $X$, then $[x]_G \subseteq \transitiveIn{G}{2}$. 
\end{corollary}
\begin{proof}
Suppose $[x]_G$ is dense in $X$. By our proof in Proposition~\ref{prop:equiv-class-not-dense-if-not-2-trans}, it follows $x \in \transitiveIn{G}{2}$. By Observation~\ref{observation:if-2/3-trans-so-is-equiv-class}, $[x]_G \subseteq \transitiveIn{G}{2}$. 
\end{proof}

\begin{figure}[h]
\[
 \begin{tikzpicture}[
       decoration = {markings,
                     mark=at position .5 with {\arrow{Stealth[length=2mm]}}},
       dot/.style = {circle, fill, inner sep=1pt, node contents={},
                     label=#1},
every edge/.style = {draw, postaction=decorate}
                        ]
%\node (a) [dot];
%\node (b) at (5, 0) [dot];
%\node (c) at (5, 5) [dot]; 
%\node (d) at (0, 5) [dot];
%\node (e) at (8, 3) {$\color{red}\text{$G = X \times \set{0} \union \set{0} \times X$}$};
%\node (g) at (6.5, 4) {${\color{red} \text{$X = [0, 1]$}}$};
%\node (f) at (6.68, 2) {${\color{red} \text{$G\of{\Fraction 1 / 2} = \set{0}$}}$};
%\node (h) at (6.55, 1) {${\color{red} \text{$G^2\of{\Fraction 1 / 2} = X$}}$}
'

\node (a) at (2, -.25) {$\color{red}\text{$G = \Graph{f_1} \union \Graph{f_2} \union \set{\p{0, x_1}, \p{0, x_2}} \union A \union \Delta_X$}$};

%\draw[dashed] (0, 2) -- (4, 2);
%\draw[dashed] (2, 0) -- (2, 4);
%\draw[dashed] (0, 2.75) -- (4, 2.75);
\draw (4, 0) -- (4, 4);
\draw (4, 4) -- (0, 4);
\draw (0, 0) -- (4, 0);
\draw (0, 0) -- (0, 4);
\draw[dashed] (2, 0) -- (2, 4);
\draw[dashed] (0, 2) -- (4, 2);
\draw[blue, very thick] (0, 0) -- (4, 4);
\draw[blue, very thick] (0, 0) -- (1, 2);
\draw[blue, very thick] (1, 2) -- (2, 0);
\draw[blue, very thick] (2, 2) -- (3, 4);
\draw[blue, very thick] (3, 4) -- (4, 2); 
%\draw[blue, very thick] (2, 2) -- (3, 4);
%\draw[blue, very thick] (3, 4) -- (4, 2);
%\draw[blue, very thick] (4, 2) -- (4, 0);
\node[blue] (a) at (0, 0.85) [dot=left:$x_1$];
\node[blue] (a) at (0, 2.85) [dot=left:$x_2$];
\node[blue] (a) at (2, 3) [dot];
\node[blue] (a) at (2, 1) [dot];
\node[blue] (a) at (1, 1.5) [dot];
\node[blue] (a) at (1, 0.5) [dot];
\node[blue] (a) at (0.5, 0.25) [dot];
\node[blue] (a) at (0.5, 0.75) [dot];
\node[blue] (a) at (1.5, 1.25) [dot];
\node[blue] (a) at (1.5, 1.75) [dot];
\node[blue] (a) at (3, 3.5) [dot];
\node[blue] (a) at (3, 2.5) [dot];
\node[blue] (a) at (3.5, 3.25) [dot];
\node[blue] (a) at (3.5, 3.75) [dot];
\node[blue] (a) at (2.5, 2.25) [dot];
\node[blue] (a) at (2.5, 2.75) [dot];

\node[blue] (a) at (2.25, 2.125) [dot];
\node[blue] (a) at (2.25, 2.375) [dot];
\node[blue] (a) at (2.75, 2.875) [dot];
\node[blue] (a) at (2.75, 2.625) [dot];

\node[blue] (a) at (3.25, 3.125) [dot];
\node[blue] (a) at (3.25, 3.375) [dot];
\node[blue] (a) at (3.75, 3.875) [dot];
\node[blue] (a) at (3.75, 3.625) [dot];

\node[blue] (a) at (0.25, 0.375) [dot];
\node[blue] (a) at (0.25, 0.125) [dot];
\node[blue] (a) at (0.75, 0.875) [dot];
\node[blue] (a) at (0.75, 0.625) [dot];
\node[blue] (a) at (1.25, 1.375) [dot];
\node[blue] (a) at (1.25, 1.125) [dot];
\node[blue] (a) at (1.75, 1.875) [dot];
\node[blue] (a) at (1.75, 1.625) [dot];
 \end{tikzpicture}
\]
\caption{The relation $G$ from Example~\ref{ex:stormy-campsite-map}}\label{fig:Stormy-Campsite-Map}
\end{figure}

We now provide an example, showing there is a CR-dynamical system with $\isolated{X} = \emptySet$, $\transitiveIn{G}{2} \neq \emptySet$ and $\transitiveIn{G}{2} \neq \transitiveIn{G}{3}$. 

\begin{example}\label{ex:stormy-campsite-map}
Let $X = [0, 1]$ and $D = \setOf{\Fraction k / {2^n}}:{n \in \naturals, k \in [2^n], 2 \nmid k}$ be the set of dyadic rationals in $[0, 1]$. From now on, when we write $\Fraction k / {2^n}$, implicitly $2 \nmid k$. Let
\[
A = \setOf{\p{\Fraction k / {2^n}, \Fraction 2k + 1 / {2^{n+1}}}}:{\Fraction k / {2^n} \in D \setMinus \set{0, 1}} \union  \setOf{\p{\Fraction k / {2^n}, \Fraction 2k - 1 / {2^{n+1}}}}:{\Fraction k / {2^n} \in D \setMinus \set{0, 1}}. 
\]

Define $f_1 : [0, \Fraction 1 / 2] \to [0, \Fraction 1 / 2]$ by
\[
f_1\of{t} = \begin{cases}
            2t &
            \text{if $t \in [0, \Fraction 1 / 4]$;} \\
            1 - 2t &
            \text{if $t \in [\Fraction 1 / 4, \Fraction 1 / 2]$.}
            \end{cases}
\]
Define $f_2 : [\Fraction 1 / 2, 1] \to [\Fraction 1 / 2, 1]$ by 
\[
f_2\of{t} = \begin{cases}
            2t - \Fraction 1 / 2 &
            \text{if $t \in [\Fraction 1 / 2, \Fraction 3 / 4]$;} \\
            \Fraction 5 / 2 - 2t &
            \text{if $t \in [\Fraction 3 / 4, 1]$.}
            \end{cases}
\]
Let $x_1 \in \text{tr}\of{f_1}$ and $x_2 \in \text{tr}\of{f_2}$. 

Let $G$ be the closed relation
\[
G = \Graph{f_1} \union \Graph{f_2} \union \set{\p{0, x_1}, \p{0, x_2}} \union A \union \Delta_X
\]
approximated in Figure~\ref{fig:Stormy-Campsite-Map}.  

Now, $G\of{0} = \set{0, x_1, x_2}$. Since $x_1 \in \text{tr}\of{f_1}$, $x_1$ and its iterates under $f_1$ are irrational. Similarly, $x_2 \in \text{tr}\of{f_2}$, $x_2$ and its iterates under $f_2$ are irrational. Therefore, the orbit structure of $0$ is determined by $G \setMinus A$. Hence, by similar argument in~\cite[Example $3.40$]{transitivity_CR}, it follows $0 \in \transitiveIn{G}{3} \setMinus \transitiveIn{G}{2}$. %In fact, $0 \in \transitiveIn{G}{\p{3, \omega, 2}}$. 

Now, we claim $\Fraction 1 / 2 \in \transitiveIn{G}{2}$. Clearly, $D \setMinus \set{0, 1} \subseteq T_G\of{\Fraction 1 / 2}^\ast$, where  $D \setMinus \set{0, 1}$ is dense in $X$. It follows $\Fraction 1 / 2 \in \transitiveIn{G}{3}$. If we can show all but finitely many dyadic rationals have a path to $\Fraction 1 / 2$ in $G$, it will follow $\Fraction 1 / 2 \in \transitiveIn{G}{2}$. For then $[\Fraction 1 / 2]_G$ will be dense in $X$, which implies $\Fraction 1 /2 \in [\Fraction 1 / 2]_G \subseteq \transitiveIn{G}{2}$. To see why, we notice all dyadic rationals in $\p{0, \Fraction 1 / 2}$ eventually go to $\Fraction 1 / 2$ under iteration of $f_1$, and all dyadic rationals in $\left[\Fraction 1 / 2, 1\right]$ eventually go to $\Fraction 1 / 2$ under iteration of $f_2$. Our claim follows, and we are done. 
\end{example} 

\begin{observation}
In Example~\ref{ex:stormy-campsite-map} we have $\transitiveIn{G}{2} = D \setMinus \set{0}$ and $\transitiveIn{G}{3} = D$. As both of these sets are not $G_\delta$, Theorem~\ref{thm:transitive-points-dense-Gdelta} does not generalise fully to $\transitiveIn{G}{k}$, for $k \in \set{2, 3}$. As we have seen, however, Theorem~\ref{thm:transitive-points-dense-Gdelta} generalises to $\transitiveIn{G}{0}$ and $\transitiveIn{G}{1}$ by Theorem~\ref{thm:0-transitive-points-dense-Gdelta} and Proposition~\ref{prop:ltrans1-Gdelta}, respectively. 
\end{observation}

For the remainder of this section, we expand upon the theory developed in~\cite[Section $4$]{transitivity_CR} to further distinguish $2$-transitive and $3$-transitive points.  We now go over the definitions provided in~\cite[Section $4$]{transitivity_CR}. With respect to~\cite[Definition $4.1$]{transitivity_CR},  we start by observing 
\[
\mathcal{L}_n\of{T_G\of{x}}^\ast = \mathcal{R}_n\of{x}
\]
and
\[
T_G\of{x}^\ast = \mathcal{R}_\omega\of{x}
\]
for each $x \in X$ and $n \in \naturals$. Furthermore, we note $\mathcal{L}_n\of{T_G\of{x}}^\ast$ is closed in $X$ for each $x \in X$ and $n \in \naturals$. Thus, we do not need to refer to their closures in our equivalent definition to~\cite[Definition $4.3$]{transitivity_CR}.

\begin{definition}[{Definition $4.3$~\cite{transitivity_CR}}]\label{definition:transitive-3n}
Let $\p{X, G}$ be a CR-dynamical system and $x \in \transitiveIn{G}{3} \setMinus \transitiveIn{G}{2}$. We say that $x$ is 
\begin{enumerate}
    \item {\color{cyan} \emph{$\p{3, n}$-transitive in $\p{X, G}$}}, if there is a positive integer $n$ such that 
    \[
    \mathcal{L}_n\of{T_G\of{x}}^\ast = X,
    \]
    and 
    \[
    \mathcal{L}_{n-1}\of{T_G\of{x}}^\ast \neq X. 
    \]
    We use {\color{cyan} $\transitiveIn{G}{\p{3, n}}$} to denote the set of $\p{3, n}$-transitive points in $\p{X, G}$. 
    \item {\color{cyan} \emph{$\p{3, \omega}$-transitive in $\p{X, G}$}}, if 
    \[
    x \notin \Union_{n=1}^\infinity \transitiveIn{G}{\p{3, n}}. 
    \]
    We use {\color{cyan} $\transitiveIn{G}{\p{3, \omega}}$} to denote the set of $\p{3, \omega}$-transitive points in $\p{X, G}$. 
\end{enumerate}
\end{definition}

\begin{example}[{Example $4.5$~\cite{transitivity_CR}}]
Let $\p{X, G}$ be a CR-dynamical system such that $X = [0, 1]$ and $G = \p{X \times \set{1}} \union \p{\set{0} \times X}$. Then, $0 \in \transitiveIn{G}{\p{3, 1}}$.  
\end{example}

\begin{example}\label{ex:3n-transitive}
Let $X = [0, 1]$ and for each $n \in \naturals$, define
\[
G_n = \setOf{\p{\Fraction 1 / {k+1}, \Fraction 1 / {k + 2}}}:{k \in [n]} \union \p{\set{\Fraction 1 / {n+2}} \times X} \union \p{X \times \set{0}}. 
\]
Then, $1 \in \transitiveIn{G_n}{\p{3, n+2}}$ for each $n \in \naturals$. 
\end{example} 

%We note Example $3.40$ in~\cite{transitivity_CR} claims $\transitiveIn{G}{3}$ is non-empty, yet not dense. We find this is incorrect, as $\transitiveIn{G}{3}$ is the set of dyadic rationals.  

%I was going through Example 3.40 again in the paper on transitive points in CR-Dynamical systems you published with Iztok (I have attached an image of the example).  Example 3.40 is meant to be an example where trans3(G) is non-empty, yet not dense.  Example 3.37 is definitely an example of this (and what I replaced with Example 3.40 in my thesis as the amendment).  The claim is that trans3(G) = {0}.  We can already see this is false, because 1/2 goes to 0, implying 1/2 âï¸ trans3(G) as well.  I am pretty certain the set of dyadic rationals will be 3-transitive, so in fact trans3(G) will be dense. For the dyadic rationals in the interval [0, 1/2] have a path to 0 under f1, and the dyadic rationals in [1/2, 1] have a path to 1/2 under f2 (and consequently a path to 0 under G).

Each example in~\cite{transitivity_CR} has $\transitiveIn{G}{2} = \emptySet$ whenever $\transitiveIn{G}{3} \neq \transitiveIn{G}{2}$. Example $3.40$\footnote{We note Example $3.40$ in~\cite{transitivity_CR} claims $\transitiveIn{G}{3} = \set{0}$, and is therefore non-dense.   We notice that the dyadic rationals in $[0, \Fraction 1 / 2]$ have a path to $0$ under $f_1$, and the dyadic rationals in $[\Fraction 1 / 2, 1]$ have a path to $\Fraction 1 / 2$ under $f_2$, and consequently a path to $0$ under $G$. Hence, $\transitiveIn{G}{3}$ is the set of dyadic rationals in $[0, 1]$.} and Example $3.41$ in~\cite{transitivity_CR} explicitly show $\transitiveIn{G}{2} = \emptySet$.  Indeed, it is straightforward to check for the other examples. We will explicitly show $\transitiveIn{G}{2} = \emptySet$ for~\cite[Example $3.24$]{transitivity_CR}.   Our proof leads to a property of CR-dynamical systems for which $\transitiveIn{G}{3}$ sets are dense. We first prove for each positive integer $n$, $\p{3, n}$-transitive points do not exist when there are $2$-transitive points. 

\begin{proposition}\label{prop:2-trans-implies-not-3n}
Let $\p{X, G}$ be a CR-dynamical system. If $\transitiveIn{G}{2} \neq \emptySet$, then $\transitiveIn{G}{3} \setMinus \transitiveIn{G}{2} = \transitiveIn{G}{\p{3, \omega}}$. 
\end{proposition}
\begin{proof}
Suppose $x \in \transitiveIn{G}{3} \setMinus \transitiveIn{G}{2}$ and  $y \in \transitiveIn{G}{2}$. To derive a contradiction, suppose there exists a positive integer $n$ such that $x \in \transitiveIn{G}{\p{3, n}}$. Therefore, $\mathcal{L}_n\of{T_G\of{x}}^\ast = X$. It follows $y \in  T_G\of{x}^\ast$, and so $x \in T_{G^{-1}}\of{y}^\ast$. Hence, by Proposition~\ref{prop:backwards}, $x \in T_{G^{-1}}\of{y}^\ast \subseteq \transitiveIn{G}{2}$. But $x \in \transitiveIn{G}{2}$ contradicts the fact $x \in \transitiveIn{G}{\p{3, n}}$. Thus, $\transitiveIn{G}{3} \setMinus \transitiveIn{G}{2} = \transitiveIn{G}{\p{3, \omega}}$. 
\end{proof}

\begin{corollary}\label{cor:3n-implies-not-2-trans}
Let $\p{X, G}$ be a CR-dynamical system. If $\transitiveIn{G}{\p{3, n}} \neq \emptySet$ for some positive integer $n$, then $\transitiveIn{G}{2} = \emptySet$.  
\end{corollary}
\begin{proof}
Follows directly from Proposition~\ref{prop:2-trans-implies-not-3n}. 
\end{proof}

\begin{example}[{Example $3.24$~\cite{transitivity_CR}}]\label{example:devils-stair-case}
Let $X = [0, 1]$ and $C$ be the standard ternary Cantor set in $X$. Let $f : X \to X$ be the standard Cantor function (also known as the Devil's Staircase~\cite[Page $131$, Figure $3$-$19$]{hocking_young_topology}), and $G = \p{\set{\Fraction 1 / 2} \times C} \union \Graph{f}$.  Then, $\Fraction 1 / 2 \in \transitiveIn{G}{\p{3, 2}}$. Hence, $\transitiveIn{G}{2} = \emptySet$ by Corollary~\ref{cor:3n-implies-not-2-trans}. 
\end{example}

Exploring Example~\ref{example:devils-stair-case} further, we make the observation $\transitiveIn{G}{\p{3, \omega}} = \emptySet$. Surprisingly, this fact yields an easy way to see $\transitiveIn{G}{3}$ is not dense.

\begin{proposition}\label{prop:ltrans-dense-no-2ltrans-imply-3omega-nonempty}
Let $\p{X, G}$ be a CR-dynamical system such that $\transitiveIn{G}{2} = \emptySet$. If $\transitiveIn{G}{3}$ is dense in $X$, then $\transitiveIn{G}{\p{3, \omega}} \neq \emptySet$.  
\end{proposition}
\begin{proof}
Suppose $\transitiveIn{G}{3}$ is dense in $X$. To derive a contradiction, suppose $\transitiveIn{G}{3} = \Union_{n = 1}^\infinity \transitiveIn{G}{\p{3, n}}$. Hence, there is a path from each $x \in \transitiveIn{G}{3}$ to each $y \in X$ in $G$. In particular, there is a path between every $x, y \in \transitiveIn{G}{3}$ in $G$. %It follows $\transitiveIn{G}{3} \subseteq T_G\of{x}^\ast \intersect T_{G^{-1}}\of{x}^\ast$ for each $x \in \transitiveIn{G}{3}$. Since $\transitiveIn{G}{3}$ is dense in $\legal{G}$, it follows $T_{G^{-1}}\of{x}^\ast$ is dense in $\legal{G^{-1}} = \legal{G}$. Hence, $\transitiveIn{G}{3} \subseteq \transitiveIn{G^{-1}}{3}$. 

Now, let $x \in \transitiveIn{G}{3}$, and let $\setOf{U_n}:{n \in \naturals}$ be a base for $X$. For each $n \in \naturals$, let $x_n \in U_n \intersect \transitiveIn{G}{3}$. There is a path from $x$ to $x_0$ in $G$, and a path from $x_n$ to $x_{n + 1}$ in $G$ for each $n \in \naturals$. But this means we may construct a dense infinite branch of $T_G\of{x}$, implying $x \in \transitiveIn{G}{2}$, a contradiction. 
\end{proof}

\begin{example}
Let $X = [0, 1]$ and $G = X \times X$. Then, $\transitiveIn{G}{2} = \transitiveIn{G}{3} = X$, which implies $\transitiveIn{G}{\p{3, \omega}} = \emptySet$. Thus, it is necessary for $\transitiveIn{G}{2}$ to be empty for Proposition~\ref{prop:ltrans-dense-no-2ltrans-imply-3omega-nonempty}.
\end{example}

We now make two observations on the equivalence relation $\sim_G$. 

\begin{observation}\label{observation:equiv-class-equal-union-3n-trans}
Let $\p{X, G}$ be a CR-dynamical system.  If $x \in \Union_{n = 1}^\infinity \transitiveIn{G}{\p{3, n}}$, then $[x]_G = \Union_{n = 1}^\infinity \transitiveIn{G}{\p{3, n}}$.
\end{observation}

\begin{observation}\label{observation:equiv-class-3omega}
Let $\p{X, G}$ be a CR-dynamical system.  If $x \in \transitiveIn{G}{\p{3, \omega}}$, then $[x]_G \subseteq \transitiveIn{G}{\p{3, \omega}}$. 
\end{observation}

Banic et al. in~\cite{transitivity_CR} remark below Example $4.5$ that there are CR-dynamical systems $\p{X, G}$ such that $\transitiveIn{G}{\p{3, n}} \neq \emptySet$ and $\transitiveIn{G}{\p{3, n}}$ is not dense in $X$. Example~\ref{example:devils-stair-case} is such an example. We find this is always the case, that if $\p{X, G}$ is a CR-dynamical system, then for all $n$, $\transitiveIn{G}{\p{3, n}}$ is not dense. 

\begin{proposition}
Let $\p{X, G}$ be a CR-dynamical system.  Then $\Union_{n = 1}^\infinity \transitiveIn{G}{\p{3, n}}$ is not dense in $X$. 
\end{proposition}
\begin{proof}
We may assume $\Union_{n = 1}^\infinity \transitiveIn{G}{\p{3, n}}$ is non-empty, otherwise we are done. Fix $x \in \Union_{n = 1}^\infinity \transitiveIn{G}{\p{3, n}}$. By Observation~\ref{observation:equiv-class-equal-union-3n-trans}, 
\[
[x]_G = \Union_{n = 1}^\infinity \transitiveIn{G}{\p{3, n}}. 
\]
By Corollary~\ref{cor:[x]G-dense-imply-contained-2trans},  $[x]_G$ must not be dense because $x \notin \transitiveIn{G}{2}$. Thus, the result follows. 
\end{proof}

We now explore conditions which imply every $3$-transitive point is $2$-transitive.

\begin{proposition}\label{prop:3-trans-and-isolated-implies-3-trans=backward-orbit-x}
Let $\p{X, G}$ be a CR-dynamical system such that $\isolated{X} \neq \emptySet$. If $x \in \transitiveIn{G}{3} \intersect \isolated{X}$, then $\transitiveIn{G}{3} = T_{G^{-1}}\of{x}^\ast$. 
\end{proposition}
\begin{proof}
Suppose $x \in \transitiveIn{G}{3} \intersect \isolated{X}$. By Proposition~\ref{prop:backwards}, $T_{G^{-1}}\of{x}^\ast \subseteq \transitiveIn{G}{3}$. Let $y \in \transitiveIn{G}{3}$. There is a path from $y$ to $x$ in $G$ by Proposition~\ref{prop:path-isolated-trans2}. Therefore, there is a path from $x$ to $y$ in $G^{-1}$, which implies $y \in T_{G^{-1}}\of{x}^\ast$. Thus, $\transitiveIn{G}{3} = T_{G^{-1}}\of{x}^\ast$. 
\end{proof}

\begin{proposition}
Let $\p{X, G}$ be a CR-dynamical system. If $\isolated{X} \neq \emptySet$ and $\transitiveIn{G}{3} \intersect \transitiveIn{G^{-1}}{3}$ is dense in $X$, then $\transitiveIn{G}{2} = \transitiveIn{G}{3}$. 
\end{proposition}
\begin{proof}
By Theorem~\ref{thm:fundamental-isolated}, we need only show $\transitiveIn{G}{2} \neq \emptySet$. Let $x \in \isolated{X}$. By Proposition~\ref{prop:sufficient-for-ltrans3-dense}, it follows $x \in \transitiveIn{G}{3} \intersect \transitiveIn{G^{-1}}{3}$.  By Proposition~\ref{prop:3-trans-and-isolated-implies-3-trans=backward-orbit-x}, $T_{G^{-1}}\of{x}^\ast = \transitiveIn{G}{3}$ and $T_G\of{x}^\ast = \transitiveIn{G^{-1}}{3}$. Hence, 
\[
[x]_G = T_{G^{-1}}\of{x}^\ast \intersect T_{G}\of{x}^\ast = \transitiveIn{G}{3} \intersect \transitiveIn{G^{-1}}{3}
\]
is dense in $X$. By Corollary~\ref{cor:[x]G-dense-imply-contained-2trans}, $[x]_G \subseteq \transitiveIn{G}{2}$, and we are done. 
\end{proof}

\begin{proposition}\label{prop:complement-2trans-dense-if-3-trans-neq-2-trans}
Let $\p{X, G}$ be a CR-dynamical system. If $\transitiveIn{G}{2} \neq \transitiveIn{G}{3}$, then 
\[
X \setMinus \transitiveIn{G}{2} = \intrans{G} \union \transitiveIn{G}{3} \setMinus \transitiveIn{G}{2}
\]
is dense in $X$. 
\end{proposition}
\begin{proof}
Let $x \in \transitiveIn{G}{3} \setMinus \transitiveIn{G}{2}$. Then, $T_G\of{x}^\ast$ is dense in $X$. Observe 
\[
T_G\of{x}^\ast \subseteq X \setMinus \transitiveIn{G}{2},
\]
since if $y \in T_G\of{x}^\ast \intersect \transitiveIn{G}{2}$, then $x \in T_{G^{-1}}\of{y}^\ast \subseteq \transitiveIn{G}{2}$, a contradiction. 
\end{proof}

\begin{corollary}
Let $\p{X, G}$ be a CR-dynamical system. If $\interior{\transitiveIn{G}{2}} \neq \emptySet$, then $\transitiveIn{G}{2} = \transitiveIn{G}{3}$. 
\end{corollary}
\begin{proof}
Suppose $\interior{\transitiveIn{G}{2}} \neq \emptySet$. To derive a contradiction, suppose $\transitiveIn{G}{2} \neq \transitiveIn{G}{3}$. By Proposition~\ref{prop:complement-2trans-dense-if-3-trans-neq-2-trans}, $X \setMinus \transitiveIn{G}{2}$ is dense in $X$. However, this implies $X \setMinus \transitiveIn{G}{2}$ meets every non-empty open set in $X$. By assumption,  $\interior{\transitiveIn{G}{2}} \neq \emptySet$ in $X$, which intersects trivially with $X \setMinus \transitiveIn{G}{2}$, a contradiction. Thus, $\transitiveIn{G}{2} = \transitiveIn{G}{3}$. 
\end{proof}

\begin{proposition}\label{prop:3-trans-imply-2-trans-if-contain-all-branches}
Let $\p{X, G}$ be a CR-dynamical system, and let $x \in X$. If $T_G\of{x}^\ast \subseteq \transitiveIn{G}{3}$, then $x \in \transitiveIn{G}{2}$. 
\end{proposition}
\begin{proof}
Let $x_0 := x$, and $\setOf{U_n}:{n \in \naturals}$ be a base for $X$. There exists $n_0 \in \naturals$ such that $G^{n_0}\of{x_0} \intersect U_0$ is non-empty. Let $x_1 \in G^{n_0}\of{x_0} \intersect U_0$. Since $x_1 \in T_G\of{x_0}^\ast \subseteq \transitiveIn{G}{2}$, there exists $n_1 \in \naturals$ such that $G^{n_1}\of{x_1} \intersect U_1$ is non-empty. Let $x_2 \in G^{n_1}\of{x_1} \intersect U_1$. Continuing in this fashion yields a dense infinite branch of $T_G\of{x}$, which implies $x \in \transitiveIn{G}{2}$. 
\end{proof}

\begin{corollary}
Let $\p{X, G}$ be a CR-dynamical system. If $G = G^{-1}$, then $\transitiveIn{G}{2} = \transitiveIn{G}{3}$. 
\end{corollary}
\begin{proof}
Suppose $x \in \transitiveIn{G}{3}$. Then, by Proposition~\ref{prop:backwards}, $T_{G^{-1}}\of{x}^\ast \subseteq \transitiveIn{G}{3}$. Since $G = G^{-1}$, $T_{G}\of{x}^\ast \subseteq \transitiveIn{G}{3}$. By Proposition~\ref{prop:3-trans-imply-2-trans-if-contain-all-branches}, $x \in \transitiveIn{G}{2}$, and we are done.  
\end{proof}

%\[
%[x]_G = T_{G}\of{x}^\ast \intersect T_{G^{-1}}\of{x}^\ast = T_G\of{x}^\ast,
%\]
%since $G = G^{-1}$. Since $x \in \transitiveIn{G}{3}$, $T_G\of{x}^\ast \intersect \legal{G}$ is dense in $\legal{G}$. Hence, $[x]_G \intersect \legal{G}$ is dense in $\legal{G}$, which implies $x \in [x]_G \subseteq \transitiveIn{G}{2}$. 

\begin{figure}[h]
\[
 \begin{tikzpicture}[
       decoration = {markings,
                     mark=at position .5 with {\arrow{Stealth[length=2mm]}}},
       dot/.style = {circle, fill, inner sep=1pt, node contents={},
                     label=#1},
every edge/.style = {draw, postaction=decorate}
                        ]
%\node (a) [dot];
%\node (b) at (5, 0) [dot];
%\node (c) at (5, 5) [dot]; 
%\node (d) at (0, 5) [dot];
%\node (e) at (8, 3) {$\color{red}\text{$G = X \times \set{0} \union \set{0} \times X$}$};
%\node (g) at (6.5, 4) {${\color{red} \text{$X = [0, 1]$}}$};
%\node (f) at (6.68, 2) {${\color{red} \text{$G\of{\Fraction 1 / 2} = \set{0}$}}$};
%\node (h) at (6.55, 1) {${\color{red} \text{$G^2\of{\Fraction 1 / 2} = X$}}$}

\node (a) at (2, -.25) {$\color{red}\text{$G = \Graph{f} \union \Graph{f}^{-1}$}$};

%\draw[dashed] (0, 2) -- (4, 2);
%\draw[dashed] (2, 0) -- (2, 4);
%\draw[dashed] (0, 2.75) -- (4, 2.75);
\draw (4, 0) -- (4, 4);
\draw (4, 4) -- (0, 4);
\draw (0, 0) -- (4, 0);
\draw (0, 0) -- (0, 4);
\draw[blue, very thick] (0, 0) -- (2, 4);
\draw[blue, very thick] (2, 4) -- (4, 0);
\draw[blue, very thick] (0, 4) -- (4, 2);
\draw[blue, very thick] (0, 0) -- (4, 2);
%\draw[blue, very thick] (2, 2) -- (3, 4);
%\draw[blue, very thick] (3, 4) -- (4, 2);
%\draw[blue, very thick] (4, 2) -- (4, 0);
%\node[blue] (a) at (0, 2.75) [dot=left:$x_2$];
 \end{tikzpicture}
\]
\caption{The relation $G$ from Example~\ref{ex:G=inverse-tent}}\label{fig:G=inverse-tent}
\end{figure}

\begin{example}\label{ex:G=inverse-tent}
Let $X = [0, 1]$, and $f : X \to X$ be the tent map  (see Figure~\ref{fig:tent map}). Let $G$ be the closed relation
\[
G = \Graph{f} \union \Graph{f}^{-1},
\]
depicted in Figure~\ref{fig:G=inverse-tent}. Then, since $G = G^{-1}$, it follows $\transitiveIn{G}{2} = \transitiveIn{G}{3}$. Moreover, there exist $2$-transitive points; we observe $0 \in \transitiveIn{\Graph{f}^{-1}}{3}$, since $\Union_{n \in \naturals} f^{-n}\of{0}$ is the set of all dyadic rationals in $X$, which are dense. It follows the set of dyadic rationals are $2$-transitive. Also, the tent map has transitive points (which are contained in the irrationals), which will also be $2$-transitive in $\p{X, G}$. 
\end{example}

\begin{proposition}\label{prop:ltrans3G-contains-inverse-imply-inverse-in-ltrans2}
Let $\p{X, G}$ be a CR-dynamical system. If $\transitiveIn{G^{-1}}{3} \subseteq \transitiveIn{G}{3}$, then $\transitiveIn{G^{-1}}{3} \subseteq \transitiveIn{G}{2}$. 
\end{proposition}
\begin{proof}
Suppose $x \in \transitiveIn{G^{-1}}{3}$. Then, by Proposition~\ref{prop:backwards},
\[
T_G\of{x}^\ast \subseteq \transitiveIn{G^{-1}}{3} \subseteq \transitiveIn{G}{3}. 
\]
It follows $x \in \transitiveIn{G}{2}$ by Proposition~\ref{prop:3-trans-imply-2-trans-if-contain-all-branches}. Thus, $\transitiveIn{G^{-1}}{3} \subseteq \transitiveIn{G}{2}$. 
\end{proof}

\begin{corollary}\label{cor:2trans=3trans-if-3trans=inverse}
Let $\p{X, G}$ be a CR-dynamical system. If $\transitiveIn{G}{3} = \transitiveIn{G^{-1}}{3}$, then $\transitiveIn{G}{2} = \transitiveIn{G}{3}$. 
\end{corollary}
\begin{proof}
Follows directly from Proposition~\ref{prop:ltrans3G-contains-inverse-imply-inverse-in-ltrans2}. 
\end{proof}

%\begin{proposition}
%Let $\p{X, G}$ be a CR-dynamical system such that $\legal{G}$ is second countable. Then $\transitiveIn{G}{3} \intersect \transitiveIn{G^{-1}}{3} \subseteq \transitiveIn{G}{2}$. 
%\end{proposition}
%\begin{proof}
%Suppose $x \in \transitiveIn{G}{3} \intersect \transitiveIn{G^{-1}}{3}$. Then, $\legal{G} = \legal{G^{-1}}$ and $\transitiveIn{G}{3}$ is dense in $\legal{G}$. Consider  $[x]_G = T_G\of{x}^\ast \intersect T_{G^{-1}}\of{x}^\ast$. Let $\setOf{U_n}:{n \in \omega}$ be a base for $\legal{G}$. 
%\end{proof}

We now give~\cite[Definition $4.8$]{transitivity_CR}. 

\begin{definition}[{Definition $4.8$~\cite{transitivity_CR}}]\label{definition:3omega}
Let $\p{X, G}$ be a CR-dynamical system, $x \in \transitiveIn{G}{\p{3, \omega}}$ and $n > 1$. We say that $x$ is 
\begin{enumerate}
    \item {\color{cyan} \emph{$\p{3, \omega, n}$-transitive in $\p{X, G}$}}, if there is $\collection{B} \subseteq \mathcal{B}_\infinity\of{T_G\of{x}}$ such that
    \begin{itemize}
        \item $\cardinality{\collection{B}} = n$;
        \item $\closure{\Union_{B \in \collection{B}} B^\ast} = X$; and
        \item if $\collection{B}^\prime \subseteq \mathcal{B}_\infinity\of{T_G\of{x}}$ and $\closure{\Union_{B \in \collection{B}^\prime} B^\ast} = X$, then $\cardinality{\collection{B}^\prime} \geq n$. 
    \end{itemize}
    We use {\color{cyan} $\transitiveIn{G}{\p{3, \omega, n}}$} to denote the set of $\p{3, \omega, n}$-transitive points in $\p{X, G}$. 
    \item {\color{cyan} \emph{$\p{3, \omega, \omega}$-transitive in $\p{X, G}$}}, if
    \[
    x \notin \Union_{n=1}^\infinity \transitiveIn{G}{\p{3, \omega, n}}. 
    \]
    We use {\color{cyan} $\transitiveIn{G}{\p{3, \omega, \omega}}$} to denote the set of $\p{3, \omega, \omega}$-transitive points in $\p{X, G}$. 
\end{enumerate}
\end{definition}
\begin{note}
Observe that $n > 1$ since otherwise, any $\p{3, \omega, 1}$-transitive point has a dense infinite branch and thus is a $2$-transitive point. 
\end{note}

%For us, $\p{3, \omega, n}$-transitive points correspond to $\p{3, \omega, n + 1}$-transitive points in~\cite{transitivity_CR}. Our reasoning for this, is that there do not exist any $\p{3, \omega, 1}$-transitive points given their definition, since this would mean the given point has a dense infinite branch (and thus is a $2$-transitive point). 
%\end{note}

\begin{figure}[h]
\[
 \begin{tikzpicture}[
       decoration = {markings,
                     mark=at position .5 with {\arrow{Stealth[length=2mm]}}},
       dot/.style = {circle, fill, inner sep=2.4pt, node contents={},
                     label=#1},
every edge/.style = {draw, postaction=decorate}
                        ]
\node[cyan] (x) [dot=below:$1$];
\node (x0) at (-3, 1.5) [dot=left:$x_0$];
\node (x1) at (-1, 1.5) [dot=left:$x_1$]; 
\node (xn) at (1, 1.5) [dot=right:$x_{n-1}$]; 
\node (an) at (3, 1.5) [dot=right:$x_n$];
\node (d) at (0, 1) {$\ldots$}; 
\node (d) at (0, 2.5) {$\ldots$}; 
\node (d) at (0, 4) {$\ldots$}; 
\node (d) at (0, 5.5) {$\ldots$}; 
\node (y0) at (-3, 3) [dot=left:$x_0$];
\node (y1) at (-1, 3) [dot=left:$x_1$];
\node (yn) at (1, 3) [dot=right:$x_{n-1}$];
\node (bn) at (3, 3) [dot=right:$x_{n+1}$];
\node (z0) at (-3, 4.5) [dot=left:$x_0$];
\node (z1) at (-1, 4.5) [dot=left:$x_1$];
\node (zn) at (1, 4.5) [dot=right:$x_{n-1}$];
\node (cn) at (3, 4.5) [dot=right:$x_{n+2}$];
\node (d1) at (-3, 6) {$\vdots$}; 
\node (d2) at (-1, 6) {$\vdots$}; 
\node (d3) at (1, 6) {$\vdots$}; 
\node (d4) at (3, 6) {$\vdots$}; 

\draw (x) -- (an);
\draw (an) -- (bn);
\draw (bn) -- (cn);
\draw   (x) -- (x0);
\draw   (x) -- (x1);
\draw   (x) -- (xn);
\draw (x0) -- (y0);
\draw (x1) -- (y1);
\draw (xn) -- (yn);
\draw (y0) -- (z0);
\draw (z1) -- (y1);
\draw (zn) -- (yn);
\draw (z0) -- (-3, 5.5);
\draw (z1) -- (-1, 5.5);
\draw (zn) -- (1, 5.5);
\draw (cn) -- (3, 5.5);
 \end{tikzpicture}
\]
\caption{Transitivity tree $T_{G_n}\of{1}$ in Example~\ref{ex:3-omega-n}}\label{fig:3-omega-n}
\end{figure}

\begin{example}\label{ex:3-omega-n}
Let $X = \set{0, 1} \union \setOf{x_n}:{n \in \naturals}$, where $x_n = \Fraction 1 / {n + 2}$ for each $n \in \naturals$. For each $n \in \naturals$,  let
\[
G_n = \set{\p{0, 0}} \union \setOf{\p{1, x_{k}}}:{k \in [n]} \union \setOf{\p{x_k, x_k}}:{k \in [n-1]} \union \setOf{\p{x_k, x_{k+1}}}:{k \geq n}. 
\]
Then, for each positive integer $n$, $1 \in \transitiveIn{G_n}{\p{3, \omega, n+1}}$. The transitivity tree of $1$ in $\p{X, G_n}$ is depicted in Figure~\ref{fig:3-omega-n}. 
\end{example}

\begin{proposition}\label{prop:equiv-class-contained-in-3-omega-n}
Let $\p{X, G}$ be a CR-dynamical system. Then, 
\[
[x]_G \subseteq \Union_{n = 2}^\infinity \transitiveIn{G}{\p{3, \omega, n}}
\]
for each $x \in \Union_{n = 2}^\infinity \transitiveIn{G}{\p{3, \omega, n}}$. 
\end{proposition}
\begin{proof}
Suppose $x \in \Union_{n = 2}^\infinity \transitiveIn{G}{\p{3, \omega, n}}$.  By Observation~\ref{observation:equiv-class-3omega}, $[x]_G \subseteq \transitiveIn{G}{\p{3, \omega}}$. There exists $k > 1$ such that $x \in \transitiveIn{G}{\p{3, \omega, k}}$. Hence, there is $\collection{B} \subseteq \mathcal{B}_\infinity\of{T_G\of{x}}$ such that $\cardinality{\collection{B}} = k$ and $\Union_{B \in \collection{B}} B^\ast$ is dense in $X$. If $y \in [x]_G$, there is a path from $y$ to $x$ in $G$, which implies there is $\collection{B}^\prime \subseteq \mathcal{B}_\infinity\of{T_G\of{y}}$ such that $\cardinality{\collection{B}^\prime} = k$ and $\Union_{B \in \collection{B}^\prime} B^\ast$ is dense in $X$. It follows $y \in \Union_{n=2}^\infinity \transitiveIn{G}{\p{3, \omega, n}}$. Thus, 
\[
[x]_G \subseteq \Union_{n = 2}^\infinity \transitiveIn{G}{\p{3, \omega, n}}
\]
and we are done. 
\end{proof}

\begin{note}
In light of our proof of Proposition~\ref{prop:equiv-class-contained-in-3-omega-n}, we make the following simple observation. 
\end{note}

\begin{observation}
Let $\p{X, G}$ be a CR-dynamical system, and let $n > 1$. If $x \in \transitiveIn{G}{\p{3, \omega, n}}$, then $[x]_G \subseteq \transitiveIn{G}{\p{3, \omega, n}}$
\end{observation}

\begin{proposition}\label{prop:equiv-class-contained-in-3-omega-omega}
Let $\p{X, G}$ be a CR-dynamical system. Then,
\[
[x]_G \subseteq \transitiveIn{G}{\p{3, \omega, \omega}}
\]
for each $x \in \transitiveIn{G}{\p{3, \omega, \omega}}$. 
\end{proposition}
\begin{proof}
Suppose $x \in \transitiveIn{G}{\p{3, \omega, \omega}}$. It follows $[x]_G \subseteq \transitiveIn{G}{\p{3, \omega}}$, by Observation~\ref{observation:equiv-class-3omega}. If $[x]_G \intersect \p{\Union_{n=2}^\infinity \transitiveIn{G}{\p{3, \omega, n}}} \neq \emptySet$, it would follow by Proposition~\ref{prop:equiv-class-contained-in-3-omega-n} that $x \in \Union_{n=2}^\infinity \transitiveIn{G}{\p{3, \omega, n}}$, which is a contradiction. Thus, it must be the case $[x]_G \subseteq \transitiveIn{G}{\p{3, \omega, \omega}}$, as desired. 
\end{proof}

\begin{proposition}
Let $\p{X, G}$ be a CR-dynamical system. Then, for each $x \in \transitiveIn{G}{\p{3, \omega, \omega}}$, 
\[
T_G\of{x}^\ast \subseteq \intrans{G} \union \transitiveIn{G}{\p{3, \omega, \omega}}. 
\]
\end{proposition}
\begin{proof}
Suppose $x \in \transitiveIn{G}{\p{3, \omega, \omega}}$. Clearly,
\[
T_G\of{x}^\ast \subseteq \intrans{G} \union \transitiveIn{G}{\p{3, \omega}}, %\legal{G} \setMinus \transitiveIn{G}{2}.
\]
since $x$ is neither $2$-transitive nor $\p{3, n}$-transitive for each $n \geq 1$. Furthermore, if $T_G\of{x}^\ast \intersect \transitiveIn{G}{\p{3, \omega, n}} \neq \emptySet$ for some $n > 1$, this would imply $x \in \Union_{k=2}^\infinity \transitiveIn{G}{\p{3, \omega, k}}$, a contradiction. 
\end{proof}

%\begin{theorem}
%Let $\p{X, G}$ be a CR-dynamical system, such that $\legal{G}$ is second countable. If $\transitiveIn{G}{2} \neq \emptySet$, then $\transitiveIn{G}{2} = \transitiveIn{G}{3}$. 
%\end{theorem}
%\begin{proof}
%Let $\setOf{U_n}:{n \in \omega}$ be a base for $\legal{G}$. To derive a contradiction, suppose there exists $x \in \transitiveIn{G}{3} \setMinus \transitiveIn{G}{2}$. Then, $\legal{G} \setMinus \transitiveIn{G}{2}$ is dense in $\legal{G}$. Observe $T_G\of{x}^\ast \intersect \legal{G} \subseteq \legal{G} \setMinus \transitiveIn{G}{2}$. 
%\end{proof}

%---------------------------------------------------------------
\section{Dense orbit transitivity}\label{section:dense-orbit-transitivity}

Dense orbit transitive CR-dynamical systems are introduced by  Banic et al.~\cite[Definition $5.2$]{transitivity_CR}. They generalise dense orbit transitive topological dynamical systems $\p{X, f}$. We firstly recall the definition of dense orbit transitive topological dynamical systems, then provide an equivalent definition to~\cite[Definition $5.2$]{transitivity_CR} (with the inclusion of \emph{type $0$ dense orbit transitive CR-dynamical} to account for $0$-transitive points). 

\begin{definition}
Let $\p{X, f}$ be a topological dynamical system. We say $\p{X, f}$ is {\color{cyan} \emph{DO-transitive}}, if $\text{tr}\of{f} \neq \emptySet$. 
\end{definition}

\begin{definition}\label{definition:dense-orbit-SV}
Let $\p{X, G}$ be a CR-dynamical system. We say that
\begin{enumerate}
    \item $\p{X, G}$ is {\color{cyan} \emph{$0$-DO-transitive}}, if $\transitiveIn{G}{0} \neq \emptySet$;
    \item $\p{X, G}$ is {\color{cyan} \emph{$1$-DO-transitive}}, if $\transitiveIn{G}{1} \neq \emptySet$;
    \item $\p{X, G}$ is {\color{cyan} \emph{$2$-DO-transitive}}, if $\transitiveIn{G}{2} \neq \emptySet$;
    \item $\p{X, G}$ is {\color{cyan} \emph{$3$-DO-transitive}}, if $\transitiveIn{G}{3} \neq \emptySet$;
    \item for each positive integer $n$, $\p{X, G}$ is {\color{cyan} \emph{$\p{3, n}$-DO-transitive}}, if $\transitiveIn{G}{\p{3, n}} \neq \emptySet$;
    \item $\p{X, G}$ is {\color{cyan} \emph{$\p{3, \omega}$-DO-transitive}}, if $\transitiveIn{G}{\p{3, \omega}} \neq \emptySet$;
    \item for each $n > 1$, $\p{X, G}$ is {\color{cyan} \emph{$\p{3, \omega, n}$-DO-transitive}}, if $\transitiveIn{G}{\p{3, \omega, n}} \neq \emptySet$;
    \item $\p{X, G}$ is {\color{cyan} \emph{$\p{3, \omega, \omega}$-DO-transitive}}, if $\transitiveIn{G}{\p{3, \omega, \omega}} \neq \emptySet$.
\end{enumerate}
\end{definition}

We now make a few observations, with the first two analogous to Observation $5.3$ and Observation $5.4$ in~\cite{transitivity_CR}, respectively. 

\begin{observation}
Let $\p{X, G}$ be a CR-dynamical system and $k \in \set{0, 1, 2, 3}$. If $\p{X, G}$ is $k$-DO-transitive, then $\p{X, G}$ is $\ell$-DO-transitive for each $\ell \in \set{0, 1, 2, 3}$ such that $\ell \geq k$. 
\end{observation}

\begin{observation}
Let $\p{X, f}$ be a topological dynamical system.  Then, the following are equivalent. 
\begin{enumerate}
\item $\p{X, f}$ is DO-transitive. 
\item $\p{X, \Graph{f}}$ is $0$-DO-transitive. 
\item $\p{X, \Graph{f}}$ is $1$-DO-transitive. 
\item $\p{X, \Graph{f}}$ is $2$-DO-transitive. 
\item $\p{X, \Graph{f}}$ is $3$-DO-transitive. 
\end{enumerate}
\end{observation}

\begin{observation}
Let $\p{X, G}$ be a CR-dynamical system and $n > 1$.  If $\p{X, G}$ is $\p{3, \omega, n}$-DO-transitive, then $\p{X, G}$ is $\p{3, \omega}$-DO-transitive.  
\end{observation}

\begin{observation}
Let $\p{X, G}$ be a CR-dynamical system. If $\p{X, G}$ is $\p{3, \omega, \omega}$-DO-transitive, then $\p{X, G}$ is $\p{3, \omega}$-DO-transitive.
\end{observation}

\begin{observation}
Let $\p{X, G}$ be a CR-dynamical system.  If $\p{X, G}$ is $\p{3, \omega}$-DO-transitive or $\p{3, n}$-DO-transitive for some positive integer $n$, then $\p{X, G}$ is $3$-DO-transitive. 
\end{observation}

Recall Example~\ref{ex:3n-transitive}, $1 \in \transitiveIn{G_n}{\p{3, n+2}}$ for each $n \in \naturals$, and hence $\p{X, G_n}$ is $\p{3, n+2}$-DO-transitive. However, it is not $\p{3, \omega}$-DO-transitive. 

\begin{note}
Example $3.40$ in~\cite{transitivity_CR} gives a $\p{3, \omega}$-DO-transitive CR-dynamical system which is not $\p{3, n}$-DO-transitive for each positive integer $n$. 
\end{note}

In Example~\ref{ex:3-omega-n}, for each positive integer $n$, $1 \in \transitiveIn{G_n}{\p{3, \omega, n+1}}$. 
Hence, $\p{X, G_n}$ is $\p{3, \omega, n+1}$-DO-transitive. However, it is not $\p{3, \omega, \omega}$-DO-transitive. 

\begin{note}
Example $3.41$ in~\cite{transitivity_CR} gives a $\p{3, \omega, \omega}$-DO-transitive CR-dynamical system which is not $\p{3, \omega, n}$-DO-transitive for each $n > 1$. 
\end{note}

\begin{corollary}
Let $\p{X, G}$ be a CR-dynamical system. If $\p{X, G}$ is $2$-DO-transitive, then $\p{X, G}$ is not $\p{3, n}$-DO-transitive for each positive integer $n$. 
\end{corollary}
\begin{proof}
Follows directly from Proposition~\ref{prop:2-trans-implies-not-3n}
\end{proof}

\begin{corollary}
Let $\p{X, G}$ be a CR-dynamical system such that $\isolated{X} \neq \emptySet$. If $\p{X, G}$ is $2$-DO-transitive, then $\p{X, G}$ is not
\begin{itemize}
    \item $\p{3, n}$-DO-transitive for each positive integer $n$;
    \item $\p{3, \omega}$-DO-transitive;
    \item $\p{3, \omega, n}$-DO-transitive for each $n > 1$; 
    \item $\p{3, \omega, \omega}$-DO-transitive.
\end{itemize}
\end{corollary}
\begin{proof}
Follows directly from Theorem~\ref{thm:fundamental-isolated}. 
\end{proof}

The following is a well-known fact for topological dynamical systems. 

\begin{theorem}
Let $\p{X, f}$ be a topological dynamical system. If $\isolated{X} = \emptySet$ and $\p{X, f}$ is DO-transitive, then $f$ is surjective. 
\end{theorem}

Banic et al. generalise the above result in~\cite[Theorem $5.6$]{transitivity_CR}.  

\begin{theorem}
Let $\p{X, G}$ be a CR-dynamical system, such that $X$ has no isolated points or $X$ is a singleton. Then, for each $k \in \set{0, 1, 2, 3}$, 
\[
\emph{$\p{X, G}$ is $k$-DO-transitive} \Longrightarrow \pi_0\of{G} = \pi_1\of{G} = X. 
\]
\end{theorem}
\begin{proof}
Let $k \in \set{0, 1, 2, 3}$. Suppose $\p{X, G}$ is $k$-DO-transitive. Then, $\p{X, G}$ is $3$-DO-transitive. The result now follows from~\cite[Theorem $5.6$]{transitivity_CR}. 
\end{proof}

We conclude this section by considering the case when there are isolated points. 

\begin{proposition}\label{prop:ltrans3-dense-iff-inverse-is-dense}
Let $\p{X, G}$ be a CR-dynamical system, such that $\isolated{X} \neq \emptySet$. Then, $\transitiveIn{G}{3}$ is dense in $X$ if, and only if, $\transitiveIn{G^{-1}}{3}$ is dense in $X$. 
\end{proposition}
\begin{proof}
Suppose $\transitiveIn{G}{3}$ is dense in $X$. It follows $\isolated{X} \subseteq \transitiveIn{G}{3}$. Let $x \in \isolated{X}$. Observe $T_{G^{-1}}\of{x}^\ast = \transitiveIn{G}{3}$ by Proposition~\ref{prop:3-trans-and-isolated-implies-3-trans=backward-orbit-x}. Since $\transitiveIn{G}{3}$ is dense in $X$, it follows $x \in \transitiveIn{G^{-1}}{3}$. Hence, $T_G\of{x}^\ast = \transitiveIn{G^{-1}}{3}$, since $x \in \isolated{X}$ and $x \in \transitiveIn{G^{-1}}{3}$ (applying Proposition~\ref{prop:3-trans-and-isolated-implies-3-trans=backward-orbit-x}). As $x \in \transitiveIn{G}{3}$, $T_G\of{x}^\ast$ is dense in $X$. Hence, $\transitiveIn{G^{-1}}{3}$ is dense in $X$. The converse holds similarly, and we are done. 
\end{proof}

\begin{corollary}
Let $\p{X, G}$ be a CR-dynamical system, such that $\isolated{X} \neq \emptySet$.  Then, if $\transitiveIn{G}{3}$ is dense in $X$, then $\pi_0\of{G} = \pi_1\of{G} = X$. 
\end{corollary}
\begin{proof}
Suppose $\transitiveIn{G}{3}$ is dense in $X$. By Proposition~\ref{prop:ltrans3-dense-iff-inverse-is-dense}, it follows $\transitiveIn{G^{-1}}{3}$ is dense in $X$.  By Theorem~\ref{thm:legally-transitive-point-motivation}, it follows $\pi_0\of{G} = \pi_1\of{G} = X$. 
\end{proof}

\begin{proposition}\label{prop:2trans-dense-iff-isolated-in-trans3-inverse}
Let $\p{X, G}$ be a $2$-DO-transitive CR-dynamical system. Then, $\transitiveIn{G}{2}$ is dense in $X$ if, and only if, $\isolated{X} \subseteq \transitiveIn{G^{-1}}{3}$.
\end{proposition}
\begin{proof}
If $\isolated{X} = \emptySet$, then $\transitiveIn{G}{2}$ is dense in $X$ by Proposition~\ref{prop:no-isolated-implies-transitive-dense}, and we are done. Therefore, we assume $\isolated{X} \neq \emptySet$. By Theorem~\ref{thm:fundamental-isolated}, there exists $x \in \isolated{X} \intersect \transitiveIn{G}{2}$ such that $T_G^{-1}\of{x}^\ast = \transitiveIn{G}{2} = \transitiveIn{G}{3}$. 

Now, suppose $\transitiveIn{G}{2} = \transitiveIn{G}{3}$ is dense in $X$. By Proposition~\ref{prop:ltrans3-dense-iff-inverse-is-dense}, it follows $\transitiveIn{G^{-1}}{3}$ is dense in $X$. Hence, $\isolated{X} \subseteq \transitiveIn{G^{-1}}{3}$. 

Conversely, suppose $\isolated{X} \subseteq \transitiveIn{G^{-1}}{3}$. If $x \in \isolated{X} \subseteq \transitiveIn{G^{-1}}{3}$, then  $T_{G^{-1}}\of{x}^\ast$ is dense in $X$. As $T_{G^{-1}}\of{x}^\ast = \transitiveIn{G}{2}$, it follows $\transitiveIn{G}{2}$ is dense in $X$, and we are done. 
\end{proof}

%---------------------------------------------------------------
\section{Transitivity}\label{section:transitivity}
In $1920$, G. D. Birkhoff introduced topological transitivity for flows~\cite{scholarpedia_topological_transitivity}, an important property in the study of chaos~\cite{Crannell1995TheRO, silverman_chaos}, and is well-studied in topological dynamics~\cite{akin2016variationsconcepttopologicaltransitivity, aspects_of_topological_transitivity_survey_sergiy, AKIN20122815, scholarpedia_topological_transitivity}.  Transitive CR-dynamical systems are introduced by Banic et al. in~\cite{transitivity_CR}. They generalise the following notion of a transitive topological dynamical system.  

\begin{definition}\label{definition:transitive-top-dyn-system}
Let $\p{X, f}$ be a topological dynamical system. We say $\p{X, f}$ is {\color{cyan} \emph{transitive}}, if for each pair of non-empty open sets $U$ and $V$ in $X$, there exists $n \in \naturals$ such that $f^n\of{U} \intersect V \neq \emptySet$. 
\end{definition}

It is well known $\p{X, f}$ is transitive, where $X = [0, 1]$ and $f : X \to X$ is the tent map defined by 
\[
f\of{x} = \begin{cases}
          2x &
          \text{if $x \in [0, \Fraction 1 / 2]$;} \\
          2 - 2x &
          \text{if $x \in [\Fraction 1 / 2, 1]$;}
          \end{cases}
\]
depicted in Figure~\ref{fig:tent map}. 

Banic et al. in~\cite[Definition $6.3$]{transitivity_CR} give the following natural generalisation for transitivity. 

\begin{definition}
Let $\p{X, G}$ be a CR-dynamical system. We say $\p{X, G}$ is {\color{cyan} \emph{transitive}}, if for each pair of non-empty open sets $U$ and $V$ in $X$, there exists $n \in \naturals$ such that $G^n\of{U} \intersect V \neq \emptySet$. 
\end{definition}

The following well-known classical result does not generalise for transitive CR-dynamical systems. 

\begin{theorem}\label{theorem:classical-transitive-iff-dense-Gdelta}
Let $\p{X, f}$ be a topological dynamical system.  Then, $\p{X, f}$ is transitive if, and only if, $\emph{tr}\of{f}$ is a dense $G_\delta$ set. 
\end{theorem}

We introduce two further types of transitivity for which  Theorem~\ref{theorem:classical-transitive-iff-dense-Gdelta} does generalise,  $0$-transitivity and $1$-transitivity, Propostion \ref{prop:1-transitive-iff-ltrans0-dense}. 
We rename transitive as $2$-transitive in keeping with their hierarchy.  

\begin{definition}
We say a CR-dynamical system $\p{X, G}$ is  {\color{cyan} \emph{$i$-transitive}}, if for each pair of non-empty open sets $U$ and $V$ in $X$, there exists $x \in U$ such that 
\begin{enumerate}
    \item[{\color{red}$\p{i = 0}$}] $\emptySet \neq G^n\of{x} \subseteq V$ for some $n \in \naturals$;  
    \item[{\color{red}$\p{i = 1}$}] $x \in \legal{G}$ and for each $B \in \mathcal{B}_\infinity\of{T_G\of{x}}$, $B^\ast \intersect V \neq \emptySet$;  
    \item[{\color{red}$\p{i = 2}$}] $G^n\of{x} \intersect V \neq \emptySet$ for some $n \in \naturals$.
\end{enumerate} 
\end{definition}

\begin{observation}\label{observation:transitivity-implication}
Let $\p{X, G}$ be a CR-dynamical system. If $\p{X, G}$ is $0$-transitive, then $\p{X, G}$ is $1$-transitive. If $\p{X, G}$ is $1$-transitive, then $\p{X, G}$ is $2$-transitive. 
\end{observation}

\begin{example}\label{ex:2-transitive-but-not-1-transitive}
Let $X = \set{0, 1}$ and $G = \set{\p{0, 0}, \p{1, 0}, \p{0, 1}}$. Indeed, $\p{X, G}$ is $2$-transitive, but it is not $1$-transitive, and hence not $0$-transitive. 
\end{example}

We have that $0$-transitive CR-dynamical systems exist, see Example \ref{ex-0-trans}. We do not know if there is a  $0$-transitive CR-dynamical system that is not $1$-transitive. 

\begin{question}\label{question:trans1-dense-imply-trans0-dense}
Let $\p{X, G}$ be an CR-dynamical system. If $\p{X, G}$ is 1-transitive, must $\p{X, G}$ be 0-transitive? 
\end{question}

It is not surprising that 0-transitivity is closely related to 0-transitive points,  1-transitivity is closely related to 1-transitive points, and  2-transitivity is closely related to 2-transitive points and 3-transitive points.

\begin{proposition}\label{proposition:transk-dense-then-k-trans}
Let $\p{X, G}$ be a CR-dynamical system.  Then, 
\begin{enumerate}
\item if $\transitiveIn{G}{0}$ is dense in $X$, then $\p{X, G}$ is $0$-transitive;
\item if $\transitiveIn{G}{1}$ is dense in $X$, then $\p{X, G}$ is $1$-transitive; 
\item if $\transitiveIn{G}{2}$ is dense in $X$, then $\p{X, G}$ is $2$-transitive; and
\item if $\transitiveIn{G}{3}$ is dense in $X$, then $\p{X, G}$ is $2$-transitive. 
\end{enumerate}
\end{proposition}
\begin{proof}
Suppose $\transitiveIn{G}{0}$ is dense in $X$, and $U, V$ are non-empty open sets in $X$. Then, there exists $x \in \transitiveIn{G}{0} \intersect U$. Since $x$ is $0$-transitive, there exists $n \in \naturals$ such that $G^n\of{x} \subseteq V$. Thus, $\p{X, G}$ is $0$-transitive, yielding $\p{1}$. 

Suppose $\transitiveIn{G}{1}$ is dense in $X$, and $U, V$ are non-empty open sets in $X$. Then, there exists $x \in \transitiveIn{G}{1} \intersect U$. Since $x$ is $1$-transitive, $B^\ast \intersect V \neq \emptySet$ for each $B \in \mathcal{B}_\infinity\of{T_G\of{x}}$. Thus, $\p{X, G}$ is $1$-transitive, yielding $\p{2}$. 

Suppose $\transitiveIn{G}{3}$ is dense in $X$. Let $U$ and $V$ be non-empty open sets in $X$. There exists $x \in \transitiveIn{G}{3} \intersect U$. Since $\Union_{B \in \mathcal{B}_\infinity\of{T_G\of{x}}} B^\ast$ is dense in $X$, it follows there is an $n \in \naturals$ such that $G^n\of{x} \intersect V \neq \emptySet$. Thus, $\p{X, G}$ is $2$-transitive, yielding $\p{4}$. 

Suppose $\transitiveIn{G}{2}$ is dense in $X$, then $\transitiveIn{G}{3}$ is dense in $X$, which implies $\p{X, G}$ is $2$-transitive by $\p{4}$ and we are done.  
\end{proof}

When working with $i$-transitive  CR-dynamical systems, the following Proposition allows us restrict our focus to SV-dynamical systems,  for any $i\in\{0,1,2\}$.
Furthermore, by Proposition \ref{proposition:transk-dense-then-k-trans}, we may restrict to  SV-dynamical systems $\p{X, G}$ if $\transitiveIn{G}{i}$ is dense in $X$, for any $i\in\{0,1,2,3\}$.

\begin{proposition}
Let $\p{X, G}$ be a CR-dynamical system and $i \in \set{0, 1, 2}$. If $\p{X, G}$ is $i$-transitive, then $\p{X, G}$ is an SV-dynamical system. 
\end{proposition}
\begin{proof}
Suppose $\p{X, G}$ is $i$-transitive. By Observation~\ref{observation:transitivity-implication}, $\p{X, G}$ is $2$-transitive. The result now follows from~\cite[Theorem $6.12$]{transitivity_CR}. 
\end{proof}

%\begin{proposition}\label{proposition:if-ltrans2-non-empty-then-legally-2-transitive-iff-ltrans2-dense}
%Let $\p{X, G}$ be a $2$-DO-transitive CR-dynamical system. Then, $\p{X, G}$ is $2$-transitive if, and only if, $\transitiveIn{G}{2}$ is dense in $X$. 
%\end{proposition}
%\begin{proof}
%$\p{\Longrightarrow}$ Suppose $\p{X, G}$ is $2$-transitive. If $\isolated{X} = \emptySet$, then we are done (by Proposition~\ref{prop:no-isolated-implies-transitive-dense}). Therefore, assume $\isolated{X} \neq \emptySet$. Since $\transitiveIn{G}{2} \neq \emptySet$, it follows there exists $x \in \transitiveIn{G}{2} \intersect \isolated{X}$ by Proposition~\ref{proposition:isolated-non-empty-then-transitive-isolated-points-exist}. By Proposition~\ref{prop:backwards},
%\[
%\Union_{k = 0}^\infinity G^{-k}\of{x} \subseteq \transitiveIn{G}{2}, 
%\]
%where the former is dense in $X$ because $x$ is isolated and $\p{X, G}$ is $2$-transitive (equivalence of $\p{1}$ and $\p{7}$ in Theorem~\ref{theorem:legally-2-transitive-equivalence}). Hence, $\transitiveIn{G}{2}$ is dense in $X$. 
%
%$\p{\Longleftarrow}$ If $\transitiveIn{G}{2}$ is dense in $X$, then $\transitiveIn{G}{3}$ is dense in $X$, which implies $\p{X, G}$ is $2$-transitive (by Proposition~\ref{proposition:transk-dense-then-k-trans}). 
%\end{proof}
%

 In this section we explore $i$-transitivity in CR-dynamical systems. 
The following generalises~\cite[Theorem $6.2$]{transitivity_CR}, a well-known classical result.

\begin{theorem}[{Theorem $6.5$~\cite{transitivity_CR}}]\label{theorem:legally-2-transitive-equivalence}
Let $\p{X, G}$ be a CR-dynamical system. Consider the following statements.
\begin{enumerate}
    \item $\p{X, G}$ is $2$-transitive. 
    \item For each pair of non-empty open sets $U$ and $V$ in $X$, there exists a positive integer $n$ such that
    \[
    G^n\of{U} \intersect V \neq \emptySet. 
    \]
    \item For each non-empty open set $U$ in $X$, 
    \[
    \Union_{k=0}^\infinity G^k\of{U}
    \]
    is dense in $X$.
    \item For each non-empty open set $U$ in $X$, 
    \[
    \Union_{k=1}^\infinity G^k\of{U}
    \]
    is dense in $X$.
    \item For each pair of non-empty open sets $U$ and $V$ in $X$, there exists non-negative integer $n$ such that
    \[
    G^{-n}\of{U} \intersect V \neq \emptySet.
    \]
    \item For each pair of non-empty open sets $U$ and $V$ in $X$, there exists positive integer $n$ such that
    \[
    G^{-n}\of{U} \intersect V \neq \emptySet.
    \]
    \item For each non-empty open set $U$ in $X$, 
    \[
    \Union_{k=0}^\infinity G^{-k}\of{U}
    \]
    is dense in $X$.
    \item For each non-empty open set $U$ in $X$, 
    \[
    \Union_{k=1}^\infinity G^{-k}\of{U}
    \]
    is dense in $X$.
\end{enumerate}
Then the following holds. 
\begin{itemize}
    \item $\p{1}$, $\p{3}$, $\p{5}$ and $\p{7}$ are equivalent. 
    \item $\p{2}$, $\p{4}$, $\p{6}$ and $\p{8}$ are equivalent. 
    \item If $\isolated{X}$ is empty, all statements are equivalent.  
\end{itemize}
\end{theorem}

By the equivalence of (1) and (5) we have, $\p{X, G}$ is $2$-transitive if, and only if, $\p{X, G^{-1}}$ is $2$-transitive~\cite[Observation $6.6$]{transitivity_CR}.

We now give a similar result to Theorem~\ref{theorem:legally-2-transitive-equivalence} for $0$-transitivity. 

\begin{theorem}\label{theorem:legally-1-transitive-equivalence}
Let $\p{X, G}$ be an SV-dynamical system, where $X$ has a countable base $\collection{U}$. Consider the following statements. 
\begin{enumerate}
    \item $\p{X, G}$ is $0$-transitive.
    \item For each pair of non-empty open sets $U$ and $V$ in $X$, there exists positive integer $n$ and $x \in U$ such that $G^n\of{x} \subseteq V$. 
  %  \item For each pair of non-empty open sets $U$ and $V$ in $\legal{G}$, there exists non-negative integer $n$ and $x \in U$ such that $G^n\of{x} \subseteq V$. 
  %  \item For each non-empty open set $U$ in $\legal{G}$, 
  %  \[
  %  \p{\Union_{k=0}^\infinity G^k[U]} \intersect \legal{G}
  %  \]
  %  is dense in $\legal{G}$. 
  %  \item For each non-empty open set $U$ in $\legal{G}$, 
  %  \[
  %  \p{\Union_{k=1}^\infinity G^k[U]} \intersect \legal{G}
  %  \]
  %  is dense in $\legal{G}$. 
    \item For each pair of non-empty open sets $U$ and $V$ in $X$, there exists non-negative integer $n$ such that $G^{-n}[U] \intersect V \neq \emptySet$. 
    \item For each pair of non-empty open sets $U$ and $V$ in $X$, there exists positive integer $n$ such that $G^{-n}[U] \intersect V \neq \emptySet$. 
    \item For each non-empty open set $U$ in $X$, 
    \[
    \Union_{k=0}^\infinity G^{-k}[U] 
    \]
    is dense in $X$. 
    \item For each non-empty open set $U$ in $X$, 
    \[
    \Union_{k=1}^\infinity G^{-k}[U] 
    \]
    is dense in $X$. 
    \item $\transitiveIn{G}{0}$ is a dense $G_\delta$ set in $X$. 
    \item $\Intersection_{U \in \collection{U}} \p{\Union_{k = 1}^\infinity  G^{-k}[U]}$ is a dense $G_\delta$ set in $X$. 
    \item $\transitiveIn{G}{0}$ is non-empty. 
    \item $\Intersection_{U \in \collection{U}} \p{\Union_{k = 1}^\infinity G^{-k}[U]}$ is non-empty. 
\end{enumerate}
Then the following holds. 
\begin{itemize}
    \item $\p{1}$, $\p{3}$, $\p{5}$ and $\p{7}$ are equivalent. 
    \item $\p{2}$, $\p{4}$, $\p{6}$ and $\p{8}$ are equivalent. 
    \item If $\isolated{X}$ is empty, all statements are equivalent. 
\end{itemize}
\end{theorem}
\begin{proof}
$\p{\p{1} \Longrightarrow \p{3}}$. Suppose $U$ and $V$ are non-empty open sets in $X$. Then, there exists $x \in V$ and $n \in \naturals$ such that $G^n\of{x} \subseteq U$. That is to say, $x \in  G^{-n}[U] \intersect V$. 

$\p{\p{3} \Longrightarrow \p{5}}$. Suppose $U$ and $V$ are non-empty open sets in $X$. Then, there exists a non-negative integer $n$ such that $G^{-n}[U] \intersect V \neq \emptySet$. Hence, 
\[
\p{\Union_{k=0}^\infinity G^{-k}[U]} \intersect V \neq \emptySet,
\]
which implies $\Union_{k=0}^\infinity  G^{-k}[U]$ is dense in $X$. 

$\p{\p{5} \Longrightarrow \p{7}}$. By Theorem~\ref{thm:0-transitive-points-dense-Gdelta}
\[
\transitiveIn{G}{0} = \Intersection_{U \in \collection{U}} \p{\Union_{k=0}^\infinity  G^{-k}[U]},
\]
which is the countable intersection of dense open sets in $X$. As $X$ is Baire, it follows $\transitiveIn{G}{0}$ is a dense $G_\delta$ set. 

$\p{\p{7} \Longrightarrow \p{1}}$. Follows from Proposition~\ref{proposition:transk-dense-then-k-trans}.

$\p{\p{2} \Longrightarrow \p{4}}$. Suppose $U$ and $V$ are non-empty open sets in $X$. Then, there exists $x \in V$ and positive integer $n$ such that $G^n\of{x} \subseteq U$. That is to say, $x \in G^{-n}[U] \intersect V$.  

$\p{\p{4} \Longrightarrow \p{6}}$. Suppose $U$ and $V$ are non-empty open sets in $X$. Then, there exists a positive integer $n$ such that $G^{-n}[U] \intersect V \neq \emptySet$. Hence, 
\[
\p{\Union_{k=1}^\infinity  G^{-k}[U]} \intersect V \neq \emptySet,
\]
which implies $\Union_{k=1}^\infinity  G^{-k}[U]$ is dense in $X$. 

$\p{\p{6} \Longrightarrow \p{8}}$. Observe 
\[
\Intersection_{U \in \collection{U}} \p{\Union_{k = 1}^\infinity  G^{-k}[U]}
\]
is a countable intersection of dense open sets in $X$. As $X$ is Baire, it follows $\Intersection_{U \in \collection{U}} \p{\Union_{k = 1}^\infinity  G^{-k}[U]}$ is a dense $G_\delta$ set. 

$\p{\p{8} \Longrightarrow \p{2}}$. Suppose $U, V$ are non-empty open sets in $X$. Then, there exists 
\[
x \in \p{\Intersection_{W \in \collection{U}} \p{\Union_{k = 1}^\infinity  G^{-k}[W]}} \intersect U. 
\]
There exists $W \in \collection{U}$ such that $W \subseteq V$, and a positive integer $n$ such that $x \in  G^{-n}[W]$. That is to say, $G^n\of{x} \subseteq W \subseteq V$, and we are done. 

\smallskip

For the remainder of our proof, we assume $\isolated{X} = \emptySet$. 

\smallskip

$\p{\p{1} \Longleftrightarrow \p{2}}$. Clearly, if $\p{2}$ holds, then $\p{1}$ holds. So, suppose $\p{X, G}$ is $0$-transitive. Suppose $U$ and $V$ are open sets in $X$. Since there are no isolated points in $X$, there exists $x \neq y \in X$ such that $x \in U$ and $y \in V$.  There are disjoint open nhoods $W_x$ and $W_y$ of $x$ and $y$, respectively, in $X$. Take $U_x = U \intersect W_x$ and $V_y = V \intersect W_y$. There exists $z \in U_x \subseteq U$ and non-negative integer $n$  such that $G^n\of{z} \subseteq V_y \subseteq V$. As $U_x$ and $V_y$ are disjoint, $n \neq 0$, and thus $\p{2}$ follows. 

$\p{\p{1} \Longleftrightarrow \p{9}}$. As $\p{1}$ is equivalent to $\p{7}$, if  $\p{X, G}$ is $0$-transitive, then $\transitiveIn{G}{0}$ is dense (and therefore non-empty). Conversely, as there are no isolated points, $\transitiveIn{G}{0} \neq \emptySet$ implies $\transitiveIn{G}{0}$ is dense in $X$ (see Proposition~\ref{prop:no-isolated-implies-transitive-dense}), and so by Proposition~\ref{proposition:transk-dense-then-k-trans}, $\p{X, G}$ is $0$-transitive. 

$\p{\p{2} \Longleftrightarrow \p{10}}$. As $\p{2}$ is equivalent to $\p{8}$, if $\p{2}$ holds, then 
\[
\Intersection_{U \in \collection{U}} \p{\Union_{k = 1}^\infinity  G^{-k}[U]}
\]
is dense (and therefore non-empty). Conversely, suppose there exists 
\[
x \in \Intersection_{U \in \collection{U}} \p{\Union_{k = 1}^\infinity  G^{-k}[U]}. 
\]
Let $U$ and $V$  be non-empty open sets in $X$. Then, there exists open $W_1, W_2\in \collection{U}$ such that $W_1 \subseteq U$ and $W_2 \subseteq V$. There exists a positive integer $n_1$ such that $x \in  G^{-n_1}[W_1]$. That is to say, $G^{n_1}\of{x} \subseteq W_1 \subseteq U$. As there are no isolated points, $G^{n_1}\of{x} \subseteq \transitiveIn{G}{0}$ (see Proposition~\ref{proposition:transitive-branches}). Now, let $y \in G^{n_1}\of{x}$. Then, there exists positive integer $n_2$ such that $y \in G^{-n_2}[W_2]$. That is to say, $G^{n_2}\of{y}  \subseteq W_2 \subseteq V$. Hence, $\p{2}$ follows, and we are done. 
\end{proof}

\begin{example}\label{ex-0-trans}
Let $\p{X, G}$ be the CR-dynamical system from Example~\ref{ex:0-transitive-interesting-example}. Since $\transitiveIn{G}{0} \neq \emptySet$ and $\isolated{X} = \emptySet$, $\p{X, G}$ is $0$-transitive by Theorem~\ref{theorem:legally-1-transitive-equivalence}. 
\end{example}

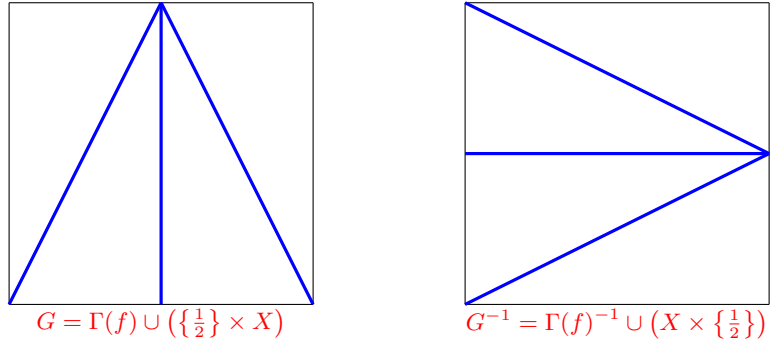
\begin{figure}[h]
\[
 \begin{tikzpicture}[
       decoration = {markings,
                     mark=at position .5 with {\arrow{Stealth[length=2mm]}}},
       dot/.style = {circle, fill, inner sep=1pt, node contents={},
                     label=#1},
every edge/.style = {draw, postaction=decorate}
                        ]
%\node (a) [dot];
%\node (b) at (5, 0) [dot];
%\node (c) at (5, 5) [dot]; 
%\node (d) at (0, 5) [dot];
%\node (e) at (8, 3) {$\color{red}\text{$G = X \times \set{0} \union \set{0} \times X$}$};
%\node (g) at (6.5, 4) {${\color{red} \text{$X = [0, 1]$}}$};
%\node (f) at (6.68, 2) {${\color{red} \text{$G\of{\Fraction 1 / 2} = \set{0}$}}$};
%\node (h) at (6.55, 1) {${\color{red} \text{$G^2\of{\Fraction 1 / 2} = X$}}$}

%\draw[dashed] (0, 2) -- (4, 2);
%\draw[dashed] (2, 0) -- (2, 4);
%\draw[dashed] (0, 2.75) -- (4, 2.75);
\draw (4, 0) -- (4, 4);
\draw (4, 4) -- (0, 4);
\draw (0, 0) -- (4, 0);
\draw (0, 0) -- (0, 4);
\draw[blue, very thick] (0, 0) -- (2, 4);
\draw[blue, very thick] (2, 4) -- (4, 0);
\draw[blue, very thick] (2, 4) -- (2, 0);
%\draw[blue, very thick] (2, 2) -- (3, 4);
%\draw[blue, very thick] (3, 4) -- (4, 2);
%\draw[blue, very thick] (4, 2) -- (4, 0);
%\node[blue] (a) at (0, 2.75) [dot=left:$x_2$];

\node (a) at (2, -.25) {$\color{red}\text{$G = \Graph{f} \union \p{\set{\Fraction 1/2} \times X}$}$};
\node (b) at (8, -.25) {$\color{red}\text{$G^{-1} = \Graph{f}^{-1} \union \p{X \times \set{\Fraction 1 / 2}}$}$};

\draw (10, 0) -- (10, 4);
\draw (10, 4) -- (6, 4);
\draw (6, 0) -- (10, 0);
\draw (6, 0) -- (6, 4);
\draw[blue, very thick] (6, 0) -- (10, 2);
\draw[blue, very thick] (6, 4) -- (10, 2);
\draw[blue, very thick] (6, 2) -- (10, 2);
 \end{tikzpicture}
\]
\caption{The relations $G$ and $G^{-1}$ in Example~\ref{ex:1-trans-inverse-not}}\label{fig:1-trans-inverse-not}
\end{figure}

It is not the case that $\p{X, G}$ is $0$-transitive if, and only if, $\p{X, G^{-1}}$ is $0$-transitive.
\begin{example}\label{ex:1-trans-inverse-not}
%We show that an analogous statement to Corollary~\ref{cor:2-trans-iff-inverse} for $0$-transitive CR-dynamical systems does not hold. 
Let $X = [0, 1]$ and $f : X \to X$ be the tent-map (see Figure~\ref{fig:tent map}). Define $G$ by
\[
G = \Graph{f} \union \p{\set{\Fraction 1/2} \times X}, 
\]
giving us the CR-dynamical system $\p{X, G}$ (where $G$ is depicted on the left in Figure~\ref{fig:1-trans-inverse-not}). 

We claim $\p{X, G}$ is $0$-transitive, but $\p{X, G^{-1}}$ is not. To see why $\p{X, G}$ is $0$-transitive,  let $x \in \text{tr}\of{f}$. Then, $f^k\of{x} \neq 0$ for each $k \in \naturals$, which implies $f^k\of{x} \neq \Fraction 1 / 2$ for each $k \in \naturals$. Therefore, $G^k\of{x} = f^k\of{x}$ for each $k \in \naturals$. It follows $x \in \transitiveIn{G}{0}$, which implies $\transitiveIn{G}{0} \neq \emptySet$. Since there are no isolated points in $X$, it follows from Theorem~\ref{theorem:legally-1-transitive-equivalence} that $\p{X, G}$ is $0$-transitive.  

Now, consider $\p{X, G^{-1}}$, where
\[
G^{-1} = \Graph{f}^{-1} \union \p{X \times \set{\Fraction 1 / 2}}
\]
is depicted on the right in Figure~\ref{fig:1-trans-inverse-not}. 

To see why $\p{X, G^{-1}}$ is not $0$-transitive, we show $\transitiveIn{G^{-1}}{0} = \emptySet$. We do this by observing the stronger fact that $\transitiveIn{G^{-1}}{1} = \emptySet$, since $\sequence{x, \Fraction 1 / 2, \Fraction 1 / 2, \ldots}$ is a trajectory of $x$ for each $x \in X$. 
\end{example}

\begin{note}
In Example~\ref{ex:2-transitive-but-not-1-transitive}, we gave a finite $2$-transitive CR-dynamical system that is not $0$-transitive. Such examples need not be finite.

Let $\p{X, G}$ be the CR-dynamical system from Example~\ref{ex:1-trans-inverse-not}. We found $\p{X, G}$ is $0$-transitive, implying that it is $2$-transitive, and hence, $\p{X, G^{-1}}$ is $2$-transitive. However, we found $\p{X, G^{-1}}$ is not $0$-transitive.  Thus, $\p{X, G^{-1}}$ is a $2$-transitive CR-dynamical system that is not $0$-transitive. 

\end{note}

We now show Theorem \ref{theorem:classical-transitive-iff-dense-Gdelta} holds for 0-transitivity and 1-transitivity. We require the following theorem.

\begin{theorem}[{Theorem $6.12$~\cite{transitivity_CR}}]\label{thm:legally-2-transitive-implies-surjective}
Let $\p{X, G}$ be a CR-dynamical system. If $\p{X, G}$ is $2$-transitive, then 
\[
\pi_0\of{G}  = \pi_1\of{G} =  X.
\]
\end{theorem}

\begin{proposition}\label{prop:1-transitive-iff-ltrans0-dense}
Let $\p{X, G}$ be a CR-dynamical system and $k \in \set{0, 1}$. Then, $\p{X, G}$ is $k$-transitive if, and only if, $\transitiveIn{G}{k}$ is a dense $G_\delta$ set in $X$.
\end{proposition}
\begin{proof}
Suppose $\p{X, G}$ is $0$-transitive. Then, $\p{X, G}$ is $2$-transitive, implying $\p{X, G}$ is an SV-dynamical system (by Theorem~\ref{thm:legally-2-transitive-implies-surjective}). Since $X$ is second countable, it has a countable base $\collection{U}$. By Theorem~\ref{thm:0-transitive-points-dense-Gdelta}, it follows
\[
\transitiveIn{G}{0} = \Intersection_{U \in \collection{U}} \p{\Union_{n = 0}^\infinity G^{-n}[U]}
\]
is a $G_\delta$ set. 

Suppose $U \in \collection{U}$. We claim $D := \Union_{n = 0}^\infinity  G^{-n}[U]$ is dense in $X$. To see why, suppose $V$ is a non-empty open set in $X$. Since $\p{X, G}$ is $0$-transitive, there exists $k \in \naturals$ such that $G^k\of{x} \subseteq U$ for some $x \in V$. It follows $x \in  G^{-k}[U]$. Hence, $x \in V \intersect D$, which establishes our claim. 

Now, $\transitiveIn{G}{0}$ is the countable intersection of dense open sets. Since $X$ is a compact metric space, it is a Baire space. Thus, $\transitiveIn{G}{0}$ is dense in $X$. The converse holds by Proposition~\ref{proposition:transk-dense-then-k-trans}.

Suppose $\p{X, G}$ is $1$-transitive. Then, $\p{X, G}$ is $2$-transitive, implying $\p{X, G}$ is an SV-dynamical system (by Theorem~\ref{thm:legally-2-transitive-implies-surjective}). Since $X$ is second countable, it has a countable base $\collection{U}$. By Proposition~\ref{prop:ltrans1-Gdelta}, it follows
\[
\transitiveIn{G}{1} = \Intersection_{U \in \collection{U}} \setOf{x \in X}:{\forall B \in \mathcal{B}_\infinity\of{T_G\of{x}}, B^\ast \intersect U \neq \emptySet}
\]
is a $G_\delta$ set. 

Suppose $U \in \collection{U}$. We claim $D := \setOf{x \in X}:{B^\ast \intersect U \neq \emptySet, \forall B \in \mathcal{B}_\infinity\of{T_G\of{x}}}$ is dense in $X$. To see why, suppose $V$ is a non-empty open set in $X$. Since $\p{X, G}$ is $1$-transitive, $x \in V$ such that $B^\ast \intersect V \neq \emptySet$ for each $B \in \mathcal{B}_\infinity\of{T_G\of{x}}$. Hence, $x \in V \intersect D$, which establishes our claim.  

Now, $\transitiveIn{G}{1}$ is the countable intersection of dense open sets. Since $X$ is a compact metric space, it is a Baire space. Thus, $\transitiveIn{G}{1}$ is dense in $X$. The converse holds by Proposition~\ref{proposition:transk-dense-then-k-trans}, and we are done. 
\end{proof}

\begin{corollary}\label{cor:legal-second-countable-legally-1-trans-then-ltransk-dense}
Let $\p{X, G}$ be a CR-dynamical system and $k \in \set{0, 1}$. If $\p{X, G}$ is $k$-transitive, then $\transitiveIn{G}{\ell}$ is dense in $X$ for each $\ell \in \set{0, 1, 2, 3}$ with $\ell \geq k$. 
\end{corollary}
\begin{proof}
Follows from Proposition~\ref{prop:1-transitive-iff-ltrans0-dense}.
\end{proof}

Since an SV-dynamical system is $i$-transitive if, and only if, $\transitiveIn{G}{i}$ is dense, where $i \in \set{0, 1}$, we can rephrase Question~\ref{question:trans1-dense-imply-trans0-dense}.

\begin{question}\label{question:trans1-dense-imply-trans0-dense}
Let $\p{X, G}$ be an SV-dynamical system. If $\transitiveIn{G}{1}$ is dense, must $\transitiveIn{G}{0}$ be dense? 
\end{question}

 If $\isolated{X} = \emptySet$, Proposition~\ref{prop:no-isolated-implies-transitive-dense} tells us we can simply ask whether $\transitiveIn{G}{1} \neq \emptySet$ implies $\transitiveIn{G}{0} \neq \emptySet$. %While this remains unresolved, we do know not every $2$-transitive CR-dynamical system is $1$-transitive, as demonstrated by the following example.  %We have an example with $\transitiveIn{G}{1} \neq \emptySet$ and $\transitiveIn{G}{0} = \emptySet$, yet $\isolated{X} = \emptySet$ in this case, and $\transitiveIn{G}{0}$

%As we have seen, Theorem~\ref{theorem:classical-transitive-iff-dense-Gdelta} fails to generalise for $2$-transitive CR-dynamical systems. 
It is not true in general that if  an SV-dynamical system  $\p{X, F}$ is 2-transitive, then $\transitiveIn{F}{3}$ is a dense for $X$. 
However, we do obtain the following for SV-dynamical systems, when adding an additional assumption.

\begin{proposition}\label{prop:SV-2-transitive-iff-3trans-dense-Gdelta}
Let $\p{X, F}$ be an SV-dynamical system such that $F = \Graph{f}$ for some continuous set-valued function $f : X \to 2^X$. Then, $\p{X, F}$ is $2$-transitive if, and only if, $\transitiveIn{F}{3}$ is a dense $G_\delta$ set. 
\end{proposition}
\begin{proof}
Suppose $\p{X, F}$ is $2$-transitive. Let $\collection{U}$ be a countable base for $X$. By Theorem~\ref{theorem:legally-2-transitive-equivalence},  $\Union_{k = 0}^\infinity F^{-k}\of{U}$ is dense for each $U \in \collection{U}$. Since $F$ is the graph of a continuous (in particular, lsc) set-valued function on $X$, it follows $\Union_{k = 0}^\infinity F^{-k}\of{U}$ is open for each $U \in \collection{U}$. By Proposition~\ref{prop:3trans-characterisation}, 
\[
\transitiveIn{F}{3} = \Intersection_{U \in \collection{U}} \p{\Union_{k = 0}^\infinity F^{-k}\of{U}},
\]
which means $\transitiveIn{F}{3}$ is the countable intersection of dense open sets. As $X$ is a Baire space, it follows $\transitiveIn{F}{3}$ is a dense $G_\delta$ set. The converse holds by Proposition~\ref{proposition:transk-dense-then-k-trans}.
\end{proof}

\begin{lemma}[{Lemma $6.4$~\cite{transitivity_CR}}]\label{lemma:legal-no-isolated-2trans-iff-union-dense}
Let $\p{X, G}$ be a CR-dynamical system. If $X$ has no isolated points, then $\p{X, G}$ is $2$-transitive if, and only if, $\Union_{k = 1}^\infinity G^{-k}\of{U}$ is dense in $X$ for each non-empty open set $U$ in $X$. 
\end{lemma}

%\begin{corollary}[Observation $6.6$~\cite{transitivity_CR}]\label{cor:2-trans-iff-inverse}
%Let $\p{X, G}$ be a CR-dynamical system. Then, $\p{X, G}$ is $2$-transitive if, and only if, $\p{X, G^{-1}}$ is $2$-transitive. 
%\end{corollary}
%\begin{proof}
%Follows directly by the equivalence of $\p{1}$ and $\p{5}$ in Theorem~\ref{theorem:legally-2-transitive-equivalence}. 
%\end{proof}
%
%

\begin{theorem}[{Theorem $6.15$~\cite{transitivity_CR}}]
Let $\p{X, G}$ be a CR-dynamical system such that $\transitiveIn{G}{2} \neq \emptySet$. If $X$ has no isolated points, then $\p{X, G}$ is $2$-transitive. 
\end{theorem}

%\begin{observation}\label{observation:stormy-legally-2-transitive}
%Let $\p{X, G}$ be a CR-dynamical system, where $X = [0, 1]$ and $G$ is the relation from Example~\ref{ex:stormy-campsite-map} (see Figure~\ref{fig:Stormy-Campsite-Map}).  Then, $\p{X, G}$ is $2$-transitive from Proposition~\ref{proposition:if-ltrans2-non-empty-then-legally-2-transitive-iff-ltrans2-dense}. 
%\end{observation}

\begin{corollary}\label{corollary:2-DO-then-2-trans-iff-trans2-dense}
Let $\p{X, G}$ be a $2$-DO-transitive CR-dynamical system. Then the following are equivalent. 
\begin{enumerate}
    \item $\p{X, G}$ is $2$-transitive. 
    \item $\transitiveIn{G}{2}$ is dense in $X$.
    \item $\transitiveIn{G}{3}$ is dense in $X$.
\end{enumerate}
\end{corollary}
\begin{proof}
Suppose $\p{X, G}$ is $2$-transitive. If $\isolated{X} = \emptySet$, then we are done (by Proposition~\ref{prop:no-isolated-implies-transitive-dense}). Therefore, assume $\isolated{X} \neq \emptySet$. Since $\transitiveIn{G}{2} \neq \emptySet$, it follows there exists $x \in \transitiveIn{G}{2} \intersect \isolated{X}$ by Proposition~\ref{proposition:isolated-non-empty-then-transitive-isolated-points-exist}. By Proposition~\ref{prop:backwards},
\[
\Union_{k = 0}^\infinity G^{-k}\of{x} \subseteq \transitiveIn{G}{2}, 
\]
where the former is dense in $X$ because $x$ is isolated and $\p{X, G}$ is $2$-transitive (equivalence of $\p{1}$ and $\p{7}$ in Theorem~\ref{theorem:legally-2-transitive-equivalence}). Hence, $\transitiveIn{G}{2}$ is dense in $X$.

 Indeed, the implication from $\p{2}$ to $\p{3}$ is trivial, and the implication from $\p{3}$ to $\p{1}$ was proved in Proposition~\ref{proposition:transk-dense-then-k-trans}. 
\end{proof}

\begin{note}
Let $\p{X, G}$ be the CR-dynamical system where $X = [0, 1]$ and $G = \p{X \times \set{0}} \union \p{\set{0} \times X}$. We have $\transitiveIn{G}{1} = \emptySet$ and $\transitiveIn{G}{2} = X$, which means this result does not extend to $\transitiveIn{G}{1}$. 

Even if $\transitiveIn{G}{1} \neq \emptySet$, $\p{X, G}$ can be $2$-transitive with $\transitiveIn{G}{1}$ not dense in $X$. Take $X = \set{0, 1, 2}$ and $G = \set{\p{0, 1}, \p{1, 2}, \p{2, 0}, \p{2, 2}}$. Then, $\transitiveIn{G}{1} = \set{0}$, which is not dense. Furthermore, $\p{X, G}$ is $2$-transitive, since $\transitiveIn{G}{2} = X$. 

We also observe the CR-dynamical system $\p{X, G}$ in~\cite[Example $6.9$]{transitivity_CR} is $2$-transitive with $\transitiveIn{G}{1} = \set{2}$ not dense. 
\end{note}

\begin{proposition}\label{proposition:interior-3-trans-nonempty-then-legally-2-transitive-iff-ltrans3-dense}
Let $\p{X, G}$ be a CR-dynamical system such that $\interior{\transitiveIn{G}{3}}$ is non-empty. Then, $\p{X, G}$ is $2$-transitive if, and only if, $\transitiveIn{G}{3}$ is dense in $X$. 
\end{proposition}
\begin{proof}
Suppose $\p{X, G}$ is $2$-transitive. Let $U := \interior{\transitiveIn{G}{3}}$. Then, $U$ is a non-empty open set in $X$. Observe
\[
\Union_{k = 0}^\infinity G^{-k}\of{U} \subseteq \transitiveIn{G}{3}, 
\]
where the former is dense in $X$ because $\p{X, G}$ is $2$-transitive (by the equivalence of $\p{1}$ and $\p{7}$ in Theorem~\ref{theorem:legally-2-transitive-equivalence}). Hence, $\transitiveIn{G}{3}$ is dense in $X$. The converse is proved in Proposition~\ref{proposition:transk-dense-then-k-trans}. 
\end{proof}

\begin{proposition}\label{proposition:legally-2-trans-then-isolated-in-3trans}
Let $\p{X, G}$ be a CR-dynamical system. If $\p{X, G}$ is $2$-transitive, then $\isolated{X} \subseteq \transitiveIn{G}{3}$. 
\end{proposition}
\begin{proof}
Suppose $\p{X, G}$ is $2$-transitive, and that $x \in \isolated{X}$. Then, by Theorem~\ref{thm:legally-2-transitive-implies-surjective}
\[
\Union_{B \in \mathcal{B}_\infinity\of{T_G\of{x}}} B^\ast = T_G\of{x}^\ast = \Union_{k=0}^\infinity G^k\of{x},
\]
which is dense in $X$ by the equivalence of $\p{1}$ and $\p{3}$ in Theorem~\ref{theorem:legally-2-transitive-equivalence}. Hence, $x \in \transitiveIn{G}{3}$.  Thus, it follows $\isolated{X} \subseteq \transitiveIn{G}{3}$. 
\end{proof}

\begin{corollary}
Let $\p{X, G}$ be a CR-dynamical system such that $\isolated{X} \neq \emptySet$. Then, $\p{X, G}$ is $2$-transitive if, and only if, $\transitiveIn{G}{3}$ is dense in $X$. 
\end{corollary}
\begin{proof}
Suppose $\p{X, G}$ is $2$-transitive. By Proposition~\ref{proposition:legally-2-trans-then-isolated-in-3trans}, $\isolated{X} \subseteq \transitiveIn{G}{3}$.  Hence, $\interior{\transitiveIn{G}{3}}$ is non-empty, which implies $\transitiveIn{G}{3}$ is dense in $X$ by Proposition~\ref{proposition:interior-3-trans-nonempty-then-legally-2-transitive-iff-ltrans3-dense}. The converse is proved in Proposition~\ref{proposition:transk-dense-then-k-trans}. 
\end{proof}

\begin{figure}[h]
\[
 \begin{tikzpicture}[
       decoration = {markings,
                     mark=at position .5 with {\arrow{Stealth[length=2mm]}}},
       dot/.style = {circle, fill, inner sep=1pt, node contents={},
                     label=#1},
every edge/.style = {draw, postaction=decorate}
                        ]
%\node (a) [dot];
%\node (b) at (5, 0) [dot];
%\node (c) at (5, 5) [dot]; 
%\node (d) at (0, 5) [dot];
%\node (e) at (8, 3) {$\color{red}\text{$G = X \times \set{0} \union \set{0} \times X$}$};
%\node (g) at (6.5, 4) {${\color{red} \text{$X = [0, 1]$}}$};
%\node (f) at (6.68, 2) {${\color{red} \text{$G\of{\Fraction 1 / 2} = \set{0}$}}$};
%\node (h) at (6.55, 1) {${\color{red} \text{$G^2\of{\Fraction 1 / 2} = X$}}$}

\node (a) at (2, -.25) {$\color{red}\text{$G = \Graph{f_1} \union \Graph{f_2} \union \set{\p{0, x_2}, \p{1, x_1}}$}$};

\draw[dashed] (0, 2) -- (4, 2);
\draw[dashed] (2, 0) -- (2, 4);
\draw[dashed] (0, 2.75) -- (4, 2.75);
\draw[dashed] (0, 1.35) -- (4, 1.35);
\draw (4, 0) -- (4, 4);
\draw (4, 4) -- (0, 4);
\draw (0, 0) -- (4, 0);
\draw (0, 0) -- (0, 4);
\draw[blue, very thick] (0, 0) -- (1, 2);
\draw[blue, very thick] (1, 2) -- (2, 0);
\draw[blue, very thick] (2, 2) -- (3, 4);
\draw[blue, very thick] (3, 4) -- (4, 2);
%\draw[blue, very thick] (4, 2) -- (4, 0);
\node[blue] (b) at (4, 1.35) [dot=right:$x_1$];
\node[blue] (a) at (0, 2.75) [dot=left:$x_2$];
 \end{tikzpicture}
\]
\caption{The relation $G$ from~\cite[Example $6.13$]{transitivity_CR}}\label{fig:2-transitive-tent}
\end{figure}
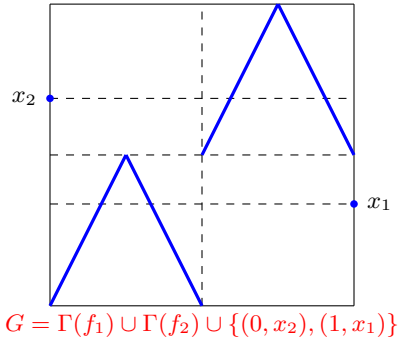

Banic et al. in~\cite[Example $6.13$]{transitivity_CR} give an example of a $2$-transitive SV-dynamical system $\p{X, G}$, such that $X$ has no isolated points and $\transitiveIn{G}{2} = \emptySet$. More precisely, the CR-dynamical system is defined as follows. Let $X = [0, 1]$. Define $f_1 : [0, \Fraction 1 / 2] \to [0, \Fraction 1 / 2]$ by
\[
f_1\of{t} = \begin{cases}
            2t &
            \text{if $t \in [0, \Fraction 1 / 4]$;} \\
            1 - 2t &
            \text{if $t \in [\Fraction 1 / 4, \Fraction 1 / 2]$;}
            \end{cases}
\]
for each $t \in [0, \Fraction 1 / 2]$. Define $f_2 : [\Fraction 1 / 2, 1] \to [\Fraction 1 / 2, 1]$ by 
\[
f_2\of{t} = \begin{cases}
            2t - \Fraction 1 / 2 &
            \text{if $t \in [\Fraction 1 / 2, \Fraction 3 / 4]$;} \\
            \Fraction 5 / 2 - 2t &
            \text{if $t \in [\Fraction 3 / 4, 1]$;}
            \end{cases}
\]
for each $t \in [\Fraction 1 / 2, 1]$. Let $x_1 \in \text{tr}\of{f_1}$ and $x_2 \in \text{tr}\of{f_2}$, and let $G$ be defined by
\[
G = \Graph{f_1} \union \Graph{f_2} \union \set{\p{0, x_2}, \p{1, x_1}},
\]
which is depicted in Figure~\ref{fig:2-transitive-tent}. Their example proves~\cite[Theorem $9$]{sang_salleh_top_transitive_set_valued} is false.  It is shown that 
\begin{itemize}
    \item $\p{X, G}$ is $2$-transitive;
    \item $\transitiveIn{G}{3} \neq \emptySet$;  
    \item $\transitiveIn{G}{2} = \emptySet$; and
    \item $\isolated{X} = \emptySet$.  
\end{itemize}
We note further that 
\begin{itemize}
\item $\closure{\transitiveIn{G}{3}} \neq X$; and
\item $\interior{\transitiveIn{G}{3}} = \emptySet$;
\end{itemize}
as it is easy to check $1 \in \transitiveIn{G}{3}$, and $\transitiveIn{G}{3} \intersect [0, \Fraction 1 / 2] = \emptySet$.  Since $\transitiveIn{G}{3} \neq \emptySet$,  it is natural to ask if $\p{X, G}$ is $2$-transitive, must $\transitiveIn{G}{3} \neq \emptySet$? The following example shows that this is not the case.

\begin{figure}[h]
\[
 \begin{tikzpicture}[
       decoration = {markings,
                     mark=at position .5 with {\arrow{Stealth[length=2mm]}}},
       dot/.style = {circle, fill, inner sep=1pt, node contents={},
                     label=#1},
every edge/.style = {draw, postaction=decorate}
                        ]
%\node (a) [dot];
%\node (b) at (5, 0) [dot];
%\node (c) at (5, 5) [dot]; 
%\node (d) at (0, 5) [dot];
%\node (e) at (8, 3) {$\color{red}\text{$G = X \times \set{0} \union \set{0} \times X$}$};
%\node (g) at (6.5, 4) {${\color{red} \text{$X = [0, 1]$}}$};
%\node (f) at (6.68, 2) {${\color{red} \text{$G\of{\Fraction 1 / 2} = \set{0}$}}$};
%\node (h) at (6.55, 1) {${\color{red} \text{$G^2\of{\Fraction 1 / 2} = X$}}$}

\node (a) at (2, -.25) {$\color{red}\text{$G = \Graph{f_1} \union \Graph{f_2} \union \set{\p{\Fraction 1 / 4, x_2}, \p{1, x_1}}$}$};

\draw[dashed] (0, 2) -- (4, 2);
\draw[dashed] (2, 0) -- (2, 4);
\draw[dashed] (0, 2.75) -- (4, 2.75);
\draw[dashed] (0, 1.35) -- (4, 1.35);
\draw (4, 0) -- (4, 4);
\draw (4, 4) -- (0, 4);
\draw (0, 0) -- (4, 0);
\draw (0, 0) -- (0, 4);
\draw[blue, very thick] (0, 0) -- (1, 2);
\draw[blue, very thick] (1, 2) -- (2, 0);
\draw[blue, very thick] (2, 2) -- (3, 4);
\draw[blue, very thick] (3, 4) -- (4, 2);
%\draw[blue, very thick] (4, 2) -- (4, 0);
\node[blue] (b) at (4, 1.35) [dot=right:$x_1$];
\node[blue] (a) at (1, 2.75) [dot];
\node (a) at (-.2, 2.75) {$x_2$};
 \end{tikzpicture}
\]
\caption{The relation $G$ from Example~\ref{ex:2-transitive-tent-with-no-3trans-points}}\label{fig:2-transitive-tent-with-no-3trans-points}
\end{figure}
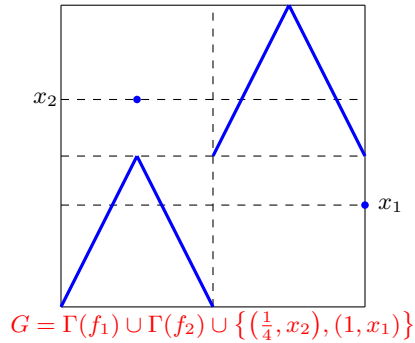

\begin{example}\label{ex:2-transitive-tent-with-no-3trans-points}
Let $X = [0, 1]$.  Define $f_1$ and $f_2$ as above. Let $x_1 \in \text{tr}\of{f_1}$ and $x_2 \in \text{tr}\of{f_2}$, and let $G$ be defined by
\[
G = \Graph{f_1} \union \Graph{f_2} \union \set{\p{\Fraction 1 / 4, x_2}, \p{1, x_1}},
\]
which is depicted in Figure~\ref{fig:2-transitive-tent-with-no-3trans-points}. Now, it is straightforward to check $\p{X, G}$ is $2$-transitive, using a similar argument as seen in~\cite[Example $6.13$]{transitivity_CR}. We show $\transitiveIn{G}{3} = \emptySet$, and hence $\transitiveIn{G}{2} = \emptySet$.  Let $U = \p{0, \Fraction 1 / 4}$ and $V = \p{\Fraction 3 / 4, 1}$. Observe 
\begin{align*}
&& \p{\Union_{k = 0}^\infinity G^{-k}\of{U}} \intersect \p{\Union_{k = 0}^\infinity G^{-k}\of{V}}  &\subseteq \p{\left[0, \Fraction 1 / 2\right] \union T_{G^{-1}}\of{1}^\ast} \intersect \p{\left[\Fraction 1 / 2, 1\right] \union T_{G^{-1}}\of{\Fraction 1 / 4}^\ast} && \\
&& &= \set{\Fraction 1 / 2} \union \p{T_{G^{-1}}\of{1}^\ast \union T_{G^{-1}}\of{\Fraction 1 / 4}^\ast}. &&
\end{align*}
By Proposition~\ref{prop:3trans-characterisation}, it follows $\transitiveIn{G}{3} \subseteq \set{\Fraction 1 / 2} \union T_{G^{-1}}\of{1}^\ast \union T_{G^{-1}}\of{\Fraction 1 / 4}^\ast$. It is obvious $\Fraction 1 / 2 \notin \transitiveIn{G}{3}$, and straightforward to check $\transitiveIn{G}{3} \intersect \p{T_{G^{-1}}\of{1}^\ast \union T_{G^{-1}}\of{\Fraction 1 / 4}^\ast}$ is empty. 
\end{example}

%To derive a contradiction, suppose there exists $x \in \transitiveIn{G}{3}$. It is clear $x \notin \text{tr}\of{f_1} \union \text{tr}\of{f_2}$. For if $x \in \text{tr}\of{f_1}$, then $\Fraction 1 / 4 \notin T_G\of{x}^\ast$, implying $T_G\of{x}^\ast \subseteq [0, \Fraction 1 / 2]$ (which is not dense); the other case follows similarly. This implies $x \notin \set{\Fraction 1 / 4, 1}$, as it is easy to see this would imply either $x_1$ or $x_2$ are $3$-transitive in $\p{X, G}$. 

%By Proposition~\ref{prop:backwards}, $T_{G^{-1}}\of{x}^\ast \subseteq \transitiveIn{G}{3}$, which means $T_{G^{-1}}\of{x}^\ast$ and $\text{tr}\of{f_1} \union \text{tr}\of{f_2}$ are disjoint. We may assume $G\of{x}$ is not a singleton, which means $x \in \set{\Fraction 1/ 4, \Fraction 1 / 2, 1}$. It is clear $\Fraction 1 / 2 \notin \transitiveIn{G}{3}$, as $T_G\of{\Fraction 1 / 2}^\ast = \set{0, \Fraction 1 / 2}$. 

\bibliographystyle{plain}
\bibliography{references}

%---------------------------------------------------------------
\end{document}